\documentstyle{article}

\newcommand{\N}{{\rm I\kern-.5ex N}}
\newcommand{\Z}{{\sf \vrule height 1.55ex depth-1.2ex width.03em\kern-.11em Z
\kern-.9ex Z\kern-.11em\vrule height 0.3ex depth0ex width.03em}}
\newcommand{\Q}{{\rm\kern.2ex\vrule height1.55ex depth-.05ex width.03em\kern-.7ex Q}}
\newcommand{\R}{{\rm I\kern-.5ex R}}
\newcommand{\Rvar}{{\rm I\kern-.5ex R}}
\newcommand{\C}{{\rm\kern.3ex\vrule height1.55ex depth-.05ex width.03em\kern-.7ex C}}
\newcommand{\Cvar}{{\, \rm\kern.3ex\vrule height1.1ex depth-.05ex width.03em\kern-.7ex C}}

\newcommand{\hoed}{\hat{\rule{0ex}{1.7ex}}\,}
\newcommand{\golf}{\tilde{\rule{0ex}{1.7ex}}\,}

\newcommand{\spat}{\hspace{4ex}}

\newcommand{\flip}{ \chi }
\newcommand{\op}{B(H)}
\newcommand{\cop}{B_0(H)}

\newcommand{\ahh}{\hat{A} \hspace{-.55ex}\hat{\rule{0ex}{2.0ex}}\hspace{.55ex}}
\newcommand{\dhh}{\hat{\de} \hspace{-.95ex}\hat{\rule{0ex}{2.05ex}}\hspace{.95ex}}

\newcommand{\ar}{A_r}
\newcommand{\arh}{\hat{A}_r}
\newcommand{\ah}{\hat{A}}
\newcommand{\deh}{\hat{\Delta}}
\newcommand{\psih}{\hat{\psi}}
\newcommand{\lah}{\hat{\la}}
\newcommand{\pih}{\hat{\pi}}

\newcommand{\nab}{\nabla}

\newcommand{\cL}{{\cal L}}
\newcommand{\cF}{{\cal F}}
\newcommand{\cG}{{\cal G}}

\newcommand{\cN}{{\cal N}}
\newcommand{\cM}{{\cal M}}
\newcommand{\cP}{{\cal P}}
\newcommand{\cQ}{{\cal Q}}

\newcommand{\uU}{U}
\newcommand{\aU}{X}

\newcommand{\od}{\odot}
\newcommand{\ot}{\otimes}

\newcommand{\la}{\Lambda}

\newcommand{\om}{\omega}
\newcommand{\io}{\iota}
\newcommand{\vfi}{\varphi}
\newcommand{\vep}{\varepsilon}

\newcommand{\al}{\alpha}
\newcommand{\be}{\beta}

\newcommand{\sde}{\delta}
\newcommand{\de}{\Delta}
\newcommand{\et}{\eta}
\newcommand{\th}{\theta}

\newcommand{\si}{\sigma}
\newcommand{\Mfi}{{\cal M}_{\vfi}}
\newcommand{\Nfi}{{\cal N}_{\vfi}}

\newcommand{\lafi}{\la_\vfi}

\newcommand{\pifi}{\pi_\vfi}

\newcommand{\Mfiu}{{\cal M}_{\vfi_u}}
\newcommand{\Nfiu}{{\cal N}_{\vfi_u}}
\newcommand{\Mpsiu}{{\cal M}_{\psi_u}}

\newcommand{\cU}{{\cal U}}
\newcommand{\cV}{{\cal V}}

\newcommand{\cK}{{\cal K}}

\newcommand{\Nfir}{{\cal N}_{\vfi_r}}
\newcommand{\Mfir}{{\cal M}_{\vfi_r}}

\newcommand{\Mpsir}{{\cal M}_{\psi_r}}
\newcommand{\lar}{\la_r}

\newcommand{\text}[1]{\mbox{#1}}

\newcommand{\cst}{\text{C}$\hspace{0.1mm}^*$}

\newcommand{\qed}{\ \hfill \rule{2mm}{2mm}}
\newenvironment{demo}{\medskip\noindent\bf Proof :\ \  \rm}{\qed\bigskip\par }

\newtheorem{definition}{Definition}[section]
\newtheorem{proposition}[definition]{Proposition}
\newtheorem{lemma}[definition]{Lemma}
\newtheorem{corollary}[definition]{Corollary}
\newtheorem{remark}[definition]{Remark}
\newtheorem{theorem}[definition]{Theorem}

\newtheorem{notation}[definition]{Notation}
\newtheorem{result}[definition]{Result}

\begin{document}
\begin{center}
\huge\bf
Universal $\protect\boldmath C^*$-algebraic quantum groups arising from algebraic quantum
groups.

\end{center}

\bigskip

\begin{center}
\rm J. Kustermans  \footnote{Research Assistant of the
National Fund for Scientific Research (Belgium)}

Institut for Matematik og Datalogi

Odense Universitet

Campusvej 55

5230 Odense M

Denmark

\bigskip

\bf April 1997 \rm
\end{center}
\bigskip

\subsection*{Abstract}
In this paper, we construct a universal \cst-algebraic quantum group out of an algebraic
one. We show that this universal \cst-algebraic quantum has the same rich structure as its
reduced companion (see \cite{Kus}). This universal \cst-algebraic quantum group also
satisfies an upcoming definition of Masuda, Nakagami \& Woronowicz except for the possible
non-faithfulness of the left Haar weight.

\section*{Introduction.}
In 1979, Woronowicz proposed the use of the \cst-algebra language in the field of
quantum groups \cite{Wor4}. Since then, a lot of work has been done in this area but there is still
no satisfactory definition of a quantum group in the \cst-algebra framework.

However, the following subjects are better understood:
\begin{itemize}
\item compact \& discrete quantum groups and their duality theory
\item some examples: quantum $SU(2)$, quantum Heisenberg, quantum $E(2)$,
quantum Lorentz, duals of locally compact groups,...
\end{itemize}

A good definition of a quantum group should satisfy the following requirements:
\begin{itemize}
\item It should incorporate all the known examples and well understood parts of the theory.
\item There will have to be a good balance between the theory which can be extracted from
      the definition and the non complexity of the definition.
\item The definition has to allow a consistent duality theory.
\end{itemize}
The most difficult part of finding a satisfactory axiom scheme for quantum groups seems
to be the ability to prove the existence and uniqueness of a left Haar weight
from the proposed definition.

At the moment, Masuda, Nakagami and Woronowicz are working on a quasi-final definition of
a quantum group in the \cst-algebra framework (see also \cite{MasNak}).
It is quasi-final in the sense that there is some hope to find simpler axioms which
imply their current definition. Later, we will give a preview of this definition.

\medskip

In \cite{VD6}, A. Van Daele introduced the notion of Multiplier Hopf-algebras. They are natural
generalizations of Hopf algebras to the case of non-unital algebras.
Recently, A. Van Daele has looked into the case where such a Multiplier Hopf algebra
possesses a non-zero left invariant functional (Haar functional) and found some very
interesting properties \cite{VD1}. It should be noted
that everything in this theory is of an algebraic nature.

This category of Multiplier Hopf algebras with a Haar functional behaves very well in
different ways:
\begin{itemize}
\item This category includes the compact \& discrete quantum groups.
\item It is possible to construct the dual within this category.
\item The category is closed under the double construction of Drinfeld (see \cite{Drab}).
\end{itemize}
The last and the first property imply that this category contains the Lorentz quantum group
\cite{PW}.
However, this new category will not exhaust all quantum groups. It is not so difficult
to find many classical groups which are not Multiplier Hopf algebras. Also, quantum $E(2)$
will not fit in this scheme. Nevertheless, we still have a nice class of algebraic
quantum groups.

\medskip

In \cite{Kus}, we constructed a reduced \cst-algebraic quantum group (in the sense of
Masuda, Nakagmi \& Woronowicz) out of a Multiplier Hopf$\,^*$-algebra which possesses a
positive left invariant functional.

The main purpose of this paper is the construction of a universal \cst-algebraic quantum
group (in the sense of Masuda, Nakagmi \& Woronowicz) out of such a Multiplier
Hopf$\,^*$-algebra. The resulting universal \cst-algebraic quantum group fits into the
definition of Masuda, Nakagali \& Woronowicz except for the possible non-faithfulness of
the left Haar weight. We will see that this universal \cst-algebraic quantum group has the
same rich structure as its reduced companion. In order to prove this, we use the results
about the algebraic quantum group as well as the results about the reduced \cst-algebraic
quantum group.

\medskip

It will probably also be possible to construct some sort of universal quantum group out of
an \cst-algebraic quantum group according to Masuda, Nakagami \& Woronowicz. We believe
that this paper will give some ideas how to go around this in the general case. However,
some algebraic arguments will have to be modified to analytical ones.

\medskip

In a first section, we wil give an overview of the results of A. Van Daele about such a
Multiplier Hopf$\,^*$-algebra $(A,\de)$. In a second section, we recapitulate some results
concerning the reduced \cst-algebraic quantum group $(A_r,\de_r)$ arising from it.
However, most of the results concerning this reduced \cst-algebraic quantum group will be
introduced when we need them.

\medskip

In the third section, we introduce the universal \cst-algebra $A_u$ together with the
universal comultiplication $\de_u$. From there on, we gradually prove that this pair
$(A_u,\de_u)$  almost in the scheme of Masuda, Nakagami \& Woronowicz (except for the
faithfulness of the left Haar weight and some minor detail).

\medskip

In section \ref{art5}, we prove that every unitary corepresentation of the reduced quantum
group $(A_r,\de_r)$ is of an algebraic nature. The same is also true for the unitary
corepresentations of the universal quantum group $(A_u,\de_u)$.

The results of section \ref{art6} imply that there is a natural bijection between the
unitary corepresentations of $(A_r,\de_r)$ and of $(A_u,\de_u)$. The same is true for
bi-\cst-isomorphisms on the reduced and universal quantum group (section \ref{art5} and
\ref{art10}).

\medskip

In \cite{VD1},A. Van Daele constructs the dual algebraic quantum group $(\ah,\deh)$ of
$(A,\de)$. We can in the same way get a universal \cst-algebraic quantum group
$(\ah_u,\deh_u)$ out of $(\ah,\deh)$.

In section \ref{art12}, we introduce the universal corepresentation of the universal
quantum group $(A_u,\de_u)$ and use it to get a bijective correspondence between unitary
corepresentations of $(A_u,\de_u)$ and non-degenerate $^*$-homomorphisms on the universal
\cst-algebra $\ah_u$.

\medskip

Section \ref{app} fixes some notations about weights and gives a thorough overview of
slice weights and is necessary to understand some of the results of this paper.

\medskip

The notations of this paper will be a little bit different from the ones in \cite{Kus}.
Mainly because we have two \cst-algebraic quantum groups. As a rule, objects associated to
the reduced \cst-algebraic quantum group will get the subscript \lq r\rq, whereas objects
associated to the universal \cst-algebraic quantum group will get the subscript \lq u\rq.
Objects without a subscript will be connected to the algebraic quantum group.

\bigskip

First, we introduce some notations and conventions. We will always use the minimal tensor
product between \cst-algebras and use the symbol $\ot$ for this complete tensor product.
For any \cst-algebra $A$, we denote the Multiplier \cst-algebra by $M(A)$. The flip map
between two \cst-algebras will be denoted by $\flip$.

For the algebraic tensor products of vector spaces and linear mappings, we use the symbol
$\od$. The algebraic dual of a vector space $V$ will be denoted by $V'$.

Let $H$ be a Hilbert space. Then $B(H)$ will denote the \cst-algebra of bounded
operators on $H$, whereas $B_0(H)$ will denote the \cst-algebra of compact operators on $H$.
Consider vectors $v,w \in H$, then $\om_{v,w}$ is the element in $B_0(H)^*$
such that $\om_{v,w}(x)=\langle x v, w \rangle$ for all $x \in B_0(H)$.

The domain of an unbounded operator $T$ on $H$ is denoted by $D(T)$.

The domain of an element $\al$ which is affiliated with some \cst\ , will be denoted
by ${\cal D}(\al)$ (we will use the same notation for closed mappings in a \cst\ which arise
from one-parameter groups).

Whenever we say that an unbounded operator is positive, it is included that this operator
is also selfadjoint. The same rules apply to elements affiliated with a \cst-algebra.
Let $\al$ be an element affiliated with a \cst\ $A$, then $\al$ is called strictly positive
if $\al$ is positive and has dense range (it is then automatically injective).

A one-parameter group $\si$ on a \cst\ $A$ is called norm-continuous if and only if
for every $a \in A$, the mapping $\R \rightarrow A:t \mapsto \si_t(a)$ is norm-continuous.

\medskip

For a very readable introduction to Hilbert-\cst-modules, we refer to \cite{Lan}. Let $E$
and $F$ be two Hilbert-\cst-modules over a \cst-algebra $A$. Then the set of adjointable
mappings from $E$ into $F$ will be denoted by $\cL(E,F)$ whereas the set of compact
operators from $E$ into $F$ will be denoted by $\cK(E,F)$.

For the notion of regular operators between to Hilbert-\cst-modules, we refer to
\cite{Lan}.

\bigskip

As promised, we give now a preview of the definition of a \cst-algebraic quantum group
according to Masuda, Nakagami \& Woronowicz:

\medskip

Let $B$ be a \cst-algebra and $\de$ a non-degenerate $^*$-homomorphism from $B$
into $M(B\ot B)$ such that
\begin{enumerate}
\item $\de$ is coassociative, i.\ e.\ $(\de \ot \io)\de = (\io \ot \de)\de$.
\item $\de$ satisfies the following density conditions:
      $\de(B)(B \ot 1)$ and $\de(B)(1 \ot B)$ are dense  subsets of $B \ot B$.
\end{enumerate}
Furthermore, we assume the existence of the following objects:
\begin{enumerate}
\item a KMS-weight $\vfi$ on $B$ with modular group $\si$,
\item a norm continuous one parameter group $\tau$ on $B$,
\item an involutive $^*$-anti-automorphism $R$ on $B$,
\end{enumerate}
which satisfy the following properties:
\begin{enumerate}
\item For every $a \in \Mfi$, we have that $\de(a)$ belongs to ${\overline{\cal M}}_{\io \ot       \vfi}$ and $(\io \ot \vfi)\de(a) = \vfi(a) 1$.
\item Consider $a,b \in \Nfi$. Let $\om \in B^*$ such that
      $\om R \tau_{-\frac{i}{2}}$ is bounded and call $\th$ the unique element in $B^*$ which
      extends $\om R \tau_{-\frac{i}{2}}$. Then
      $$\vfi(\,b^* \, (\om \ot \io)\de(a)\,) = \vfi(\,(\th  \ot \io)(\de(b^*)) \,a \,) \ .$$
\item \begin{itemize}
      \item $\vfi$ is invariant under $\tau$.
      \item $\vfi$ commutes with $\vfi R$.
      \end{itemize}
\item \begin{itemize}
      \item For every $t \in \R$, we have that $R \tau_t = \tau_t R$.
      \item We have that $\de \tau_t = (\tau_t \ot \tau_t)\de$ for all $t \in \R$.
      \item $\de R = \flip (R \ot R) \de$
      \end{itemize}
\end{enumerate}
Then, we call $(B,\de,\vfi,\tau,R)$ a \cst-algebraic quantum group.

\medskip

We call $\vfi$ the left Haar weight, $R$ the anti-unitary antipode and $\tau$ the scaling
group of our quantum group. We put $\kappa = R \tau_{-\frac{i}{2}}$, then $\kappa$ plays
the role of the antipode of our quantum group.

This definition seems to be a \cst-version of a definition of a quantum group in the von
Neumann algebra setting (see \cite{MasNak}), which in turn was a generalization of the
(too restrictive) definition of a Kac-algebra (see \cite{E}). We are not sure that this
will be the ultimate definition of a \cst-algebraic quantum group proposed by Masuda,
Nakagami \& Woronowicz, but we expect that this one gives a fairly good idea of it.

A possible drawback of this definition is the complexity of the axioms. However, we will
show that the \cst-algebraic versions of Van Daele's algebraic objects fit almost in this
scheme. The only difference lies in the fact that we only can prove that $\vfi$ is
relatively invariant with respect to $\tau$ in stead of invariant. It is not clear at the
moment whether this definition of Masuda, Nakagami \& Woronowicz should be modified in
this respect.

\section{Algebraic quantum groups} \label{alg}

In this first section, we will introduce the notion of an algebraic
quantum group as can be found in \cite{VD1}. Moreover, we will give an
overview of the properties of this algebraic quantum group. The proofs
of these results can be found in the same paper \cite{VD1}. After this
section, we will construct a universal \cst-algebraic quantum group
out of this algebraic one, thereby heavily depending on the material
gathered in this section and \cite{Kus}. We will first introduce some
terminology.

\medskip

We call a $^*$-algebra $A$ non-degenerate if and only if we have for
every $a\in A$ that :
$$(\forall b \in A: a b =0) \Rightarrow a=0 \hspace{1.5 cm} \text{ and  } \hspace{1.5cm}
(\forall b \in A: b a =0) \Rightarrow a=0 .$$

For a non-degenerate $^*$-algebra $A$, you can define the multiplier
algebra $M(A)$. This is a unital $^*$-algebra in which $A$ sits as a
selfadjoint ideal (the definition of this multiplier algebra is the
same as in the case of \cst-algebras).

\medskip

If you have two non-degenerate $^*$-algebras $A,B$ and a
multiplicative linear mapping $\pi$ from $A$ to $M(B)$, we call $\pi$
non-degenerate if and only if the vector spaces $\pi(A) B$ and $B
\pi(A)$ are equal to $B$. Such a non-degenerate multiplicative linear
map has a unique multiplicative linear extension to $M(A)$. This
extension will be denoted by the same symbol as the original mapping.
Of course, we have similar definitions and results for
antimultiplicative mappings. If we work in an algebraic setting, we
will always use this form of non-degeneracy as opposed to the non
degeneracy of $^*$-homomorphisms between \cst-algebras!

\medskip

For a linear functional $\om$ on a non-degenerate $^*$-algebra $A$ and any $a \in M(A)$ we
define the linear functionals $\om a$ and $a \om$ on $A$ such that
$(a \om)(x) = \om(x a)$ and $(\om a)(x) = \om(a x)$ for every $x \in  A$.

\medskip

You can find some more information about non-degenerate algebras in the appendix of
\cite{VD6}.

\bigskip

Let $\om$ be a linear functional on a $^*$-algebra $A$, then :
\begin{enumerate}
\item $\om$ is called positive if and only if
$\om(a^* a)$ is positive for every $a \in A$.
\item If $\om$ is positive, then $\om$ is called faithful if and only if
for every $a \in A$, we have that $$\om(a^* a)=0 \Rightarrow a=0.$$
\end{enumerate}

\medskip

We have now gathered the necessary information to understand the following definition

\begin{definition}
Consider a non-degenerate $^*$-algebra $A$ and a non-degenerate $^*$-homomorphism $\de$
from $A$ into $M(A \od A)$ such that
\begin{enumerate}
\item $(\de \od \io)\de = (\io \od \de)\de$.
\item The linear mappings $T_1$, $T_2$ from $A \od A$ into $M(A \od A)$
such that $$T_1(a \ot b) = \de(a)(b \ot 1)
\hspace{1cm} \text{ and  } \hspace{1cm} T_2(a \ot b) = \de(a)(1 \ot b)$$
for all $a,b \in A$, are bijections from $A \od A$ to $A \od A$.
\end{enumerate}
Then we call $(A,\de)$ a Multiplier Hopf$\,^*$-algebra.
\end{definition}

In \cite{VD6}, A. Van Daele proves the existence of a unique non-zero
$^*$-homomorphism $\vep$ from $A$ to $\C$ such that $$(\vep \od
\io)\de = (\io \od \vep)\de =\io \ .$$ Furthermore, he proves the existence of a
unique anti-automorphism $S$ on $A$ such that
$$m(S \od \io)(\de(a)(1 \ot b)) = \vep(a) b \hspace{1cm} \text{ and } \hspace{1cm}
m(\io \od S)((b \ot 1)\de(a)) = \vep(a) b $$
for every $a,b \in A$  (here, $m$ denotes the multiplication map from $A \od A$ to $A$).
As usual, $\vep$ is called the counit and $S$ the antipode of $(A,\de)$.
Moreover, $S(S(a^*)^*) = a$ for all $a \in A$. Also, $\flip(S \od S)\de = \de S$.

\medskip

Let $\om$ be a linear functional on $A$ an $a$ an element in $A$. We define the element
$(\om \od \io)\de(a)$ in $M(A)$ such that
\begin{itemize}
\item $(\om \od \io)(\de(a)) \, \, b = (\om \od \io)(\de(a)(1 \ot b))$ \item
$b \, \, (\om \od \io)(\de(a)) = (\om \od \io)((1 \ot b)\de(a))$
\end{itemize}
for every $b \in A$.

In a similar way, the multiplier $(\io \od \om)\de(a)$ is defined.

\bigskip

Let $\om$ be a linear functional on $A$. We call $\om$ left invariant
(with respect to $(A,\de)$), if and only if $(\io \od \om)\de(a) =
\om(a)\,1$ for every $a \in A$. Right invariance is defined in a
similar way.

\begin{definition}
Consider a Multiplier Hopf$\,^*$-algebra $(A,\de)$ such that there exists a
non-zero positive linear functional $\vfi$ on $A$ which is left invariant.  Then we call
$(A,\de)$ an algebraic quantum group.
\end{definition}

For the rest of this paper, we will fix an algebraic quantum group $(A,\de)$ together with
a non-zero left invariant positive linear functional $\vfi$ on it.

An important feature of such an algebraic quantum group is the
faithfulness and uniqueness of left invariant functionals :
\begin{enumerate}
\item Consider a left invariant linear functional $\om$ on $A$, then there exists a unique
      element $c \in \C$ such that $\om = c \, \vfi$.
\item Consider a non-zero left invariant linear functional $\om$ on $A$, then $\om$ is
      faithful.
\end{enumerate}
In particular, $\vfi$ is faithful.

\medskip

A first application of this uniqueness result concerns the antipode :
Because $\vfi S^2$ is left invariant, there exists a unique complex
number $\mu$ such that $\vfi S^2 = \mu \vfi$ (in \cite{VD1}, our $\mu$
is denoted by $\tau$!). It is not so difficult to prove in an
algebraic way that $|\mu|=1$. The question remains open if there
exists an example of an algebraic quantum group (in our sense) with
$\mu \neq 1$.

\medskip

It is clear that $\vfi S$ is a non-zero right invariant linear
functional on $A$. However, in general, $\vfi S$ will not be positive.
In \cite{Kus}, we use  the \cst-algebra approach to prove the
existence of a non-zero positive right invariant linear functional on
$A$.

Of course, we have similar faithfulness and uniqueness results about right invariant linear
functionals.

\medskip

In this paper, we will need frequently the following formula :
\begin{equation}
(\io \od \vfi)(\,(1 \ot a)\de(b)\,)
= S(\,(\io \od \vfi)(\de(a)(1 \ot b))\,)   \label{eq1.3}  \end{equation}
for all $a,b \in A$. A proof of this result can be found in proposition
3.11 of \cite{VD1}. It is in fact nothing else but an algebraic form of the strong left
invariance in the definition of Masuda, Nakagami \& Woronowicz.

\medskip

Another non-trivial property about $\vfi$ is the existence of a unique automorphism $\rho$
on $A$ such that $\vfi(a b) = \vfi(b \rho(a))$ for every $a,b \in A$. We call this the
weak KMS-property of $\vfi$ (In \cite{VD1}, our mapping $\rho$ is denoted by $\sigma$!).

This weak KMS-property is crucial to extend $\vfi$ to a weight on the \cst-algebra level.
We have moreover that $\rho(\rho(a^*)^*) = a$ for every $a \in A$.

As usual, there exists a similar object $\rho'$ for the right
invariant functional $\vfi S$, i.e. $\rho'$ is an automorphism on $A$
such that $(\vfi S)(a b) = (\vfi S)(b \rho'(a))$ for every $a,b \in
A$.

Using the antipode, we can connect $\rho$ and $\rho'$ via the formula
$S \rho' = \rho S$. Furthermore, we have that $S^2$ commutes with
$\rho$ and $\rho'$. The interplay between $\rho$,$\rho'$ and $\de$ is
given by the following formulas :
$$\de \rho = (S^2 \od \rho)\de \hspace{1cm} \text{ and } \hspace{1cm} \de \rho' =
(\rho' \od S^{-2})\de.$$

\medskip

It is also possible to introduce the modular function of our algebraic quantum group. This
is an invertible element $\sde$ in $M(A)$ such that $$(\vfi \od \io)(\de(a)(1 \ot b)) =
\vfi(a) \, \sde  \, b $$ for every $a,b \in A$.

Concerning the right invariant functional, we have that $$(\io \od \vfi S)(\de(a)(b \ot
1)) = (\vfi S)(a) \, \sde^{-1} \, b$$ for every $a,b \in A$.

This modular function is, like in the classical group case, a one
dimensional (generally unbounded) corepresentation of our algebraic
quantum group :
$$ \de(\sde) = \sde \od \sde \hspace{2cm} \vep(\sde)= 1 \hspace{2cm}
S(\sde) = \sde^{-1} .$$

As in the classical case, we can relate the left invariant functional
to our right invariant functional via the modular function : we have
for every $a \in A$ that
$$ \vfi(S(a)) =  \vfi(a \sde) = \mu \, \vfi(\sde a) .$$
If we apply this equality two times and use the fact that $S(\sde) =
\sde^{-1}$, we get that $ \vfi(S^2(a)) = \vfi(\sde^{-1} a \sde)$ for every $a \in A$.

\medskip

Not surprisingly, we have also that $\rho(\sde) = \rho'(\sde) =
\mu^{-1} \sde$.

\medskip

Another connection between $\rho$ and $\rho'$ is given by the equality
$\rho'(a) = \sde \rho(a) \sde^{-1}$ for all $a \in A$.

\bigskip

We have also a property which says, loosely speaking, that every
element of $A$ has compact support :

Consider $a_1,\ldots\!,a_n \in A$. Then there exists an element $c$ in $A$ such that
$c \, a_i = a_i \, c = a_i$ for every $i \in \{1,\ldots\!,n\}$.

\bigskip\bigskip

In a last part, we are going to say something about  duality.

We define the subspace $\ah$ of $A'$ as follows :
$$\ah = \{\, \vfi a \mid a \in A \,\} = \{\, a \vfi \mid a \in A \,\}.$$
Like in the theory of Hopf$\,^*$-algebras, we turn $\ah$ into a
non-degenerate $^*$-algebra :
\begin{enumerate}
\item For every $\om_1,\om_2 \in \ah$ and $a \in A$, we have that $(\om_1 \om_2)(a) =
(\om_1 \od \om_2)(\de(a))$.
\item For every $\om \in \ah$ and $a \in A$, we have that $\om^*(a) = \overline{\om(S(a)^*)}$.
\end{enumerate}
We should remark that a little bit of care has to be taken by defining the product and
the $^*$-operation in this way.

Also, a comultiplication $\deh$ can be defined on $\ah$
such that $\deh(\om)(x \ot y) = \om(x y)$ for every $\om \in \ah$ and $x,y \in A$.

Again, this has to be made more precise. This can be done by embedding
$M(A \od A)$ into $(A \od A)'$ in the right way but we will not go
into this subject (see \cite{JK3} for more information about this). A
definition of the comultiplication $\deh$ without the use of such an
embedding can be found in definition $4.4$ of \cite{VD1}.

In this way, $(\ah,\deh)$ becomes a Multiplier Hopf$\,^*$-algebra. The counit $\hat{\vep}$
and the antipode $\hat{S}$ are such that
\begin{enumerate}
\item For every $\om \in \ah$, we have that $\hat{\vep}(\om) = \om(1)$.
\item For every $\om \in \ah$ and every $a \in A$, we have that $\hat{S}(\om)(a) = \om(S(a))$.
\end{enumerate}

For any $a \in A$, we define $\hat{a} = a \vfi \in \ah$. The mapping $A \rightarrow \ah :
a \mapsto \hat{a}$ is a bijection, which is in fact nothing else but the Fourier transform.

\medskip

Next, we define the linear functional $\psih$ on $\ah$ such that
$\psih(\hat{a}) = \vep(a)$ for every $a \in A$. It is possible to prove that $\psih$
is right invariant.

Also, we have that $\psih(\hat{a}^* \hat{a}) = \vfi(a^* a)$ for every $a \in A$. This implies
that $\psih$ is a non-zero positive right invariant linear functional on $\ah$.

\medskip

From theorem 9.9 of \cite{Kus}, we know  that $(A,\de)$ possesses a non-zero positive right
invariant linear functional.  In a similar way, this functional will
give rise to a non-zero positive left invariant linear functional on
$\ah$. This will imply that $(\ah,\deh)$ is again an algebraic quantum
group.

\bigskip\bigskip

In definition 6.5 of \cite{JK3}, we introduced the universal corepresentation $\aU$ of $A$
(it was denoted by $U$ in \cite{JK3}). This is an element  of $M(A \od \ah)$ such that
$(\de \od \io)(\aU) = \aU_{13}\,\aU_{23}$ and $(\io \od \deh)(\aU) = \aU_{12}\,\aU_{13}$.

This element $\aU$ serves as a bridge between unitary corepresentations
of $(A,\de)$ and non-degenerate \newline $^*$-homomorphisms on $\ah$.

\bigskip

Consider a unitary corepresentation $\cU$ of $(A,\de)$ on a $^*$-algebra $C$. Let $\om
\in A'$ and $a \in A$.

By proposition 4.3 and result 5.9, we have an element $(a \om \od
\io)(\cU)$ in $M(C)$ such that the following holds :
\begin{enumerate}
\item We have for every $c \in C$ that $\cU (a \ot c)$ belongs to $A \od C$ and
$(a \om \od \io)(\cU) \, c = (\om \od \io)(\cU (a \ot c))$.
\item We have for every $c \in C$ that $(1 \ot c) \, \cU (a \ot 1)$ belongs to $A \od C$ and
\newline $c \, (a \om \od \io)(\cU)  = (\om \od \io)((1 \ot c)\,\cU (a \ot
1))$.
\end{enumerate}
Notice that in this way, we defined $(\rho \od \io)(\cU)$ for every
$\rho \in \ah$.

\medskip

\begin{proposition}
Consider a non-degenerate $^*$-algebra $C$. Let $\cU$ be a unitary
corepresentation of $(A,\de)$ on $C$. Define the mapping $\th$ from
$\ah$ into $M(C)$ such that $\th(\om) = (\om \od \io)(\cU)$ for every
$\om \in \ah$. Then $\th$ is a non-degenerate $^*$-homomorphism from
$\ah$ into $M(C)$ such that $\th(a \om) = (a \om \od \io)(\cU)$ for
every $a \in A$ and $\om \in A'$.
\end{proposition}

\bigskip

\begin{theorem}
Consider a non-degenerate $^*$-algebra $C$. Let $\cU$ be a unitary
corepresentation of $(A,\de)$ on $C$.  Define the mapping $\th$ from
$\ah$ into $M(C)$ such that $\th(\om) = (\om \od \io)(\cU)$ for every
$\om \in \ah$. Then $\th$ is the  unique non-degenerate
$^*$-homomorphism from $\ah$ into $M(C)$ such that $(\th \od \io)(\aU)
= \cU$.
\end{theorem}

We have also a converse of this :

\begin{theorem}
Consider a non-degenerate $^*$-algebra $C$. Let $\th$ be a
non-degenerate $^*$-homomorphism of $\ah$ into $M(C)$. Define $\cU =
(\io \od \th)(\aU)$. Then $\cU$ is the unique unitary corepresentation
of $(A,\de)$ on $C$ satisfying $\th(\rho) = (\rho \od \io)(\cU)$  for
every $\rho \in \ah$.
\end{theorem}

\section{The reduced \cst-algebraic quantum group arising from $(A,\de)$}
\label{art1}

In \cite{Kus}, we constructed a reduced \cst-algebraic quantum group
$(A_r,\de_r)$ in the sense of Masuda, Nakagami \& Woronowicz out of
the algebraic quantum group $(A,\de)$. The aim of this paper to
construct the universal \cst-algebraic quantum group out of $(A,\de)$.
We will show that this resulting universal \cst-algebraic quantum
group has the same rich structure as its reduced companion.

The proofs of the results about the universal case depend  heavily on
the results concerning the reduced case. Therefore we first
recapitulate  some results about the reduced
\cst-algebraic quantum group. This section gives an overview of the results
of section 2 of \cite{Kus} .

\bigskip

For the rest of this paper, we fix a GNS-pair $(H,\la)$ of the left Haar functional $\vfi$
on $A$. This means that $H$ is a Hilbert space and $\la$ is a linear mapping from $A$ into
$H$ such that
\begin{enumerate}
\item The set $\la(A)$ is dense in $H$.
\item We have for every $a,b \in A$ that $\langle \la(a) , \la(b) \rangle
= \vfi(b^* a)$.
\end{enumerate}

\medskip

As usual, we can associate a multiplicative unitary to $(A,\de)$ :

\begin{definition}
We define $W$ as the unique unitary element in $B(H \ot H)$ such that
\newline $W \, (\la \od \la)(\de(b)(a \ot 1)) = \la(a) \ot \la(b)$ for every $a,b \in A$.
The element $W$ is called the fundamental unitary associated to $(A,\de)$.
\end{definition}

The coassociativity of $\de$ on the algebra level implies that $W$ is
multiplicative : $W_{12} W_{13} W_{23} = W_{23} W_{12}$.

\bigskip

The GNS-pair $(H,\la)$ allows us to represent $A$ by bounded operators
on $H$ :

\begin{definition}
We define $\pi_r$ as the unique $^*$-homomorphism from $A$ into $B(H)$ such that $\pi_r(a)
\la(b) = \la(a b)$ for every $a,b \in A$. We have also that $\pi_r$ is injective.
\end{definition}

We would like to mention that $\pi_r$ is denoted by $\pi$ in
\cite{Kus}, but we will reserve the symbol $\pi$ for another object in
this paper!

\medskip

Notice also that it is not immediate that $\pi_r(x)$ is a bounded
operator on $H$ (because $\vfi$ is merely a functional, not a weight),
but the boundedness of $\pi_r(x)$ is connected with the following
equality :

\medskip

We have for every $a,b \in A$ that
\begin{equation}
\pi_r\bigl((\io \ot \vfi)( \de(b^*)(1 \ot a))\bigr)
= (\io \ot \om_{\la(a),\la(b)})(W)      \label{eq1.1}
\end{equation}

The mapping $\pi_r$ makes it possible to define our reduced
\cst-algebra :

\begin{definition}
We define $\ar$ as the closure of $\pi_r(A)$ in $\op$. So $\ar$ is a non-degenerate
sub-\cst-algebra of $B(H)$.
\end{definition}

Equation \ref{eq1.1} implies that
\begin{equation}
\ar = \text{\ \ closure of \ \ }
\{\, (\io \ot \om)(W) \mid \om \in \cop^* \,\}  \text{\ \ in \  }  \op \ . \label{eq1.2}
\end{equation}

As usual, we use the fundamental unitary to define a comultiplication
on $A_r$. We will denote it by $\de_r$ (in \cite{Kus}, it is denoted
by $\de$).

\begin{definition}
We define the mapping $\de_r$ from $\ar$ into $B(H \ot H)$ such that
$\de_r(x) =  W^* (1 \ot x) W$ for all $x \in \ar$. Then $\de_r$ is an
injective $^*$-homomorphism.
\end{definition}

It is not so difficult to show that $\de_r$ on the \cst-algebra level
is an extension of $\de$ on the $^*$-algebra level :

\begin{result}
We have for all $a \in A$ and $x \in A \od A$ that $(\pi_r
\od \pi_r)(x) \, \de_r(\pi_r(a)) = (\pi_r \od \pi_r)(x \, \de(a))$
and $\de_r(\pi_r(a)) \, (\pi_r \od \pi_r)(x) = (\pi_r \od \pi_r)(\de(a)\,x).$
\end{result}

Using this result, it is easy to prove formulas like :
$$ \de_r(\pi_r(a)) \, (1 \ot \pi_r(b)) = (\pi_r \od \pi_r)(\de(a)(1 \ot b))$$
for all $a,b \in A$.

\medskip

Using the above results, it is not so hard to prove the following
theorem :

\begin{theorem}
We have that $A_r$ is a non-degenerate sub-\cst-algebra of $\op$ and
$\de_r$ is a non-degenerate injective $^*$-homomorphism from $A_r$
into $M(A_r \ot A_r)$ such that :
\begin{enumerate}
\item $(\de_r \ot \io)\de_r = (\io \ot \de_r)\de_r$
\item The vector spaces $\de_r(A_r)(A_r \ot 1)$ and
 $\de_r(A_r)(1 \ot A_r)$ are dense subspaces of $A_r \ot A_r$.
\end{enumerate}
\end{theorem}

\bigskip\bigskip

It is also possible to represent the dual multiplier Hopf$\,^*$-algebra $\ah$ on $H$.

Remember from section \ref{alg} that we have a non-zero positive right invariant linear
functional  $\psih$ on $\ah$. We have moreover that $\psih(\hat{b}^*\,\hat{a}) = \vfi(b^*
a)$ for all $a,b \in A$.

\medskip

We define the linear map $\lah$ from $\ah$ into $H$ such that
$\lah(\hat{a})= \la(a)$ for every $a \in A$. Then
\begin{enumerate}
\item $\lah(\ah)$ is dense in $H$
\item $\langle \lah(a) , \lah(b) \rangle = \psih(b^* a)$  for every $a,b \in \ah$.
\end{enumerate}

\bigskip

The multiplicative unitary $W$ can also be expressed in terms of
$\lah$ :

\begin{result}
We have for every $\om_1,\om_2 \in \ah$ that $W\,(\lah(\om_1) \ot \lah(\om_2)) = (\lah \od
\lah)(\deh(\om_1)(1 \ot \om_2))$
\end{result}

With this expression in hand, we can do the same things for the dual
$\ah$ as we did for $A$ itself.

\begin{definition}
We define the mapping $\pih_r$ from $\ah$ into $\op$ such that $\pih_r(\om) \lah(\th) =
\lah(\om \th)$ for all $\om,\th \in \ah$. Then $\pih_r$ is an injective $^*$-homomorphism.
\end{definition}

Furthermore,  we have for every $\th,\et \in \ah$ that
\begin{equation}
(\om_{\lah(\th),\lah(\et)} \ot \io)(W) = \pih_r\bigl((\psih \od \io)((\et^* \ot
1)\deh(\th))\bigr).    \label{eq1.4}
\end{equation}

\begin{definition}
We define $\arh$ as the closure of $\pih_r(\ah)$ in $\op$. So $\arh$ is a non-degenerate
sub-\cst-algebra of $B(H)$.
\end{definition}

Equation \ref{eq1.4} implies that
$$\arh = \text{\ \ closure of \ \ } \{\, (\om \ot \io)(W)
\mid \om \in \cop^* \,\} \text{\ \ in\ } \op,$$ which is again something familiar.

\begin{definition}
We define the mapping $\deh_r$ from $\arh$ into $B(H \ot H)$ such that
$\deh_r(x) = W (x \ot 1) W^*$ for all $x \in \arh$. Then $\deh_r$ is
an injective $^*$-homomorphism.
\end{definition}

We will of course also have the following theorem.

\begin{theorem}
We have that $\ah_r$ is a non-degenerate sub-\cst-algebra of $\op$ and
$\deh_r$ is a non-degenerate injective $^*$-homomorphism from $\ah_r$
into $M(\ah_r \ot \ah_r)$ such that :
\begin{enumerate}
\item $(\deh_r \ot \io)\deh_r = (\io \ot \deh_r)\deh_r$
\item The vector spaces $\deh_r(\ah_r)(\ah_r \ot 1)$ and $\deh_r(\ah_r)(1 \ot \ah_r)$
are dense subspaces of $\ah_r \ot \ah_r$.
\end{enumerate}
\end{theorem}

\bigskip\bigskip

Lemma 2.21 of \cite{Kus} and definition 6.5 of \cite{JK3} imply the
next proposition.

\begin{proposition} \label{prop1.3}
We have for every $x \in A_r \od \ah_r$ that
$W \, (\pi_r \od \hat{\pi}_r)(x) = (\pi_r \od \hat{\pi}_r)(\aU \, x)$
\newline and $(\pi_r \od \hat{\pi}_r)(x) \, W = (\pi_r \od \hat{\pi}_r)(x \, \aU)$.
\end{proposition}

\begin{corollary}
The element $W$ belongs to $M(A_r \ot \ah_r)$.
\end{corollary}

\hspace*{1.5cm}

An important object associated to $(A_r,\de_r)$ is the left Haar
weight on $A_r$. It is determined by the following theorem (see
\cite{Kus},  theorem 6.12, proposition 6.2 and the remarks after it) :

\begin{theorem}
There exists a unique closed linear map $\lar$ from  within $A_r$ into $H$ such that
$\pi_r(A)$ is a core for $\lar$ and $\lar(\pi_r(a)) = \la(a)$ for every $a \in A$.

There exists moreover a unique weight $\vfi_r$ on $A_r$ such that
$(H,\lar,\io)$ is a GNS-construction for $\vfi_r$.

We have also that $\pi(A) \subseteq \Mfir$ and that $\vfi_r(\pi(a)) =
\vfi(a)$ for every $a \in A$.
\end{theorem}

\begin{proposition}  \label{prop1.1}
The weight $\vfi_r$ is a faithful KMS-weight. We denote the modular group of $\vfi_r$ by
$\si_r$.
\end{proposition}

We denote the modular operator of $\vfi_r$ by $\nab$ and the modular
conjugation of $\vfi_r$ by $J$ (both with respect to the
GNS-representation $(H,\lar,\io)$\ ).

\medskip

Put $T = J \nab^{\frac{1}{2}}$. By the remarks after proposition 3.2
of \cite{Kus}, we know that $\la(A)$ is a core for $T$ and that $T
\la(a) = \la(a^*)$ for every $a \in A$.
We have moreover that  $\la(A) \subseteq D(T^*)$ and that $T^*
\la(a) = \la(\rho(a^*))$ for every $a \in A$.

\medskip

So $\la(A) \subseteq D(\nab)$ and $\nab \la(a) = \la(\rho(a))$ for
every $a \in A$.

This implies for every element $a \in A$ that $\la(a)$ is analytic
with respect to $\nab$ and that $\nab^n \la(a) = \la(\rho^n(a))$ for
every $n \in \Z$.

\bigskip

The left invariance on the $^*$-algebra level is then transferred to
the left invariance on the \cst-algebra level. For used notations, we
refer to the appendix.

\begin{theorem}
Consider $x \in \Mfir$, then $\de_r(x)$ belongs to
$\overline{{\cal M}}_{\io \ot \vfi_r}$ and
$(\io \ot \vfi_r)\de_r(x) = \vfi_r(x) \, 1$.
\end{theorem}

For some more detailed information about the left invariance, we refer
to section 6 of \cite{Kus} and \cite{JK4}.

We will need the following natural formula for $W$. For notations, we refer to section
\ref{app2}

\begin{proposition} \label{prop1.2}
Consider $a,b \in \Nfir$. Then $\de(b)(a \ot 1)$ belongs to
$\cN_{\vfi_r \ot \vfi_r}$ and \newline $W (\la_r \ot \la_r)(\de(b)(a
\ot 1)) = \la_r(a) \ot \la_r(b)$.
\end{proposition}

This proposition follows from the definition of $W$ and  the
definition of $\la_r$. We will give an explicit proof of a similar
result in a later section (see proposition \ref{prop9.1}).

\section{The universal bi-\cst-algebra}

In this section, we will introduce the universal enveloping
\cst-algebra of the $^*$-algebra $A$ and extend the comultiplication
to this \cst-algebra. In order to do so, we will need the following
crucial lemma.

\medskip

\begin{lemma} \label{lem2.1}
Consider a \cst-algebra  $C$ and a sub-$^*$-algebra $B$ of $M(C)$ such that $B C$ is dense
in $C$. Let $\phi$ be a $^*$-homomorphism from $\ah$ into $M(B)$ such that $\phi(\ah) B =
B$.

Take  $\om \in \ah$ and define $\th$ to be the unique element in $A_r^*$ such that $\th
\!  \circ \! \pi_r = \om$. Then there exists a unique element $T$ in $M(C)$ such that
$\phi(\om)\,b = T\,b$ and $b\,\phi(\om) = b\,T$ for every $b \in B$. Furthermore,
$\|T\| \leq \|\th\|$.
\end{lemma}

\begin{remark} \rm
We should be careful with the terminology above. The $^*$-algebra $B$ has only an
algebraic character and the mapping $\phi$ is an algebraically non-degenerate
$^*$-homomorphism from $\ah$ into $M(B)$. This implies that $\phi(\om)$ is an element of
the algebraic multiplier $^*$-algebra $M(B)$. In general, there is no reason for
$\phi(\om)$ to be an element of $M(C)$. Nevertheless, the previous lemma guarantees that
this is true in this special case.
\end{remark}

\medskip\noindent\bf Proof of the lemma :\ \  \rm
For the moment, consider $\pi_r$ as an algebraically non-degenerate $^*$-homomorphism from
$A$ into $\pi_r(A)$ and $\phi$ as an algebraically non-degenerate $^*$-homomorphism from
$\ah$ into $M(B)$ \ (The non-degeneracy of $\pi_r$ follows from the fact that $A^2=A$).

Then the remarks at the end  of section \ref{alg} imply the existence
of a unitary element $\cV$ in the algebraic multiplier $^*$-algebra
$M(A \od B)$ such that $\phi(\rho) = (\rho \od \io)(\cV)$ for every
$\rho \in \ah$. Put $\cP = (\pi_r \od \io)(\cV)$, then $\cP$ is a unitary
element in the algebraic multiplier $^*$-algebra $M(\pi_r(A) \od B)$.

\medskip

It is easy to see that
$$\bigl([\cP(a_2 \ot b_2)](1 \ot c_2)\bigr)^*
\, \bigl([\cP(a_1 \ot b_1)](1 \ot c_1)\bigr)
= (a_2 \ot b_2 c_2)^* \, (a_1  \ot b_1 c_1)   $$
for every $a_1,a_2 \in \pi_r(A)$, $b_1,b_2 \in B$, and $c_1,c_2 \in C$. Using the fact that
$\pi_r(A)$ is dense in $A_r$ and the fact that $B C$ is dense in $C$, this implies the
existence of a unique unitary element $\cQ \in M(A_r \ot C)$ such that
$$\cQ (a \ot b c) = [\cP(a \ot b)](1 \ot c)$$ for every $a \in \pi_r(A)$,  $b \in B$ and
$c \in C$.

This implies immediately that $\cQ(a \ot b) = \cP(a \ot b)$ for every $a \in
\pi_r(A)$ and $b \in B$.

Define $T = (\th \ot \io)(\cQ)$, then $T$ is an element in $M(C)$ and
$\|T\| \leq \|\th\|$.

\medskip

Because $\om$ belongs to $\ah$, there exists $a \in A$ such that $\om = a \om$. Because $\th
\! \circ \! \pi_r = \om$, this implies that $\th = \pi_r(a)  \th$. So we get for every
$b \in B$ that
\begin{eqnarray*}
T\,b & = & (\th \ot \io)(\cQ) \, b = (\pi_r(a) \th \ot \io)(\cQ) \, b
=  (\th \ot \io)(\cQ(\pi_r(a) \ot b)) \\
& = & (\th \ot \io)(\cP(\pi_r(a) \ot b)) = (\th \od \io)(\cP(\pi_r(a) \ot b)) \\
& = & (\th \od \io)( (\pi_r \od \io)(\cV)\,(\pi_r(a) \ot b))
= (\th \od \io)\bigl( (\pi_r \od \io)(\cV(a \ot b))\bigr) \\
& = & (\om \od \io)(\cV(a \ot b))
=  (a \om \od \io)(\cV) \,  b = (\om \od \io)(\cV) \,  b
= \phi(\om) \,  b \ .
\end{eqnarray*}
Using the associativity in $M(B)$ and $M(C)$, we get from this result also that
$b\,\phi(\om) = b \, T$ for every $b \in B$.
\qed\bigskip\par

\begin{result}
Consider a \cst-algebra $C$ and a $^*$-homomorphism $\phi$ from $\ah$ into $M(C)$ such
that $\phi(\ah) C$ is dense in $C$. Let $\om$ be an element in $\ah$ and define $\th$ to
be the unique element in $A_r^*$ such that $\th \! \circ \! \pi_r = \om$. Then
$\|\phi(\om)\|
\leq \|\th\|$.
\end{result}
\begin{demo}
Put $B = \phi(\ah)$. We have by assumption that $B C$ is dense in $C$. It is also clear
that $\phi$ is a $^*$-homomorphism from $\ah$ into $B$. Because $\ah^2 = \ah$, we have
that $\phi(\ah) B = B$.

Therefore, the previous lemma implies the existence of a unique element $T \in M(C)$ such
that $\|T\| \leq \|\th\|$ and $T \, b = \phi(\om) \, b$ for every $b \in B$. Because $B C$
is dense in $C$, this implies that $T = \phi(\om)$.
\end{demo}

This result implies the following one :

\begin{corollary}
Consider a \cst-algebra $C$ and a $^*$-homomorphism $\phi$ from $\ah$ into $C$. Let $\om$
be an element in $\ah$ and define $\th$ to be the unique element in $A_r^*$ such that $\th
\! \circ \! \pi_r = \om$. Then $\|\phi(\om)\| \leq \|\th\|$.
\end{corollary}

As a consequence, we get the following boundedness property.

\begin{corollary}
Consider $\om \in \ah$. Then there exists a positive number M such that $\|\phi(\om)\|
\leq M$ for every \cst-algebra $C$ and every $^*$-homomorphism $\phi$ from $\ah$ into
$C$.
\end{corollary}

\bigskip

Another application of lemma \ref{lem2.1} can be found in the
following result :

\begin{result}
Consider two \cst-algebras $C_1$, $C_2$ and $^*$-homomorphisms $\phi_1$ from $\ah$ into
$M(C_1)$ and $\phi_2$ from $\ah$ into $M(C_2)$ such that $\phi_1(\ah)C_1$ is dense in
$C_1$ and $\phi_2(\ah)C_2$ is dense in $C_2$.

Then there exists a unique $^*$-homomorphism $\phi$ from $\ah$ into $M(C_1 \ot C_2)$ such
that $\phi(\om)\,(\phi_1(\om_1) \ot \phi_2(\om_2)) = (\phi_1
\od \phi_2)(\deh(\om)(\om_1 \od \om_2))$ and $(\phi_1(\om_1) \ot \phi_2(\om_2))\,\phi(\om) =
(\phi_1 \od \phi_2)((\om_1 \od \om_2)\deh(\om))$ for every $\om,\om_1,\om_2 \in \ah$. We have
moreover that $\phi(\ah)(C_1 \ot C_2)$ is dense in $C_1 \ot C_2$.
\end{result}
\begin{demo}
Define $B_1 = \phi_1(\ah)$ and $B_2 = \phi_2(\ah)$. Then $B_1$ is a sub$^*$-algebra of
$M(C_1)$ such that $B_1 C_1$ is dense in $C_1$ and $B_2 C_2$ is dense in $C_2$. This
implies that $B_1 \od B_2$ is a sub$^*$-algebra of $M(C_1 \ot C_2)$ such that $(B_1 \od
B_2)(C_1 \ot C_2)$ is dense $C_1 \ot C_2$.

Consider $\phi_k$ as an algebraically non-degenerate $^*$-homomorphism from $\ah$ into
$B_k$ ($k=1,2$) and  define $\tilde{\phi} = (\phi_1 \od \phi_2)\de$ in the algebraic way.
So $\tilde{\phi}$ is an algebraically non-degenerate $^*$-homomorphism from $\ah$ into
$M(B_1 \od B_2)$.

By lemma \ref{lem2.1}, there exists a unique mapping $\phi$ from $\ah$ into $M(C_1 \ot
C_2)$ such that $\phi(\om) \, b = \tilde{\phi}(\om) \, b$ and $b \, \phi(\om) = b \,
\tilde{\phi}(\om)$ for every $b \in B_1 \od B_2$. Because $\tilde{\phi}$ is a
$^*$-homomorphism, we get easily that $\phi$ is a $^*$-homomorphism.

We have also for every $\om_1,\om_2,\om \in \ah$ that
\begin{eqnarray*}
& & \phi(\om)\,(\phi_1(\om_1) \ot \phi_2(\om_2)) = \tilde{\phi}(\om)
\,(\phi_1(\om_1) \ot \phi_2(\om_2)) \\
& & \spat = (\phi_1 \od \phi_2)(\deh(\om))\,(\phi_1(\om_1) \ot \phi_2(\om_2))
= (\phi_1 \od \phi_2)(\deh(\om)(\om_1 \od \om_2)) \ .
\end{eqnarray*}
The other equality is proven in a similar way.

From this, we get easily that $\phi(\ah)(B_1 \od B_2) = B_1 \od B_2$ which implies that
$\phi(\ah)(C_1 \ot C_2)$ is dense in $C_1 \ot C_2$.
\end{demo}

\bigskip\bigskip

By duality, we get the following results :

\begin{proposition}  \label{prop2.1}
Consider $a \in A$. Then there exists a positive number M such that $\|\phi(a)\| \leq M$
for every \cst-algebra $C$ and every $^*$-homomorphism $\phi$ from $A$ into $C$.
\end{proposition}

\begin{proposition} \label{prop2.2}
Consider \cst-algebras $C_1$, $C_2$ and $^*$-homomorphisms $\phi_1$ from $A$ into $M(C_1)$
and $\phi_2$ from $A$ into $M(C_2)$ such that $\phi_1(A)C_1$ is dense in $C_1$ and
$\phi_2(A)C_2$ is dense in $C_2$.

Then there exists a unique $^*$-homomorphism $\phi$ from $A$ into $M(C_1 \ot C_2)$ such that
$(\phi_1(a_1) \ot \phi_2(a_2)) \, \phi(a) = (\phi_1 \od \phi_2)((a_1 \ot a_2)\de(a))$ and
$\phi(a)\,(\phi_1(a_1) \ot \phi_2(a_2)) = (\phi_1 \od \phi_2)(\de(a)(a_1 \ot a_2))$ for every
$a \in A$ and $a_1,a_2 \in A$. We have moreover that $\phi(A)(C_1 \ot C_2)$ is dense in $C_1
\ot C_2$.
\end{proposition}

\bigskip

We are now in a position to define the universal \cst-algebra associated to the
$^*$-algebra $A$ together with a comultiplication on it.

\begin{definition}
We define the norm $\|.\|_u$ on $A$ such that
$$\|a\|_u = \sup \{ \, \|\phi(a)\| \mid C \text{ a \cst-algebra and } \phi
\text{ a $^*$-homomorphism from } A \text{ into } C \, \} \ .$$
Then $\|.\|_u$ is a \cst-norm on $A$, we define $A_u$ to be the completion of $A$ with
respect to $\|.\|_u$,  so $A_u$ is a \cst-algebra.
\end{definition}

\begin{remark} \rm
Proposition \ref{prop2.1} implies that the norm is finite. The regular representation
$\pi_r$ guarantees that $\|.\|_u$ is a norm and not merely a semi-norm.
\end{remark}

\begin{notation}
We define $\pi_u$ to be the identity mapping from $A$ into $A_u$. Hence, $\pi_u$ is an
injective $^*$-homomorphism from $A$ into $A_u$ such that $\pi_u(A)$ is dense in $A_u$.
\end{notation}

Restating the definition of the norm in terms of the mapping $\pi_u$, we have the following
equality.

\begin{proposition}
Consider $a \in A$. Then we have the following equality :
$$\|\pi_u(a)\|  =  \sup\ \{ \, \|\phi(a)\| \mid C \text{ a \cst-algebra and } \phi \text{ a
$^*$-homomorphism from } A \text{ into } C \, \} \ .$$
\end{proposition}

By the definition of the norm, we have the following universality property.

\begin{proposition} \label{prop2.3}
Consider a \cst-algebra $C$ and a $^*$-homomorphism  $\phi$ from $A$
into $C$. Then there exists a unique $^*$-homomorphism $\th$ from
$A_u$ into $C$ such that $\th \! \circ \! \pi_u = \phi$.
\end{proposition}

\bigskip

If  we apply proposition \ref{prop2.2} with $\phi$ equal to $\pi_u$
and combine this with the previous universality property, we get our
comultiplication on $A_u$ :

\begin{definition} \label{def2.1}
There exists a unique non-degenerate $^*$-homomorphism $\de_u$ from $A_u$ into $M(A_u \ot
A_u)$ such that $(\pi_u \od \pi_u)(x) \, \de_u(\pi_u(a)) = (\pi_u \od \pi_u)(x \, \de(a))$
and $\de_u(\pi_u(a))\,(\pi_u \od \pi_u)(x)= (\pi_u \od \pi_u)(\de(a)\,x)$ for every $a
\in A$ and $x \in A \od A$.
\end{definition}

As in section 2 of \cite{Kus}, this definition implies for instance that $(\pi_u(a) \ot
1)\de_u(\pi_u(b)) = (\pi_u \od \pi_u)((a \ot 1)\de(b))$ for every $a,b \in A$.

\medskip

As before, this implies the next theorem.

\begin{theorem}
The mapping $\de_u$ is a  non-degenerate $^*$-homomorphism from $A_u$
into $M(A_u \ot A_u)$ such that :
\begin{enumerate}
\item $(\de_u \ot \io)\de_u = (\io \ot \de_u)\de_u$
\item The vector spaces $\de_u(A_u)(A_u \ot 1)$ and $\de_u(A_u)(1 \ot A_u)$ are dense subspaces
of $A_u \ot A_u$.
\end{enumerate}
\end{theorem}

We call $(A_u,\de_u)$ the universal bi-\cst-algebra associated to $(A,\de)$.

\bigskip

Thanks to the universality property of $A_u$ (proposition
\ref{prop2.3}), the counit on $A$ can also be extended to $A_u$ :

\begin{definition}
We define $\vep_u$ as the unique $^*$-homomorphism from $A_u$ into $\C$ such that $\vep_u \!
\circ \! \pi_u = \vep$.
\end{definition}

It is not difficult to extend the counit property from $A$ to $A_u$ :

\begin{proposition}
We have for every $a \in A_u$ that $(\vep_u \ot \io)\de_u(a) = (\io \ot \vep_u)\de_u(a) = a$.
\end{proposition}

This implies immediately that the mapping $\de_u$ is injective.

\bigskip

In the rest of this paper, we will need regularly to switch between
$A_u$ and $A_r$. The regular representation $\pi_r$ together with the
universality property (proposition \ref{prop2.3}) allows us to make
the bridge :

\begin{definition}
We define $\pi$ to be the unique $^*$-homomorphism from $A_u$ into $A_r$ such that $\pi \!
\circ \! \pi_u = \pi_r$.
\end{definition}

Notice that $\pi$ here is something different than in \cite{Kus} !

\medskip

By definition \ref{def2.1} and lemma 2.8 of \cite{Kus}, we have of course that $(\pi \ot
\pi)\de_u = \de_r \,\pi$.

\section{The left regular corepresentation}

In this section, we introduce the left regular corepresentation of
$(A_u,\de_u)$. This corepresentation will be essential to transform
objects connected with $(A_r,\de_r)$ to objects connected with
$(A_u,\de_u)$. The first application of this principle can be found in
the next section.

\bigskip

Remember that $A_u \ot H$ is a Hilbert-\cst-module in a natural way
and let $\cL(A_u \ot H)$ denote the set of adjointable mappings from
$A_u \ot H$ into $A_u \ot H$. Remember also that $\cL(A_u \ot H) =
M(A_u \ot B_0(H))$.

\begin{definition}
There exists a unique unitary $V \in \cL(A_u \ot H)$ such that $V (\pi_u \od \la)(\de(b)(a
\ot 1)) = \pi_u(a) \ot \la(b) $ for every $a,b \in A$. We call $V$ the left regular
corepresentation of $(A_u,\de_u)$.
\end{definition}

This definition is possible because $\de(A)(A \ot 1) = A \od A$ and because of the left
invariance of $\vfi$.

\bigskip

In the same way as in proposition 2.2 of \cite{Kus}, we get the
following alternative formula for $V$ :

\begin{result}
We have for every $a,b \in A$ that $V (\pi_u(a) \ot \la(b)) = (\pi_u \od \la)\bigl((S^{-1}
\od \io)(\de(b)) \, (a \ot 1)\bigr)$.
\end{result}

\bigskip

Equation (1) of \cite{Kus} has its obvious variant in the universal
case :

\begin{lemma}  \label{lem3.1}
We have for every $a,b \in A$ that
$(\io \od \om_{\la(a),\la(b)})(V) = \pi_u\bigl( (\io \od \vfi)(\de(b^*)(1 \ot a))\bigr)$.
\end{lemma}
\begin{demo}
We have for every $c,d \in A$ that
\begin{eqnarray*}
& & \pi_u(d)^* \, (\io \od \om_{\la(a),\la(b)})(V) \, \pi_u(c)
= \langle V (\pi_u(c) \ot \la(a)) , \pi_u(d) \ot \la(b) \rangle \\
& & \spat = \langle \pi_u(c) \ot \la(a) , V^* (\pi_u(d) \ot \la(b)) \rangle
 = \langle \pi_u(c) \ot \la(a) , (\pi_u \od \la)(\de(b)(d \ot 1)) \rangle \\
& & \spat = (\pi_u \od \vfi)((d^* \ot 1)\de(b^*)(c \ot a))
= \pi_u \bigl( (\io \od \vfi)((d^* \ot 1)\de(b^*)(c \ot a))\bigr) \\
& & \spat = \pi_u(d^* (\io \od \vfi)(\de(b^*)(1 \ot a)) \,c)
= \pi_u(d)^* \, \pi_u\bigl((\io \od \vfi)(\de(b^*)(1 \ot a))\bigr) \, \pi_u(c) \ .
\end{eqnarray*}
This implies that $(\io \od \om_{\la(a),\la(b)})(V) = \pi_u\bigl((\io
\od \vfi)(\de(b^*)(1 \ot a))\bigr)$.
\end{demo}

As a consequence, we find the following result.

\begin{corollary}
The set  $\{\, (\io \od \om)(V) \mid \om \in B_0(H)^* \,\}$ is a dense
subspace of $A_u$.
\end{corollary}

\bigskip

Lifting $V$ to $H$ gives $W$ :

\begin{result}
We have that $(\pi \ot \io)(V) = W$.
\end{result}
\begin{demo}
Choose $a,b \in A$. Using lemma \ref{lem3.1} and equation (1) of \cite{Kus}, we get that
\begin{eqnarray*}
& & (\io \ot \om_{\la(a),\la(b)})((\pi \ot \io)(V))
= \pi( (\io \ot \om_{\la(a),\la(b)})(V))
= \pi\bigl(\,\pi_u(\,(\io \od \vfi)(\de(b^*)(1 \ot a))\,)\,\bigr) \\
& & \spat = \pi_r\bigl((\io \od \vfi)(\,\de(b^*)(1 \ot a))\bigr)
= (\io \ot \om_{\la(a),\la(b)})(W)\ .
\end{eqnarray*}
Consequently, $(\pi \ot \io)(V) = W$.
\end{demo}

\medskip

\begin{notation} \label{not3.1}
For every $p,q \in H$, we denote $\om_{p,q}' \in A_u^*$ such that $\om_{p,q}'(x) = \langle
\pi(x) p , q \rangle $ for every $x \in A_u$.
\end{notation}

Because $(\pi \ot \io)(V) = W$, lemma 2.20 of \cite{Kus} implies the
following equality :

\begin{result}
We have for every $a,b \in A$ that $(\om_{\la(a),\la(b)}' \ot \io)(V)
= \pih_r(a \vfi b^*)$.
\end{result}

\begin{corollary} \label{cor3.1}
We have that $\ah_r$ is a subset of the closure of $\{\, (\om \ot \io)(V) \mid \om \in
A_u^* \,\}$ in $B_0(H)$.
\end{corollary}

\bigskip

The left regular corepresentation lifts the algebraic universal corepresentation $\aU$ from $A
\od \ah$ to $A_u \ot \ah_r$ :

\begin{proposition} \label{prop3.1}
We have for every $a,b \in A$ and $x \in A \od \ah$ that $V \, (\pi_u \od \pih_r)(x) =
(\pi_u \od  \pih_r)(\aU  \, x)$ and $(\pi_u \od \pih_r)(x) \, V  = (\pi_u \od \pih_r)(x \,
\aU)$
\end{proposition}

The proof of the first equality is completely the same as the proof of lemma 2.21 of
\cite{Kus}. The second one follows from the first.

\begin{corollary}
We have the equalities $V \,(\pi_u(A) \od \pih_r(\ah)) = \pi_u(A) \od \pih_r(\ah)$ and
$(\pi_u(A) \od \pih_r(\ah)) \, V = \pi_u(A) \od \pih_r(\ah)$.
\end{corollary}

\begin{corollary}
The element $V$ belongs to $M(A_u \ot \ah_r)$.
\end{corollary}

\medskip

Because $(\de \od \io)(\aU) = \aU_{13} \, \aU_{23}$ and $(\io \od \deh)(\aU) =
\aU_{12} \, \aU_{13}$,  proposition \ref{prop3.1} implies the following
result :

\begin{proposition} \label{prop3.2}
We have that $(\de_u \ot \io)(V) = V_{13} \, V_{23}$ and $(\io \od
\deh_r)(V) = V_{12} \, V_{23}$.
\end{proposition}

So $V$ is a corepresentation of $(A_u,\de_u)$ on $A_r$ and a corepresentation of
$(A_r,\de_r)$ on $A_u$.

\bigskip

As usual, $V$ implements the comultiplication :

\begin{proposition} \label{prop3.3}
We have for every $a \in A_u$ that $(\io \ot \pi)\de_u(a) = V^* (1 \ot \pi(a)) V $.
\end{proposition}
\begin{demo}
Choose $x \in A$. Take $b,c \in A$. Then there exists $e \in A$ such
that $(e \ot 1)\de(c)(b \ot 1) = \de(c)(b \ot 1)$. This implies that
\begin{eqnarray*}
& &  V^* (1 \ot \pi(\pi_u(x))) \, (\pi_u(b) \ot \la(c)) = V^* (1 \ot
\pi_r(x)) (\pi_u(b) \ot \la(c)) \\
 & & \spat = V^* (\pi_u(b)  \ot \la(x c))
= (\pi_u \od \la)(\de(x c) (b \ot 1)) \\
& & \spat = (\pi_u \od \la)(\de(x) \de(c)(b \ot 1))
= (\pi_u \od \la)(\de(x) (e \ot 1) \de(c) (b \ot 1)) \\
& & \spat = (\pi_u \od \pi_r)(\de(x)(e \ot 1))  \, (\pi_u \od \la)(\de(c)(b \ot 1)) \\ & &
\spat = (\io \ot \pi)\bigl( (\pi_u \od \pi_u)(\de(x)(e \ot 1)) \bigr) \, (\pi_u \od
\la)(\de(c)(b \ot 1)) \\ & & \spat = (\io \ot \pi)(\de_u(\pi_u(x))(\pi_u(e) \ot 1)) \, (\pi_u
\od \la)(\de(c)(b \ot 1)) \\ & & \spat = (\io \ot \pi)\de_u(\pi_u(x)) \, (\pi_u(e) \ot 1) \,
(\pi_u \od \la)(\de(c)(b \ot 1))  \\ & & \spat = (\io \ot \pi)\de_u(\pi_u(x)) \, (\pi_u \od
\la)(\de(c)(b \ot 1))
= (\io \ot \pi)\de_u(\pi_u(x)) \, V^* (\pi_u(b) \ot \la(c)) \ .
\end{eqnarray*}
So we see that $V^* (1 \ot \pi(\pi_u(x))) = (\io \ot \pi)\de_u(\pi_u(x)) \, V^*$. Because
$\pi_u(A)$ is dense in $A_u$, we get that $V^* (1 \ot \pi(a)) = (\io \ot \pi)(\de_u(a))
V^*$.
\end{demo}

\bigskip

In the rest of this section, we prove a formula which formally says that $(\io \ot \hat{S})(V)
= V^*$. It will prove useful in the next section. First, we need a lemma.

\begin{lemma}
Consider $a,b \in A$ and $\om \in A_u^*$. Then
$$(\pi_u(a) \om \pi_u(b) \ot \io)(V^*)\,\la(x)
= \la\bigl((\om \! \circ \! \pi_u \od \io)((b \ot 1)\de(x)(a \ot 1))\bigr)$$ for every
$x \in A$.
\end{lemma}
\begin{demo}
We have for every $y \in A$ that
\begin{eqnarray*}
& & \langle (\pi_u(a) \om  \pi_u(b) \ot \io)(V^*)\, \la(x) , \la(y) \rangle
= \om(\langle V^* (\pi_u(a) \ot \la(x)) , \pi_u(b^*) \ot \la(y) \rangle) \\
& & \spat = \om( \langle (\pi_u \od \la)(\de(x)(a \ot 1)) ,  \pi_u(b^*) \ot \la(y) \rangle )
= \om\bigl( (\pi_u \od \vfi)( (b \ot y^*) \de(x)(a \ot 1) )\bigr) \\
& & \spat = \vfi\bigl( (\om \! \circ \! \pi_u \od \io)((b \ot y^*)\de(x)(a \ot 1))\bigr)
= \vfi\bigl( y^* (\om \! \circ \! \pi_u \od \io)((b \ot 1)\de(x)(a \ot 1))\bigr) \\
& & \spat = \langle \la\bigl((\om \! \circ \! \pi_u \od \io)((b \ot 1)\de(x)(a \ot 1))\bigr) ,
\la(y) \rangle \ .
\end{eqnarray*}
This implies that $(\pi_u(a) \om  \pi_u(b) \ot \io)(V^*) \, \la(x)
= \la\bigl((\om \! \circ \! \pi_u \od \io)((b \ot 1)\de(x)(a \ot 1))\bigr)$.
\end{demo}

\begin{proposition}
Consider $p \in D(\nab^\frac{1}{2})$ and $q \in D(\nab^{-\frac{1}{2}})$. Then $(\io \ot
\om_{J \nab^\frac{1}{2} p , J \nab^{-\frac{1}{2}} q})(V^*) = (\io \ot \om_{q,p})(V)$.
\end{proposition}
\begin{demo}
Choose $a,b \in A$ and $\om \in A_u^*$. Take $x \in A$. We know from
the remarks after proposition \ref{prop1.1} that $\la(x)$ belongs to
$D(\nab^\frac{1}{2})$ and $J \nab^\frac{1}{2} \la(x) = \la(x^*)$.
Therefore
\begin{eqnarray*}
& & \hspace{-.5cm} (\pi_u(a) \om  \pi_u(b^*))\bigl((\io \ot \om_{J \nab^\frac{1}{2}
\la(x) , J \nab^{-\frac{1}{2}} q})(V^*)\bigr)
=  (\pi_u(a) \om  \pi_u(b^*))\bigl((\io \ot \om_{\la(x^*), J \nab^{-\frac{1}{2}} q})(V^*)
\bigr) \\
& & =  \langle (\pi_u(a) \om  \pi_u(b^*) \ot \io)(V^*) \, \la(x^*)
,  J \nab^{-\frac{1}{2}} q \rangle
=  \langle \la\bigl((\om \! \circ \! \pi_u \od \io)((b^* \ot 1)\de(x^*)(a \ot 1))\bigr) ,
J \nab^{-\frac{1}{2}} q \rangle  \hspace{0.65cm} (*)
\end{eqnarray*}
where we used the previous lemma in the last equality. By the remarks after proposition
\ref{prop1.1}, we know that also $\la\bigl((\overline{\om} \! \circ
\! \pi_u \od \io)((a^* \ot 1)\de(x)(b \ot 1))\bigr)$ belongs to $D(\nab^\frac{1}{2})$ and
$$J \nab^\frac{1}{2} \la\bigl((\overline{\om} \! \circ \! \pi_u \od \io)((a^* \ot 1)
\de(x)(b \ot 1))\bigr) = \la\bigl((\om \! \circ \! \pi_u \od \io)((b^* \ot 1)\de(x^*)
(a \ot 1))\bigr) \ . $$
Substituting this in equality (*), we get that
$$  (\pi_u(a) \om  \pi_u(b^*))\bigl((\io \ot \om_{J \nab^\frac{1}{2} \la(x) ,
J \nab^{-\frac{1}{2}} q})(V^*)\bigr)
= \langle J \nab^\frac{1}{2} \la\bigl((\overline{\om} \! \circ \! \pi_u \od \io)
((a^* \ot 1)\de(x)(b \ot 1))\bigr) , J \nab^{-\frac{1}{2}} q \rangle  \ .$$
This implies that
\begin{eqnarray*}
& &  (\pi_u(a) \om  \pi_u(b^*))\bigl((\io \ot \om_{J \nab^\frac{1}{2} \la(x) , J
\nab^{-\frac{1}{2}} q})(V^*)\bigr) \\
& & \spat = \langle q ,  \la\bigl((\overline{\om} \! \circ \! \pi_u \od \io)
((a^* \ot 1)\de(x)(b \ot 1))\bigr)\rangle
= \langle q , (\pi_u(b) \overline{\om} \, \pi_u(a^*) \ot \io)(V^*) \la(x) \rangle \\
& & \spat = \langle (\pi_u(a) \om  \pi_u(b^*) \ot \io)(V)\,q , \la(x) \rangle
= (\pi_u(a) \om  \pi_u(b^*))((\io \ot \om_{q,\la(x)})(V)) \ .
\end{eqnarray*}
From this all, we conclude that
$$(\io \ot \om_{J \nab^\frac{1}{2} \la(x) , J \nab^{-\frac{1}{2}} q})(V^*)
= (\io \ot \om_{q,\la(x)})(V) \ . $$
Because $\la(A)$ is a core for $\nab^\frac{1}{2}$ (see again the
remarks after proposition \ref{prop1.1}), the proposition follows.
\end{demo}

\begin{corollary}  \label{cor3.2}
Consider $p \in D(\nab^\frac{1}{2})$ and $q \in D(\nab^{-\frac{1}{2}})$. Then $(\io \ot
\om_{p,q})(V^*) = (\io \ot \om_{J \nab^{-\frac{1}{2}} q , J \nab^\frac{1}{2} p})(V)$.
\end{corollary}

\section{Lifting bi-\cst-isomorphisms from the reduced to the universal \cst-algebra}
\label{art4}

In \cite{Kus}, we saw that the reduced quantum group $(A_r,\de_r)$
possesses a very rich structure :  a left and right Haar weight which
are KMS, a polar decomposition of the antipode, ...

We will show that this rich structure on $(A_r,\de_r)$ can be transported to
$(A_u,\de_u)$. In this section, we will provide the most important method to do so.

\bigskip

For the rest of this section, we  fix $^*$-automorphisms $\al,\be$ on $A_r$ such that
$(\al \ot \be) \de_r = \de_r  \, \al$.

\begin{lemma} \label{lem4.2}
We have that $(\be \ot \be) \de_r = \de_r \, \be$.
\end{lemma}
\begin{demo}
We have that
\begin{eqnarray*}
& & (\al \ot \de_r \, \be)\de_r = (\io \ot \de_r)\de_r \,\al = (\de_r \ot \io)\de_r \, \al \\
& & \spat = (\de_r \ot \io)(\al \ot \be)\de_r = (\al \ot \be \ot \be)(\de_r \ot \io)\de_r
=  (\al \ot (\be \ot \be)\de_r)\de_r
\end{eqnarray*}
which implies that $\de_r \, \be = (\be \ot \be)\de_r$.
\end{demo}

Hence, using proposition 7.1 of \cite{Kus}, we get the following
result.

\begin{corollary}
There exists a unique strictly positive number $r$ such that $\vfi_r \, \be
= r \, \vfi_r $.
\end{corollary}

This relative invariance can be extended to $\al$ :

\begin{proposition} \label{prop4.2}
We have also that $\vfi_r \, \al = r \, \vfi_r$.
\end{proposition}
\begin{demo}
Choose $a \in \Mfir^+$.

Take $\om \in (A_r)_+^*$. Then
$$(\om \ot \io)\de_r(\al(a)) = \be( (\om \al \ot \io)\de_r(a))\ . \hspace{1cm} (*)$$
Because $a$ belongs to $\Mfir^+$ and $\vfi_r$ is left invariant, we get that $(\om \al \ot
\io)\de_r(a)$ belongs to $\Mfir^+$ and
$$\vfi_r((\om \al \ot \io)\de_r(a)) = \om(\al(1)) \, \vfi_r(a) = \om(1) \, \vfi_r(a) \ . $$
Because $\vfi_r \, \be = r \, \vfi_r$, this implies that $\be((\om \al \ot \io)\de_r(a))$
belongs to $\Mfir^+$ and
$$\vfi_r\bigl(\be((\om \al \ot \io)\de_r(a))\bigr) = r \, \vfi_r((\om \al \ot \io)\de_r(a)) =
r \, \om(1) \, \vfi_r(a) \ .$$
Looking at equality (*), this implies that $(\om \ot \io)\de_r(\al(a))$ belongs to $\Mfir^+$
and $$\vfi_r\bigl((\om \ot \io)\de_r(\al(a))\bigr) = r \, \om(1) \, \vfi_r(a) \ .$$

By theorem of 3.11 of \cite{JK4}, we get that $\al(a)$ belongs to $\Mfir^+$.

\medskip

Using the left invariance of $\vfi_r$, we get moreover for every $\th \in (A_r)_+^*$ that
$$\th(1) \, \vfi_r(\al(a)) = \vfi_r\bigl((\th \ot \io)\de(\al(a))\bigr)
= r \, \th(1) \, \vfi_r(a) $$
which implies that $\vfi_r(\al(a)) = r \, \vfi_r(a)$.
\end{demo}

\begin{notation}
We define the unitary elements $u,v \in B(H)$ such that $u \la_r(a) = r^{-\frac{1}{2}}
\la_r(\al(a))$ and $v \la_r(a) = r^{-\frac{1}{2}} \la_r(\be(a))$ for every $a \in \Nfir$.
\end{notation}

Then we have the following results :

\begin{result}
We have for every $x \in A_r$ that $\al(x) = u x u^*$ and
$\be(x) = v x v^*$.
\end{result}

\begin{result} \label{res4.3}
We have the following commutation relations :
\begin{enumerate}
\item $J u = u J$ and $\nab u = u \nab$
\item $J v = v J$ and $\nab v = v \nab$
\end{enumerate}
\end{result}

\bigskip

The equality $(\al \ot \be)\de_r = \de_r \, \al$ implies the following commutation relations
between $u$, $v$ and $W$.

\begin{proposition} \label{prop4.1}
The following commutation relations hold :
\begin{enumerate}
\item $W(u \ot v) = (u \ot u)W$
\item $W(v \ot v) = (v \ot v)W$
\end{enumerate}
\end{proposition}
\begin{demo}
\begin{enumerate}
\item Choose $a,b \in \Nfir$.

Then $\al(a), \al(b)$ belong to $\Nfir$ and $u \la_r(a) = r^{-\frac{1}{2}} \, \la_r(\al(a))$
and $u \la_r(b) = r^{-\frac{1}{2}} \,
\la_r(\al(b))$.

Because $\al(a), \al(b)$ belong to $\Nfir$, proposition \ref{prop1.2}
implies that $\de(\al(b))(\al(a) \ot 1)$ belongs to $\cN_{\vfi_r \ot
\vfi_r}$ and
$$ (\la_r \ot \la_r)(\de(\al(b))(\al(a) \ot 1))
= W^* (\la_r(\al(a)) \ot \la_r(\al(b)))
= r \,\, (W^* (u \ot u)) (\la_r(a) \ot \la_r(b)) \ . \hspace{1cm} \text{(a)} $$

Because $a, b$ belong to $\Nfir$,  proposition \ref{prop1.2} implies
that also $\de(b)(a \ot 1)$ belongs to $\cN_{\vfi_r \ot \vfi_r}$  and
$$ (\la_r \ot \la_r)(\de(b)(a \ot 1)) = W^* (\la_r(a) \ot \la_r(b)) \ . $$
This implies that $(\al \ot \be)(\de(b)(a \ot 1))$ belongs to $\cN_{\vfi_r \ot \vfi_r}$ and
\begin{eqnarray*}
(\la_r \ot \la_r)\bigl((\al \ot \be)(\de(b)(a \ot 1))\bigr) & = &
r \,\, (u \ot v) (\la_r \ot \la_r)(\de(b)(a \ot 1)) \\
& = & r \,\, ((u \ot v) W^*) (\la_r(a) \ot \la_r(b)) \ .
\hspace{1cm} \text{(b)}
\end{eqnarray*}
Because  $(\al \ot \be)(\de(b)(a \ot 1)) = \de(\al(b))(\al(a) \ot 1)$, equalities (a) and (b)
imply that
$$(W^* (u \ot u)) (\la_r(a) \ot \la_r(b)) = ((u \ot v) W^*) (\la_r(a) \ot \la_r(b)) \ .$$
So we get that $W^* (u \ot u) = (u \ot v) W^*$.
\item The second equality is proven in the same way as the first one.
\end{enumerate}
\end{demo}

\begin{corollary} \label{cor4.1}
We have for every $\om \in B_0(H)^*$ that $\al((\io \ot \om)(W)) =
(\io \ot v \om u^*)(W)$ \newline and $\be((\io \ot \om)(W)) = (\io \ot
v \om v^*)(W)$.
\end{corollary}

\bigskip

Using these commutation relations, we get the following result.

\begin{result} \label{res4.1}
The following equalities hold :
\begin{enumerate}
\item $\ah_r = u^* \ah_r v = v^* \ah_r u$
\item $\ah_r = u^* \ah_r u = v^* \ah_r v$
\end{enumerate}
\end{result}
\begin{demo}
By proposition \ref{prop4.1}, we have for every $\om \in B_0(H)^*$ that
$$u^* (\om \ot \io)(W) v = (\om \ot \io)((1 \ot u^*) W (1 \ot v))
= (\om \ot \io)((u \ot 1)W (u^* \ot 1))
= (u^* \om u \ot \io)(W) \ .$$
Because $\ah_r$ is the closure of $\{\, (\om \ot \io)(W) \mid \om \in B_0(H)^* \,\}$, this
equality implies that $u^* \ah_r v = \ah_r$.

Taking the adjoint of this equation gives us that $v^* \ah_r u =
\ah_r$.

We have moreover that
$$u^* (\ah_r)^2 u = (u^* \ah_r v)(v^* \ah_r u) = (\ah_r)^2 \ . $$
Because $(\ah_r)^2$ is a dense subset of $\ah_r$, this implies that
$u^* \ah_r u = \ah_r$. In a similar way, we get that $v^* \ah_r v = \ah_r$.
\end{demo}

Using this last result, it is not difficult to infer the following one.

\begin{result}  \label{res4.2}
The following equalities hold :
\begin{enumerate}
\item $M(\ah_r) = u^* M(\ah_r) v = v^* M(\ah_r) u$
\item $M(\ah_r) = u^* M(\ah_r) u = v^* M(\ah_r) v$
\end{enumerate}
\end{result}

\bigskip

This last result will imply the following lemma.

\begin{lemma}  \label{lem4.1}
Consider $\om,\th \in B_0(H)^*$ such that $(\io \ot \om)(V) = (\io \ot \th)(V)$. Then
$(\io \ot v \om u^*)(V) = (\io \ot v \th u^*)(V)$ and $(\io \ot v \om v^*)(V) = (\io \ot v
\th v^*)$.
\end{lemma}
\begin{demo}
We have for every $\eta \in A_u^*$ that
$$\om((\eta \ot \io)(V)) = \eta((\io \ot \om)(V)) = \eta((\io \ot \th)(V))
= \th((\eta \ot \io)(V)) \ . $$
Therefore, corollary \ref{cor3.1} implies that $\om(x) = \th(x)$ for every
$x \in \ah_r$. Therefore, $\om(x) = \th(x)$ for every $x \in M(\ah_r)$.

\medskip

Take $\eta \in A_u^*$. Because $V$ belongs to $M(A_u \ot \ah_r)$, we
get that $(\eta \ot \io)(V)$ belongs to $M(\ah_r)$. Therefore, the previous result
implies that $u^* (\eta \ot \io)(V) v$ belongs to $M(\ah_r)$.

By the first part of the proof, we know that
$$\om(u^* (\eta \ot \io)(V) v) = \th(u^* (\eta \ot \io)(V) v) \ ,$$
implying that $\eta((\io \ot v \om u^*)(V)) = \eta((\io \ot v \th u^*)(V))$.

Consequently, $(\io \ot v \om u^*)(V) = (\io \ot v \th u^*)(V)$.
Similarly, $(\io \ot v \om v^*)(V) = (\io \ot v \th v^*)(V)$.
\end{demo}

\medskip

This allows us to prove the following proposition :

\begin{proposition}
There exists a unique linear mapping $F$ from $\pi_u(A)$ into $A_u$ such
that $F((\io \ot \om_{p,q})(V))$ \newline
$= (\io \ot \om_{v p , u q})(V)$ for every $p,q \in \la(A)$.
\end{proposition}
\begin{demo}
By the previous lemma, we have for every $p_1,\ldots\!,p_m , q_1,\ldots\!,q_m
\ , \ r_1,\ldots\!,r_n , s_1,\ldots\!,s_n \in \la(A)$ such that
$$\sum_{i=1}^m (\io \ot \om_{p_i,q_i})(V) =
\sum_{j=1}^n (\io \ot \om_{r_j,s_j})(V)$$
that  $$\sum_{i=1}^m (\io \ot \om_{v p_i,u q_i})(V) =
\sum_{j=1}^n (\io \ot \om_{v r_j,u s_j})(V) \ .$$
This implies that we can define a mapping $F$ from $\pi_u(A)$ into $A_u$ such that
$$F\bigl(\,\sum_{i=1}^m (\io \ot \om_{p_i,q_i})(V)\,\bigr) =
\sum_{i=1}^m (\io \ot \om_{v p_i,u q_i})(V) $$
for every $m \in \N$ and $p_1,\ldots\!,p_m , q_1,\ldots\!,q_m \in
\la(A)$. It is easy to check that $F$ is linear.
\end{demo}

\begin{result}
The mapping $F$ is multiplicative.
\end{result}
\begin{demo}
Choose $p,q,r,s \in \la(A)$.

Define $\om \in B_0(H)^*$ such that $\om(x) = \langle W(x \ot 1)W^* (p \ot r) ,
q \ot s \rangle$ for every $x \in B_0(H)$. Using proposition \ref{prop3.2}, we get that
\begin{eqnarray*}
& & (\io \ot \om_{p,q})(V) \, (\io \ot \om_{r,s})(V)
= (\io \ot \om_{p,q} \ot \om_{r,s})(V_{12} V_{13}) \\
& & \spat = (\io \ot \om_{p,q} \ot \om_{r,s})(W_{23} V_{12} W_{23}^*)
= (\io \ot \om)(V) \ .
\end{eqnarray*}
Therefore, lemma \ref{lem4.1} implies that
$$F\bigl(\,(\io \ot \om_{p,q})(V) \, (\io \ot \om_{r,s})(V)\,\bigr)
= (\io \ot v \om u^*)(V) \ . \hspace{1cm} \text{(*)}$$
Now we have for every $x \in B_0(H)$ that
\begin{eqnarray*}
& & (v \om u^*)(x)  =  \om(u^* x v) = \langle W(u^* x v \ot 1) W^* (p \ot r) , q \ot s
\rangle \\
& & \spat = \langle W (u^* \ot v^*) (x \ot 1) (v \ot v) W^* (p \ot r) , q \ot s \rangle
= \langle W (x \ot 1)  W^* (v p \ot v r) , u q \ot u s \rangle
\end{eqnarray*}
where we used proposition \ref{prop4.1} in the last equality. Combining this result with
equality (*), we see that
\begin{eqnarray*}
& & F\bigl(\,(\io \ot \om_{p,q})(V) \, (\io \ot \om_{r,s})(V)\,\bigr)
= (\io \ot \om_{v p , u q} \ot \om_{v r , u s})(W_{23} V_{12} W_{23}^*)
= (\io \ot \om_{v p , u q} \ot \om_{v r , u s})(V_{12} V_{13}) \\
& & \spat = (\io \ot \om_{v p , u q})(V) \, (\io \ot \om_{v r , u s})(V)
= F((\io \ot \om_{p,q})(V)) \, F((\io \ot \om_{r , s})(V)) \ .
\end{eqnarray*}
The result follows by linearity.
\end{demo}

\begin{result}
The mapping $F$ is selfadjoint.
\end{result}
\begin{demo}
Choose $p,q \in \la(A)$. We know from the remarks after proposition \ref{prop1.1} that $p$
belongs to $D(\nab^{-\frac{1}{2}})$ and that $q$ belongs to $D(\nab^\frac{1}{2})$.
Therefore, corollary \ref{cor3.2} implies that
$$(\io \ot \om_{p,q})(V)^* = (\io \ot \om_{q,p})(V^*)
= (\io \ot \om_{J \nab^{-\frac{1}{2}} p , J \nab^\frac{1}{2} q})(V) \ . $$
Using lemma \ref{lem4.1}, this implies that
$$F((\io \ot \om_{p,q})(V)^*)
= (\io \ot \om_{v J \nab^{-\frac{1}{2}} p , u J \nab^\frac{1}{2} q})(V) \ . $$
Result \ref{res4.3} implies that $v p$ belongs to
$D(\nab^{-\frac{1}{2}})$, that $u q$ belongs to $D(\nab^\frac{1}{2})$
and $J \nab^{-\frac{1}{2}} v p = v J \nab^{-\frac{1}{2}} p$ and $J
\nab^\frac{1}{2} u q = u J \nab^\frac{1}{2} q$. Hence,
\begin{eqnarray*}
& & F((\io \ot \om_{p,q})(V)^*)
= (\io \ot \om_{J \nab^{-\frac{1}{2}} v p,J \nab^\frac{1}{2} u q})(V) \\
& & \spat \stackrel{(*)}{=} (\io \ot \om_{u q , v p})(V^*)
= (\io \ot \om_{v p, u q})(V)^* = F((\io \ot \om_{p,q})(V))^*
\end{eqnarray*}
where corollary \ref{cor3.2} was used once again in equality (*). The result follows by
linearity.
\end{demo}

\bigskip

Consequently, we have proven that $F$ is a $^*$-homomorphism from
$\pi_u(A)$ into $A_u$. By proposition \ref{prop2.3}, this justifies
the following definition :

\begin{definition}
There exists a unique $^*$-homomorphism $\al_u$ from $A_u$ into $A_u$ such that
$\al_u((\io \ot \om_{p,q})(V))$ \newline $= (\io \ot \om_{v p , u q})(V)$ for every $p,q
\in \la(A)$.
\end{definition}

Then we get immediately the following proposition.

\begin{result} \label{res4.4}
The equality $(\al_u \ot \io)(V) = (1 \ot u^*) V (1 \ot v)$ holds.
\end{result}

\begin{corollary} \label{cor4.2}
We have for every $\om \in B_0(H)^*$ that $\al_u((\io \ot \om)(V))
= (\io \ot v \om u^*)(V)$.
\end{corollary}

\bigskip

\begin{proposition}
The mapping $\al_u$ is an $^*$-automorphism on $A_u$.
\end{proposition}
\begin{demo}
We have of course also the equality $(\al^{-1} \ot \be^{-1})\de = \de \al^{-1}$. Therefore we
can do the same thing for $\al^{-1}$ as we did for $\al$. So we get the existence of a $^*$
homomorphism $\th$ from $A_u$ into $A_u$ such that $\th((\io \ot \om)(V))
= (\io \ot v^* \om u)(V)$ for every $\om \in B_0(H)$. Then we get immediately that
$$\al_u(\th((\io \ot \om)(V))) = \th(\al_u((\io \ot \om)(V))) = (\io
\ot \om)(V)$$ for every $\om \in B_0(H)^*$.

From this, it follows that $\al_u \! \circ \! \th = \th \! \circ \! \al_u = \io$. The
proposition follows.
\end{demo}

\bigskip

In a similar way, we can do the same things for $\be$. So we get the following result.

\begin{proposition}  \label{prop4.4}
There exists a unique $^*$-automorphism $\be_u$ on $A_u$ such that $\be_u((\io \ot \om)(V))$
\newline $= (\io \ot v \om v^*)(V)$ for every $\om \in B_0(H)^*$.
\end{proposition}

\begin{result}  \label{res4.5}
The equality $(\be_u \ot \io)(V) = (1 \ot v^*)V(1 \ot v)$ holds.
\end{result}

\medskip

Furthermore, the  commutation relation between $\al$, $\be$ and
$\de_r$ is transferred to the same commutation relation between
$\al_u$, $\be_u$ and $\de_u$ :

\begin{proposition}  \label{prop4.3}
We have that $(\al_u \ot \be_u) \de_u = \de_u \, \al_u$.
\end{proposition}
\begin{demo}
Using proposition \ref{prop3.2}, we have that
$$ ((\al_u \ot \be_u)\de_u \ot \io)(V)
= (\al_u \ot \be_u \ot \io)(V_{13} V_{23})
= (\al_u \ot \io)(V)_{13} \, (\be_u \ot \io)(V)_{23} \ .$$
Therefore, result \ref{res4.4} and \ref{res4.5} imply that
\begin{eqnarray*}
& & ((\al_u \ot \be_u)\de_u \ot \io)(V)
= [(1 \ot u^*)V(1 \ot v)]_{13} \, [(1 \ot v^*)V(1 \ot v)]_{23} \\
& & \spat = (1 \ot 1 \ot u^*) V_{13} (1 \ot 1 \ot v) \, (1 \ot 1 \ot v^*) V_{23}
(1 \ot 1 \ot v) = (1 \ot 1 \ot u^*)V_{13} V_{23} (1 \ot 1 \ot  v) \ .
\end{eqnarray*}
Hence, using proposition \ref{prop3.2} and result \ref{res4.4} once
more, we get that
\begin{eqnarray*}
& & ((\al_u \ot \be_u)\de_u \ot \io)(V)
= (1 \ot 1 \ot u^*)(\de_u \ot \io)(V) (1 \ot 1 \ot v) \\
& & \spat = (\de_u \ot \io)((1 \ot u^*)V(1 \ot v))
=   (\de_u \, \al_u \ot \io)(V) \ .
\end{eqnarray*}
This last equality implies for every $\om \in B_0(H)^*$ that
$$ (\al_u \ot \be_u)\bigl(\de_u((\io \ot\om)(V))\bigr)
=  \de_u\bigl(\al_u((\io \ot \om)(V))\bigr) \ .$$
From this all, we can conclude that $(\al_u \ot \be_u)\de_u = \de_u \, \al_u$.
\end{demo}

Because $(\pi \ot \io)(V) = W$,  corollaries \ref{cor4.1},
\ref{cor4.2} and  proposition \ref{prop4.4} imply easily the following
result.

\begin{result} \label{res4.6}
We have that $\pi \al_u = \al$ and $\pi \be_u = \be$.
\end{result}

\begin{remark} \rm
Later (proposition \ref{prop10.1}), we will prove that $\al_u$ and
$\be_u$ are uniquely determined by the properties $\pi \al_u = \al$
and $(\al_u \ot \be_u) \de = \de \al_u$.
\end{remark}

\bigskip

In a later section, we will also transform corepresentations from
$(A_r,\de_r)$ to corepresentations of $(A_u,\de_u)$.

\section{The algebraic nature of a unitary corepresentation of $(A_r,\de_r)$}
\label{art5}

We will prove in this section that every unitary corepresentation of $(A_r,\de_r)$ is of
an algebraic nature.

\medskip

For the  rest of this section, we  fix a \cst-algebra $C$ and a
unitary corepresentation $\cU$ of $(A_r,\de_r)$ on $C$. We want to
prove a generalization of proposition 7.7 of \cite{Kus}.  In order to
do so, we need to single out a special sub-$^*$-algebra of $C$.

\medskip

\begin{notation}
Define the set $B = \langle \, (\om_{p,q} \ot \io)(\cU) \mid p,q \in \la(A) \,\rangle$, so
$B$ is a subspace of $M(C)$.
\end{notation}

\bigskip

For every $\om \in \ah$, there exists a unique element $\th \in A_r^*$ such that $\th
\! \circ \! \pi_r = \om$ and we define $\tilde{\om} = \th$.
In this notation , we have for every $a,b \in A$ that $(a \vfi
b^*)\golf  =\om_{\la(a),\la(b)}$.

So we get that $B = \{\,(\tilde{\om} \ot \io)(\cU) \mid \om \in \ah \,\}$.

\medskip

It is not so difficult to check that $(\om_1 \om_2)\golf
= (\tilde{\om}_1 \ot \tilde{\om}_2)\de_r$ for every $\om_1,\om_2 \in \ah$.

Because $(\de_r \ot \io)(\cU) = \cU_{13} \, \cU_{23}$, this implies easily that the
mapping $\ah \rightarrow M(C) : \om \mapsto (\tilde{\om} \ot \io)(\cU)$ is an homomorphism
of algebras.

\medskip

Therefore, we get the following result :

\begin{result}
We have that $B^2 = B$.
\end{result}

\bigskip

This algebra $B$ is non-degenerate with respect to $C$ :

\begin{result}
We have that $B C$ and $C B$ are dense in $C$.
\end{result}
\begin{demo}
Choose $c \in C$. There exists $a,b \in A$ and $x \in A_r$ such that $\langle x \, \la(a) ,
\la(b)
\rangle = 1$. We have for every $d \in A$ and $y \in C$ that
$$(\om_{\la(a),\la(b)} \ot \io)(\cU (\pi_r(d) \ot y))
= (\om_{\pi_r(d) \la(a) , \la(b)} \ot \io)(\cU) \, y
= (\om_{\la(d a) , \la(b)} \ot \io)(\cU) \, y \in B C \ . $$
Because $\pi_r(A) \od C$ is dense in $A_r \ot C$, this implies for every $z \in A_r \ot C$
that $(\om_{\la(a),\la(b)} \ot \io)(\cU z)$ belongs to $\overline{B C}$. Because $\cU$ is
unitary, this implies that $(\om_{\la(a),\la(b)} \ot \io)(z)$ belongs to $\overline{B C}$
for every $z \in A_r \ot C$.

Because $c = (\om_{\la(a),\la(b)} \ot \io)(x \ot c)$, this implies that $c$ is an element
of $\overline{B C}$.

In a similar way, one proves that $C B$ is dense in $C$.
\end{demo}

\bigskip

We would like to prove that $\cU$ is an algebraic multiplier of $A \od B$. For this, we
need the following lemma.

\begin{lemma}
Consider elements $a_1,\ldots\!,a_m , b_1,\ldots\!,b_m , c_1,\ldots\!,c_m$ and
$p_1,\ldots\!,p_n ,  q_1,\ldots\!,q_n , r_1,\ldots\!,r_n$ in $A$ such that $\sum_{i=1}^m
a_i \ot b_i \rho(c_i^*) = \sum_{j=1}^n \de(r_j \rho(s_j^*))(t_j \ot 1)$. Then
$$\sum_{i=1}^m \cU \, (\pi_r(a_i) \ot (\om_{\la(b_i),\la(c_i)} \ot \io)(\cU))
= \sum_{j=1}^n \pi_r(t_j) \ot (\om_{\la(r_j),\la(s_j)} \ot \io)(\cU) \ .$$
\end{lemma}
\begin{demo}
Using the fact that $(\de_r \ot \io)(\cU) = \cU_{13} \,\cU_{23}$, we get that
\begin{eqnarray*}
\sum_{i=1}^m \cU \, (\pi_r(a_i) \ot (\om_{\la(b_i),\la(c_i)} \ot \io)(\cU))
& = & (\io \ot \om_{\la(b_i),\la(c_i)} \ot \io)(\cU_{13} \, \cU_{23}(\pi_r(a_i) \ot 1 \ot 1))
\\
& = & (\io \ot \om_{\la(b_i),\la(c_i)} \ot \io)((\de_r \ot \io)(\cU) (\pi_r(a_i) \ot 1 \ot 1))
\ . \hspace{1cm} \text{(*)}
\end{eqnarray*}

\begin{list}{}{\setlength{\leftmargin}{.4 cm}}

\item Choose $p \in A$ and $q \in C$. Then
\begin{eqnarray*}
& & \sum_{i=1}^m (\io \ot \om_{\la(b_i),\la(c_i)} \ot \io)((\de_r \ot \io)(\pi_r(p) \ot
q)(\pi_r(a_i) \ot 1 \ot 1)) \\
& & \spat = \sum_{i=1}^m (\io \ot \om_{\la(b_i),\la(c_i)} \ot \io)((\pi_r
\od \pi_r)(\de(p)(a_i \ot 1)) \ot q) \\
& & \spat = \sum_{i=1}^m (\pi_r \od \vfi)( (1 \ot c_i^*)\de(p)(a_i \ot b_i) ) \ot q \\
& & \spat = \sum_{i=1}^m (\pi_r \od \vfi)\bigl(\de(p)(a_i \ot b_i \rho(c_i^*))\bigr) \ot q \\
& & \spat = \sum_{j=1}^n (\pi_r \od \vfi)(\de(p \, r_j \rho(s_j^*))(t_j \ot 1)) \ot q \ .
\end{eqnarray*}
Therefore, the left invariance of $\vfi$ implies that
\begin{eqnarray*}
& & \sum_{i=1}^m (\io \ot \om_{\la(b_i),\la(c_i)} \ot \io)((\de_r \ot \io)(\pi_r(p) \ot
q)(\pi_r(a_i) \ot 1 \ot 1)) \\
& & \spat = \sum_{j=1}^n \pi_r(t_j) \, \vfi(p \, r_j \rho(s_j^*)) \ot q
= \sum_{j=1}^n \pi_r(t_j) \ot \vfi(s_j^* p \, r_j) \, q \\
& & \spat = \sum_{j=1}^n \pi_r(t_j) \ot \om_{\la(r_j),\la(s_j)} (\pi_r(p)) \, q
= \sum_{j=1}^n \pi_r(t_j) \ot (\om_{\la(r_j),\la(s_j)} \ot \io)(\pi_r(p) \ot  q) \ .
\end{eqnarray*}

\end{list}

Hence, using equation (*), the usual continuity arguments imply that
$$\sum_{i=1}^m \cU \, (\pi_r(a_i) \ot (\om_{\la(b_i),\la(c_i)} \ot \io)(\cU))
= \sum_{j=1}^n \pi_r(t_j) \ot (\om_{\la(r_j),\la(s_j)} \ot \io)(\cU) \ .$$
\end{demo}

This lemma implies immediately the following proposition (the second equality requires of
course a lemma similar to the one above).

\begin{proposition} \label{prop5.1}
We have that $\cU \, (\pi_r(A) \od B) = \pi_r(A) \od B$
and $(\pi_r(A) \od B) \, \cU = \pi_r(A) \od B$.
\end{proposition}

This justifies the following notation.

\begin{notation}
We define  $\tilde{\cU}$ as the unique element in $M(A \od B)$ such that $(\pi_r \od \io)(x)
\, \cU = (\pi_r \od \io)(x \, \tilde{\cU})$ and $\cU \, (\pi_r \od \io)(x) =
(\pi_r \od \io)(\tilde{\cU} \, x)$ for
every $x \in A \od B$.
\end{notation}

\begin{proposition}
The element $\tilde{\cU}$ is a non-degenerate corepresentation of
$(A,\de)$ on $B$.
\end{proposition}
\begin{demo}
We have for every $a \in A$, $x \in A \od A$  and $b \in B$ that
\begin{eqnarray*}
& & (\pi_r \od \pi_r \od \io)((\de \od \io)(\tilde{\cU}) (\de(a)\,x \ot b))
= (\pi_r \od \pi_r \od \io)((\de \od \io)(\tilde{\cU}(a \ot b))(x \ot 1)) \\
& & \spat = (\de_r \ot \io)\bigl((\pi_r \od \io)(\tilde{\cU}(a \ot b))\bigr) \,
((\pi_r \od \pi_r)(x) \ot 1)
= (\de_r \ot \io)(\cU \, (\pi_r(a) \ot b) )\, ((\pi_r \od \pi_r)(x) \ot 1) \\
& & \spat = (\de_r \ot \io)(\cU)\,(\de_r(\pi_r(a)) (\pi_r \od \pi_r)(x) \ot b)
= \cU_{13} \, \cU_{23} \, ((\pi_r \od \pi_r)(\de(a)\,x) \ot b) \\
& & \spat =  \cU_{13} \, (\pi_r \od \pi_r \od \io)(\tilde{\cU}_{23} (\de(a)\, x \ot b))
= (\pi_r \od \pi_r \od \io)(\tilde{\cU}_{13} \, \tilde{\cU}_{23} (\de(a)\,x \ot b)) \
\end{eqnarray*}
which implies that $(\de \od \io)(\tilde{\cU}) (\de(a)\,x \ot b)
= \tilde{\cU}_{13} \, \tilde{\cU}_{23} (\de(a)\,x \ot b)$.
So we get that $(\de \od \io)(\tilde{\cU}) = \tilde{\cU}_{13} \, \tilde{\cU}_{23}$.

Proposition \ref{prop5.1} of this section and  proposition 6.10 of
\cite{JK3} imply immediately that $\tilde{\cU}$ is non-degenerate.
\end{demo}

\bigskip

Next, we want to show that $B$ is a sub-$^*$-algebra of $M(C)$. For this, we will use the
next lemma.

\begin{lemma}
We have for every $a \in A$ and $b \in B$ that $$\cU^* (\pi_r(a) \ot
b) = (\pi_r \od \io)\bigl((S \od \io)((S^{-1}(a) \ot 1)
\,\tilde{\cU} (1 \ot b))\bigr) \ .$$
\end{lemma}
\begin{demo}
By proposition 5.11 of \cite{JK3} we know  that  $$\tilde{\cU} \, (S \od \io)( (S^{-1}(a)
\ot 1) \, \tilde{\cU} (1 \ot b)) = a \ot b \ . $$
Therefore $$\cU \, (\pi_r \od \io)\bigl((S \od \io)((S^{-1}(a) \ot 1) \, \tilde{\cU} (1 \ot
b))\bigr)
= (\pi_r \od \io)\bigl(\tilde{\cU} \, (S \od \io)((S^{-1}(a) \ot 1)\,\tilde{\cU}
(1 \ot b))\bigr) = \pi_r(a) \ot b  $$ which implies that
$$(\pi_r \od \io)\bigl((S \od \io)((S^{-1}(a) \ot 1)\,\tilde{\cU} (1 \ot b))\bigr)
= \cU^* (\pi_r(a) \ot b) \ .$$
\end{demo}

\begin{proposition}
The set $B$ is a sub-$^*$-algebra of $M(C)$ such that
$\overline{B C} = C$.
\end{proposition}
\begin{demo}
Because of the previous results, we only have to prove that $B$ is selfadjoint.

Choose $a,b,c \in A$. Take $y \in B$. Then
\begin{eqnarray*}
& & (\om_{\la(a),\la(c b)} \ot \io)(\cU)^* \, y
= (\om_{\la(c b),\la(a)} \ot \io)(\cU^*) \, y \\
& & \spat =  (\om_{\pi_r(c) \la(b),  \la(a)} \ot \io)(\cU^*) \, y
= (\om_{\la(b),\la(a)} \ot \io)(\cU^* (\pi_r(c) \ot y))  \\
& & \spat = (\om_{\la(b),\la(a)} \ot \io)\bigl(\,(\pi_r \od \io)
(\,(S \od \io)((S^{-1}(c) \ot 1) \, \tilde{\cU} (1 \ot y))\,)\,\bigr)
\hspace{1cm} \text{(*)}
\end{eqnarray*}
where we used the previous lemma in the last equality.

\begin{list}{}{\setlength{\leftmargin}{.4 cm}}

\item We have for every $p \in A$ and $q \in B$ that
\begin{eqnarray*}
& & (\om_{\la(b),\la(a)} \ot \io)\bigl((\pi_r \od \io)((S \od \io)(p \ot q))\bigr)
= \om_{\la(b),\la(a)}(\pi_r(S(p))) \,\, q \\
& & \spat = \vfi(a^* S(p) b) \, q = \vfi(S(S^{-1}(b)\,p\,S^{-1}(a^*))) \, q
= \vfi( S^{-1}(b) \,p\, S^{-1}(a^*) \, \sde) \, q \\
& & \spat = \om_{\la(S^{-1}(a^*)\sde),\la(S^{-1}(b)^*)}(\pi_r(p)) \, q
= (\om_{\la(S^{-1}(a^*)\sde),\la(S^{-1}(b)^*)} \ot \io)(\pi_r(p) \ot q) \ .
\end{eqnarray*}

\end{list}

Combining this with equality (*), we find that
\begin{eqnarray*}
& & (\om_{\la(a),\la(c b)} \ot \io)(\cU)^* \, y
=  (\om_{\la(S^{-1}(a^*)\sde),\la(S^{-1}(b)^*)} \ot \io)\bigl((\pi_r \od \io)
((S^{-1}(c) \ot 1)\,\tilde{\cU} (1 \ot y))\bigr) \\
& & \spat = (\om_{\la(S^{-1}(a^*)\sde),\la(S^{-1}(b)^*)} \ot
\io)((\pi_r(S^{-1}(c)) \ot 1) \,\cU (1 \ot y)) \\ & & \spat =
(\om_{\la(S^{-1}(a^*)\sde),\pi_r(S^{-1}(c)^*) \la(S^{-1}(b)^*)} \ot
\io)(\cU (1 \ot y)) \\ & & \spat =
(\om_{\la(S^{-1}(a^*)\sde),\la(S^{-1}(c b)^*)} \ot \io)(\cU) \, y \ .
\end{eqnarray*}
So we see that $(\om_{\la(a),\la(c b)} \ot \io)(\cU)^* =
(\om_{\la(S^{-1}(a^*)\sde),\la(S^{-1}(c b)^*)} \ot \io)(\cU)$ which
clearly belongs to $B$.

\medskip

Using the fat that $A^2 = A$, the proposition follows.
\end{demo}

So we have proven in effect that $\cU$ is of an algebraic nature. However, there is a
$^*$-algebra which is more natural than the $^*$-algebra $B$. We will introduce this
$^*$-algebra in the rest of this section.

\bigskip

By the previous results, we have the following proposition.

\begin{proposition}
The set $B C B$ is a dense sub-$^*$-algebra of $C$ such that $(\pi_r(A) \od B C B) \, \cU$
\newline $=\pi_r(A) \od B C B$ and $\cU \, (\pi_r(A) \od B C B) = \pi_r(A) \od B C B$.
\end{proposition}

\bigskip

The following object seems to be a natural candidate for our special $^*$-algebra connected
with $\cU$.

\begin{definition} \label{def5.1}
We define the set
$$C_\cU = \{ \, y \in C \mid \cU (\pi_r(A) \ot y) \subseteq \pi_r(A) \od C \text{ and }
(\pi_r(A) \ot y) \, \cU \subseteq \pi_r(A) \od C \,\} \ . $$
\end{definition}

It is easy to see that $C_\cU$ is a subalgebra of $C$. The previous proposition implies
that $B C B \subseteq C_\cU$ which implies that $C_\cU$ is dense in $C$. It is also easy
to see that $B C_\cU$ and $C_\cU B$ are subsets of $C_\cU$.

\begin{proposition}
We have that $\cU \, (\pi_r(A) \od C_\cU) \subseteq \pi_r(A) \od C_\cU$ and $(\pi_r(A) \od
C_\cU) \, \cU \subseteq \pi_r(A) \od C_\cU$.
\end{proposition}
\begin{demo}
Choose $a \in A$ and $x \in C_\cU$.

By definition, there exists $b_1,\ldots\!,b_n \in A$ and $y_1,\ldots\!,y_n \in C$ such that
$\cU(\pi_r(a) \ot x) = \sum_{i=1}^n \pi_r(b_i) \ot y_i$.

Using the Gramm-Schmidt orthonormalization procedure, we can find elements
$e_1,\ldots\!,e_m \in A$ such that $\vfi(e_j^* e_i) = \sde_{ij}$ for every $i,j \in
\{1,\ldots\!,m\}$ and such that $b_1,\ldots\!,b_n$ belong to $\langle e_1,\ldots\!,e_m
\rangle$.

Then there exist $z_1,\ldots\!,z_m \in C$ such that
$\cU(\pi_r(a) \ot x) = \sum_{i=1}^m \pi_r(e_i) \ot z_i$.

Fix $j \in \{1,\ldots\!,m\}$.
There exist $d \in A$ such that $e_i d = e_i$ for every $i \in \{1,\ldots\!,m\}$.

We have that
\begin{eqnarray*}
& & (\om_{\la(a d),\la(e_j)} \ot \io)(\cU) \, x
= (\om_{\pi_r(a) \la(d) , \la(e_j)} \ot \io)(\cU) \, x
= (\om_{\la(d) , \la(e_j)} \ot \io)(\cU(\pi_r(a) \ot x)) \\
& & \spat = \sum_{i=1}^m  (\om_{\la(d) , \la(e_j)} \ot \io)(\pi_r(e_i) \ot z_i)
= \sum_{i=1}^m \vfi(e_j^* e_i d) \, z_i
= \sum_{i=1}^m \vfi(e_j^* e_i) \, z_i = z_j \ .
\end{eqnarray*}
Hence we see that $z_j$ belongs to $B C_\cU$ which implies that $z_j$ belongs to $C_\cU$.

So we get that $\cU (\pi_r(a) \ot x)$ belongs to $\pi_r(A) \od C_\cU$.
Analogously, we get that $(\pi_r(a) \ot x)\,\cU$ belongs to $\pi_r(A) \od C_\cU$.
\end{demo}

\begin{proposition}
We have that $C_\cU = B C_\cU = C_\cU B = B C_\cU B$.
\end{proposition}
\begin{demo}
Choose $x \in C_\cU$. Take $a,b,c \in A$ such that $\vfi(c^* a b) =
1$. By the previous result, we know that $\cU(\pi_r(a) \ot x)$ belongs
to $\pi_r(A) \od C_\cU$. So there exist $b_1,\ldots\!,b_n \in A$ and
$y_1,\ldots\!,y_n \in C_\cU$ such that $\cU (\pi_r(a) \ot x)
= \sum_{i=1}^n \pi_r(b_i) \ot y_i$.

Therefore, we get  that $\pi_r(a) \ot x = \sum_{i=1}^n \cU^*(\pi_r(b_i) \ot y_i)$ which
implies that
\begin{eqnarray*}
x & = & \om_{\la(b),\la(c)}(\pi_r(a)) \, x
= \sum_{i=1}^n (\om_{\la(b),\la(c)} \ot \io)(\cU^* (\pi_r(b_i) \ot y_i))\\
& = &  \sum_{i=1}^n (\om_{\pi_r(b_i) \la(b) , \la(c)} \ot \io)(\cU^*(1 \ot y_i))
= \sum_{i=1}^n (\om_{\la(b_i b), \la(c)} \ot \io)(\cU^*) \, y_i  \\
& = &  \sum_{i=1}^n (\om_{\la(c),\la(b_i b)} \ot \io)(\cU)^* \, y_i
\ .
\end{eqnarray*}
Because $B$ is selfadjoint, this equation implies that $x$ belongs to $B C_\cU$.

So we have proven that $C_\cU \subseteq  B C_\cU$ which implies that
$C_\cU   = B C_\cU$. Similarly, we get that $C_\cU B= C_\cU$. From
this, we get immediately that $B C_\cU B = C_\cU$.
\end{demo}

\begin{corollary}
We have that $C_\cU = B C B$.
\end{corollary}

This follows immediately because $C_\cU = B C_\cU B \subseteq B C B \subseteq C_\cU$.

\bigskip

Consequently, we arrive at the following conclusion.

\begin{proposition}
The set $C_\cU$ is a dense sub-$^*$-algebra of $C$  such that $\cU \, (\pi_r(A) \od C_\cU) =
\pi_r(A) \od C_\cU$ and $(\pi_r(A) \od C_\cU) \, \cU = \pi_r(A) \od C_\cU$.
\end{proposition}

\medskip

\begin{remark} \rm
From the definition of $C_\cU$, it is also immediately clear that $C_\cU$ is the largest
subspace of $C$ such that $\cU (\pi_r(A) \od C_\cU) \subseteq \pi_r(A) \od C_\cU$ and
$(\pi_r(A) \od C_\cU) \, \cU \subseteq \pi_r(A) \od C_\cU$.
\end{remark}

Now we can give the following definition :

\begin{definition}
We define the element $\hat{\cU} \in M(A \od C_\cU)$ such that $(\pi_r \od \io)(x) \, \cU =
(\pi_r \od \io)(x \, \hat{\cU})$ and $\cU \, (\pi_r \od \io)(x) = (\pi_r \od \io)(\hat{\cU}
\, x)$ for every $x \in A \od C_\cU$.
Then $\hat{\cU}$ is a unitary corepresentation of $(A,\de)$ on $C_\cU$.
\end{definition}

\section{Lifting unitary corepresentations from the reduced to the universal \cst-algebra}
\label{art6}

In section \ref{art4}, we lifted certain $^*$-automorphisms on $A_r$
to $^*$-automorphisms on $A_u$. Using the results of the previous
section, we can now easily do a similar thing for unitary
corepresentations of $(A_r,\de_r)$.

\medskip

For the first part of this section, we fix a \cst-algebra $C$ and a unitary
corepresentation $\cU$ from $(A_r,\de_r)$ on $C$. This corepresentation will give rise to
a corepresentation on $(A_u,\de_u)$.

\medskip

We will use the notations $C_{\cU}$ and $\hat{\cU}$ from the previous
section. So $\hat{\cU}$ is an algebraic unitary corepresentation of
$(A,\de)$ on $C_{\cU}$. Then we can give immediately the following
definition.

\begin{definition} \label{def6.1}
We define $\cU_u$ as the unitary element in $M(A_u \ot C)$ such that $\cU_u \, (\pi_u \od
\io)(x) = (\pi_u \od \io)(\hat{\cU} \, x)$ and $(\pi_u \od
\io)(x) \, \cU_u = (\pi_u \od \io)(x \, \hat{\cU})$for every $x \in A \od C_\cU$.
\end{definition}

Then the following proposition is an easy consequence.

\begin{proposition} \label{prop6.2}
The element $\cU_u$ is a unitary corepresentation of $(A_u,\de_u)$ on
$C$ such that \newline $(\pi \ot \io)(\cU_u) = \cU$.
\end{proposition}

\bigskip

The following lemma guarantees that $\cU_u$ is uniquely determined by the above property.

\begin{lemma}  \label{lem6.1}
Consider a corepresentation $\cQ$ of $(A_u,\de_u)$ on $C$ such that $(\pi \ot \io)(\cQ) =
\cU$. Then we have that $\cQ_{13} = V_{12}^* \, \cU_{23} V_{12} \, \cU_{23}^*$.
\end{lemma}
\begin{demo}
Because $\cQ$ is a corepresentation of $(A_u,\de_u)$ on $C$, we have that $(\de_u \ot
\io)(\cQ) = \cQ_{13} \cQ_{23}$.

If we apply $\io \ot \pi \ot \io$ to this equation and use the fact that $(\pi \ot
\io)(\cQ) = \cU$, we get that $$((\io \ot \pi)\de_u \ot \io)(\cQ) = \cQ_{13} \,
\cU_{23} \ .$$

By proposition \ref{prop3.3}, we know that $(\io \ot \pi)\de_u(x) =
V^* (1 \ot \pi(x)) V$ for every $x \in A_u$. This implies that
$$((\io \ot \pi)\de_u \ot \io)(\cQ) = V_{12}^* (\pi  \ot \io)(\cQ)_{23} V_{12}
= V_{12}^* \, \cU_{23} V_{12} \ . $$
So we get that $V_{12}^* \,\cU_{23} V_{12} = \cQ_{13} \,\cU_{23}$.
\end{demo}

Using this lemma, we get the following uniqueness result.

\begin{proposition} \label{prop6.1}
The element $\cU_u$ is the unique corepresentation of $(A_u,\de_u)$ on
$C$ such that \newline $(\pi \ot \io)(\cU_u) = \cU$.
\end{proposition}

\bigskip

Reformulating lemma \ref{lem6.1} with $\cQ$ equal to $\cU_u$, the following
equality holds.

\begin{result}
We have that $(\cU_u)_{13} = V_{12}^* \, \cU_{23} V_{12} \, \cU_{23}^*$.
\end{result}

\begin{corollary} \label{cor6.1}
We have for every $\om \in B_0(H)$ that $((\io \ot \om)(V) \ot 1) \, \cU_u = (\io \ot
\om \ot \io)(\cU_{23} V_{12} \, \cU_{23}^*)$.
\end{corollary}

\medskip

By definition \ref{def6.1}, we have immediately the following result.

\begin{result} \label{res6.1}
We have that $\cU_u \, (\pi_u(A) \od C_\cU) =
\pi_u(A) \od C_\cU$ and $(\pi_u(A) \od C_\cU) \, \cU_u = \pi_u(A) \od C_\cU$.
\end{result}

Combining this with the fact that $(\pi \ot \io)(\cU_u) = \cU$ and definition
\ref{def5.1}, this implies the following result.

\begin{result} \label{res6.2}
We have that
$$C_\cU = \{ \, y \in C \mid \cU_u (\pi_u(A) \ot y) \subseteq \pi_u(A) \od C
 \text{ and } (\pi_u(A) \ot y) \, \cU_u  \subseteq \pi_u(A) \od C \,\} \ . $$
\end{result}

\bigskip\bigskip

Now we state the fact that there is a bijective correspondence between unitary
corepresentations of $(A_r,\de_r)$ and $(A_u,\de_u)$.

\begin{theorem} \label{thm6.1}
Consider a \cst-algebra $C$, then we have the following two properties.
\begin{enumerate}
\item Let $\cU$ be a unitary corepresentation of $(A_r,\de_r)$ on $C$.
Then $\cU_u$ is the unique unitary corepresentation of $(A_u,\de_u)$ on $C$ such that
$(\pi \ot \io)(\cU_u) = \cU$.
\item Let $\cU$ be a unitary corepresentation of $(A_u,\de_u)$ on $C$.
Then $(\pi \ot \io)(\cU)$ is the unique unitary corepresentation of $(A_r,\de_r)$ on $C$
such that $[(\pi \ot \io)(\cU)]_u = \cU$.
\end{enumerate}
\end{theorem}
\begin{demo}
\begin{enumerate}
\item This is just repeating proposition \ref{prop6.1}.
\item By proposition \ref{prop6.1}, we know that $[(\pi \ot \io)(\cU)]_u$ is the unique
unitary corepresentation of $(A_u,\de_u)$ on $C$ such that $(\pi \ot \io)\bigl([(\pi \ot
\io)(\cU)]_u\bigr)
= (\pi \ot \io)(\cU)$. So we must have that $[(\pi \ot \io)(\cU)]_u = \cU$.

If $\cV$ is a unitary corepresentation of $(A_r,\de_r)$ on $C$ such that $\cV_u = (\pi \ot
\io)(\cU)$, then proposition \ref{prop6.2} implies that
$\cV = (\pi \ot \io)(\cV_u) = (\pi \ot \io)(\cU)$.
\end{enumerate}
\end{demo}

\medskip

If we combine the second statement of this proposition with results \ref{res6.1} and
\ref{res6.2}, we get the followin proposition.

\begin{proposition} \label{prop6.3}
Consider a unitary corepresentation $\cU$ of $(A_u,\de_u)$ on a \cst-algebra $C$ and
define the set
$$C_\cU = \{ \, y \in C \mid \cU (\pi_u(A) \ot y) \subseteq \pi_u(A) \od C \text{ and }
(\pi_u(A) \ot y) \, \cU \subseteq \pi_u(A) \od C \,\} \ . $$
Then $C_\cU$ is a dense sub-$^*$-algebra of $C$ such that
$\cU \, (\pi_u(A) \od C_\cU) =
\pi_u(A) \od C_\cU$ and $(\pi_u(A) \od C_\cU) \, \cU = \pi_u(A) \od C_\cU$.
\end{proposition}

Of course, this proposition can also be proven in the same way as we did for unitary
corepresentations of $(A_r,\de_r)$ in the previous section.

\medskip

So we can give the following definition.

\begin{definition} \label{def6.2}
Consider a unitary corepresentation $\cU$ of $(A_u,\de_u)$ on a \cst-algebra $C$. Then we
define the element $\hat{\cU} \in M(A \od C_\cU)$ such that $(\pi_u \od
\io)(x) \, \cU = (\pi_u \od \io)(x \, \hat{\cU})$ and $\cU \, (\pi_u \od \io)(x)
= (\pi_u \od \io)(\hat{\cU} \, x)$ for every $x \in A \od C_\cU$.
Then $\hat{\cU}$ is a unitary corepresentation of $(A,\de)$ on $C_\cU$.
\end{definition}

\medskip

The first statement of the next proposition follows from result \ref{res6.2}, definition
\ref{def6.1} and the previous definition. Using theorem \ref{thm6.1}, the second one
follows from the first.

\begin{proposition}
Consider a \cst-algebra $C$, then we have the following properties.
\begin{enumerate}
\item Consider a unitary corepresentation $\cU$ of $(A_r,\de_r)$ on $C$.
Then we have that $C_{\cU_u} = C_\cU$ and $(\cU_u)\hoed = \hat{\cU}$
\item Consider a unitary corepresentation $\cU$ of $(A_u,\de_u)$ on $C$.
Then we have that $C_{(\pi \ot \io)(\cU)} = C_\cU$ and $\bigl((\pi \ot
\io)(\cU)\bigr)\hoed = \hat{\cU}$
\end{enumerate}
\end{proposition}

\section{The modular group of the left Haar weight on $(A_u,\de_u)$}

In a later section, we will introduce the left Haar weight on $(A_u,\de_u)$ using the left
Haar weight on $(A_r,\de_r)$. We want this left Haar weight on $A_u$ to be a KMS-weight
with respect to some one-parameter group. In this section, we will introduce this
one-parameter group using the techniques of section \ref{art4}.

\medskip

Remember from section \ref{art1}, that we have the left Haar weight
$\vfi_r$ on $(A_r,\de_r)$. This weight is a KMS-weight and the modular
group of $\vfi_r$ was denoted by $\si_r$. For some more notations, we
refer to the remarks after proposition \ref{prop1.1}.

Then we have that $\vfi_r$ is invariant under $\si_r$ and that $\la_r((\si_r)_t(a)) =
\nab^{it} \la_r(a)$ for every $a \in \Nfir$.

We  know that $\la(a)$ is analytic with respect to $\nab$ and that
$\nab^n \la(a) = \la(\rho^n(a))$ for every $a \in A$ and $n \in \Z$.

\bigskip

By definition 3.14 of \cite{Kus}, there exist a unique norm continuous one-parameter group
$K_r$ on $A_r$ such that
$$((\si_r)_t \ot (K_r)_t)\de_r = \de_r \, (\si_r)_t$$ for every $t \in \R$.
Notice that in \cite{Kus}, $K_r$ was denoted by $K$.  Proposition \ref{prop4.2} implies
that $\vfi_r$ is invariant with respect to $K_r$.

\medskip

In section 3 of \cite{Kus}, we introduced a positive injective
operator $P$ in $H$ which implements $K_r$, i.e. $(K_r)_t(x) = P^{it}
x P^{-it}$ for every $x \in A_r$ and $t \in \R$.

\begin{result}
We have for every $t \in \R$  and $a \in \Nfir$ that $P^{it} \la_r(a) = \la_r((K_r)_t(a))$.
\end{result}
\begin{demo}
Fix $t \in \R$ and define the unitary operator $u$ on $H$ such that $u \la_r(a) =
\la_r((K_r)_t(a))$ for every $a \in \Nfir$.

Then proposition \ref{prop4.1} implies that $(\nab^{it} \ot \nab^{it})
W = W (\nab^{it} \ot u)$. At the same time, proposition 3.12 of
\cite{Kus} implies that $(\nab^{it} \ot \nab^{it}) W = W (\nab^{it}
\ot P^{it})$. Comparing these two results, we see that $W(1 \ot
P^{it}) = W(1 \ot u)$. Consequently, $u=P^{it}$.
\end{demo}

\begin{result} \label{res7.1}
Consider $a \in A$ and $n \in \Z$. Then $\la(a)$ belongs to $D(P^n)$ and $P^n \la(a) =
\la(\sde^{-n} S^{-2n}(a) \sde^n)$.
\end{result}

The case $n=0$ is trivially true. The case $n=1$ follows from the
remark after definition 3.6 of \cite{Kus}. The case $n=-1$ follows
from the case $n=1$. The result follows now by induction.

\begin{corollary}
Consider $a \in A$. Then $\la(a)$ is analytic with respect to $P$.
\end{corollary}

\bigskip

By the results of section \ref{art4}, we have for every $t \in \R$ the
$^*$-automorphisms $(\si_t)_u$ and $(K_t)_u$ on $A_u$ such that
\begin{itemize}
\item $(\si_t)_u((\io \ot \om)(V)) = (\io \ot P^{it} \om \nab^{-it})(V)$
\item $(K_t)_u((\io \ot \om)(V)) = (\io \ot P^{it} \om P^{-it})(V)$
\end{itemize}
for every $\om \in B_0(H)^*$.

\bigskip

From the previous formulas, we get the following results :
\begin{enumerate}
\item We have for every $a \in A_u$ that the mapping
$\R \rightarrow A_u : t \mapsto (\si_t)_u(a)$ is norm continuous.
\item We have for every $s,t \in \R$ that $(\si_s)_u \, (\si_t)_u = (\si_{s+t})_u$.
\end{enumerate}

So we have the following definition.

\begin{definition}  \label{def7.1}
We define the norm continuous one-parameter group $\si_u$ on $A_u$ such that $(\si_u)_t =
(\si_t)_u$ for every $t \in \R$.
\end{definition}

Hence, $\si_u$ is determined by the property that
$$(\si_u)_t((\io \ot \om)(V)) = (\io \ot P^{it} \om \nab^{-it})(V)$$
for every $\om \in B_0(H)^*$ and $t \in \R$.

\medskip

Result \ref{res4.6} implies that $\pi \, (\si_u)_t = (\si_r)_t \, \pi$
for every $t \in \R$.

\bigskip

In a similar way, we can give the following definition.

\begin{definition}
We define the norm continuous one-parameter group $K_u$ on $A_u$ such that $(K_u)_t =
(K_t)_u$ for every $t \in \R$.
\end{definition}

Again, $K_u$ is determined by the property that
$$(K_u)_t((\io \ot \om)(V)) = (\io \ot P^{it} \om P^{-it})(V)$$
for every $\om \in B_0(H)^*$.

\medskip

We have also that $\pi \, (K_u)_t = (K_r)_t \, \pi$ for every $t \in \R$.

\bigskip

Proposition \ref{prop4.3} implies the following commutation relation.

\begin{result} \label{res7.2}
We have for every $t \in \R$ that $((\si_u)_t \ot (K_u)_t)\de_u = \de_u \, (\si_u)_t$.
\end{result}

Later on, we shall show that $\si_u$ and $K_u$ are determined by the
properties above (see proposition \ref{prop10.1}).

\bigskip\bigskip

In the rest of this section, we want to show that every element of $\pi_u(A)$ is analytic
with respect to $\si_u$.

\medskip

\begin{lemma}
Let $b,c$ be elements in $A$. Then $\rho^{-1}\bigl((\io \od \vfi)(\de(c^*)
(1 \ot b))\bigr) = (\io \od \vfi)(\de(\rho(c)^*)(1 \ot \sde S^2(b) \sde^{-1}))$.
\end{lemma}
\begin{demo}
By section \ref{alg}, we have for every $x \in A$ that $\rho^{-1}(x) =
\sde \, (\rho')^{-1}(x) \, \sde^{-1}$. This implies that
\begin{eqnarray*}
& & (\rho^{-1} \od \io)(\de(c^*)(1 \ot b)) = (\sde \ot 1) \,
((\rho')^{-1} \od \io)(\de(c^*)) \, (\sde^{-1} \ot b) \\ & & \spat =
(\sde \ot 1)\, (\io \od S^{-2})\bigl(\de((\rho')^{-1}(c^*))\bigr)\,(\sde^{-1}
\ot b)
\end{eqnarray*}
where  we used the results of section \ref{alg} in the last equality.
Using the fact that $(\rho')^{-1}(c^*) = \sde^{-1} \rho^{-1}(c^*)
\sde$, this  implies that
\begin{eqnarray*}
& & (\rho^{-1} \od \io)(\de(c^*)(1 \ot b))
= (\sde \ot 1)\,(\io \od S^{-2})(\de(\sde^{-1} \rho^{-1}(c^*) \sde))
\,(\sde^{-1} \ot b) \\
& & \spat = (1 \ot \sde^{-1})\,(\io \od S^{-2})(\de(\rho(c)^*))\,(1 \ot \sde b)
= (1 \ot \sde^{-1})\,(\io \od S^{-2})(\de(\rho(c)^*)(1 \ot \sde S^2(b) \sde^{-1}))
\,(1 \ot \sde)
\end{eqnarray*}
where we used that $\de(\sde) = \sde \ot \sde$ and $S^{-2}(\sde) =
\sde$. By section \ref{alg}, we know that $\vfi(\sde^{-1} S^{-2}(x)
\sde) = \vfi(x)$ for every $x \in A$. Therefore,
\begin{eqnarray*}
& & \rho^{-1}\bigl((\io \od \vfi)(\de(b^*)(1 \ot c))\bigr)
=  (\io \od \vfi)\bigl((\rho^{-1} \od \io)(\de(b^*)(1 \ot c))\bigr) \\
& & \spat = (\io \od \vfi)((1 \ot \sde^{-1})\,(\io \od S^{-2})(\de(\rho(c)^*)
(1 \ot \sde S^2(b) \sde^{-1}))\,(1 \ot \sde)) \\
& & \spat = (\io \od \vfi)(\de(\rho(c)^*)(1 \ot \sde S^2(b) \sde^{-1})) \ .
\end{eqnarray*}
\end{demo}

\begin{proposition}  \label{prop7.1}
Consider $a \in A$. Then $\pi_u(a)$ is analytic with respect to $\si_u$ and
$(\si_u)_{ni}(\pi_u(a)) = \pi_u(\rho^{-n}(a))$ for every $n \in \Z$.
\end{proposition}
\begin{demo}
Choose $b,c \in A$ and put $d = (\io \od \vfi)(\de(c^*)(1 \ot b))$.
Then lemma \ref{lem3.1} implies that
$$\pi_u(d) =(\io \ot \om_{\la(b),\la(c)})(V) \ .$$
Therefore, the remarks after definition \ref{def7.1} implies for every $t \in \R$ that
$$(\si_u)_t(\pi_u(d)) = (\io \ot \om_{P^{it}\la(b),\nab^{it} \la(c)})(V) \hspace{1cm} (*)$$
Because $\la(b)$ is analytic with respect to $P$ and $\la(c)$ is analytic with respect to
$\nab$, this implies easily that $\pi_u(d)$ is analytic with respect to $\si_u$.

\medskip

By equation (*), we have immediately that
$$(\si_u)_i(\pi_u(d))
= (\io \ot \om_{P^{i i} \la(b) , \nab^{i \,\overline{i}} \la(c)})(V) =
(\io \ot \om_{P^{-1} \la(b),\nab \la(c)})(V) \ . $$
Therefore, the remarks of the beginning of the section and result \ref{res7.1} imply that
$$(\si_u)_i(\pi_u(d)) = (\io \ot \om_{\la(\sde S^2(b) \sde^{-1}),\la(\rho(c))})(V) \ .$$
Hence, using lemma \ref{lem3.1} once more, we get that
$$ (\si_u)_i(\pi_u(d)) =
\pi_u\bigl((\io \od \vfi)(\de(\rho(c)^*) (1 \ot \sde S^2(b) \sde^{-1}))\bigr) \ .$$
Therefore, the previous lemma implies that
$$(\si_u)_i(\pi_u(d)) = \pi_u\bigl(\,\rho^{-1}(\,(\io \od \vfi)(\de(c^*)
(1 \ot b))\,)\,\bigr)= \pi_u(\rho^{-1}(d)) \ . $$

\medskip

From this all, we conclude for every $x \in A$ that  $\pi_u(x)$ is analytic with respect to
$\si_u$ and that $(\si_u)_i(\pi_u(x))$ \newline $=  \pi_u(\rho^{-1}(x))$. The general result
can now be proven easily.
\end{demo}

\section{The polar decomposition of the antipode on $A_u$}

In section 5 of \cite{Kus}, we arrived at a polar decomposition of the
antipode on $(A_r,\de_r)$. This polar decomposition consists of a norm
continuous one-parameter group $\tau_r$ on $A_r$ (denoted by $\tau$ in
\cite{Kus}) and an involutive anti-$^*$-automorphism $R_r$ on $A_r$
(denoted by $R$ in \cite{Kus}).

In this section, we want to transform these objects into similar objects on $A_u$ and
arrive at a polar decomposition of the antipode on $A_u$.

\bigskip

By proposition 5.8 of \cite{Kus}, we have for every $t  \in \R$ that $((\tau_r)_t \ot
(\tau_r)_t)\de_r = \de_r \, (\tau_r)_t$. In corollary 7.3 of \cite{Kus}, we proved
moreover the existence of a unique strictly positive number $\nu$ such that $\vfi
\, (\tau_r)_t = \nu^t \, \vfi$ for every $t \in \R$.

\medskip

Define the positive injective operator $Q$ in $H$ such that $Q^{it} \la_r(a) =
\nu^{-\frac{t}{2}} \, \la_r((\tau_r)_t(a))$ for every $a \in \Nfir$.
So we have that $(\tau_r)_t(x) = Q^{it} x Q^{-it}$ for every $x \in A_r$ and $t \in \R$.

We can now use again the results of section \ref{art4}. So we have for every $t
\in \R$ a $^*$-automorphism $(\tau_t)_u$ on $A_u$ such that $(\tau_t)_u((\io
\ot \om)(V)) = (\io \ot Q^{it} \om Q^{-it})(V)$ for every $t \in \R$.

\medskip

In the same way as in the previous section, this justifies the following definition.

\begin{definition}
We define the norm continuous one-parameter group $\tau_u$ on $A_u$ such that $(\tau_u)_t =
(\tau_t)_u$ for every $t \in \R$.
\end{definition}

Hence, $\tau_u$ is determined by the property that
$$(\tau_u)_t((\io \ot \om)(V)) = (\io \ot Q^{it} \om Q^{-it})(V)$$
for every $\om \in B_0(H)^*$ and $t \in \R$.

\medskip

Again, result \ref{res4.6}  implies that $\pi \, (\tau_u)_t =
(\tau_r)_t \, \pi$ for every $t \in \R$.

\medskip

From proposition \ref{prop4.3}, we get the following  commutation relation.

\begin{result}
We have for every $t \in \R$ that $((\tau_u)_t \ot (\tau_u)_t)\de_u = \de_u \, (\tau_u)_t$.
\end{result}

\medskip

We would also like to obtain a formula like in equation (8) of section
5 of \cite{Kus} for $\tau_u$.

\begin{result} \label{res8.1}
We have for every $t \in \R$ and $\om \in B_0(H)^*$ that $(\tau_u)_t((\io \ot \om)(V)) =
(\io \ot \nab^{it} \om \nab^{-it})(V)$.
\end{result}
\begin{demo}
Fix $t \in \R$ and $\om \in B_0(H)^*$. By proposition \ref{prop4.1}, we have that $(Q^{it}
\ot Q^{it})W = W (Q^{it} \ot Q^{it})$.

Hence, using equation (8) of \cite{Kus} and the fact that $\tau_r$ is implemented by $Q$,
we get that
$$(\io \ot \nab^{it} \om \nab^{-it})(W)
= (\tau_r)_t((\io \ot \om)(W)) = Q^{it} (\io \ot \om)(W) Q^{-it}
= (\io \ot Q^{it} \om Q^{-it})(W) \ .$$

The previous equality implies for every $\th \in B_0(H)^*$ that
$$(\nab^{it} \om \nab^{-it})((\th \ot \io)(W))
=(Q^{it} \om Q^{-it})((\th \ot \io)(W)) \ .$$
This implies that $(\nab^{it} \om \nab^{-it})(x) = (Q^{it} \om Q^{-it})(x)$ for every $x
\in \ah_r$. As a consequence, we get that $(\nab^{it} \om \nab^{-it})(x) =
(Q^{it} \om Q^{-it})(x)$ for every $x \in M(\ah_r)$.

\medskip

Choose $\eta \in A_u^*$. Because $V$ is an element of $M(A_u \ot
\ah_r)$, we have that $(\eta \ot \io)(V)$ belongs to $M(\ah_r)$. So the first part
implies that
$$(\nab^{it} \om \nab^{-it})((\eta \ot \io)(V)) =
(Q^{it} \om Q^{-it})((\eta \ot \io)(V))\ . $$
Hence,
$$\eta((\io \ot \nab^{it} \om \nab^{-it})(V))
= \eta((\io \ot Q^{it} \om Q^{-it})(V)) \ .$$
So we see that
$$(\io \ot \nab^{it} \om \nab^{-it})(V)
= (\io \ot Q^{it} \om Q^{-it})(V) = (\tau_u)_t((\io \ot \om)(V)) \ . $$
\end{demo}

\medskip

Using the equalities $((\si_u)_t \ot \io)(V) = (1 \ot \nab^{-it}) V (1 \ot P^{it})$ and
$((\tau_u)_t \ot \io)(V) = (1 \ot \nab^{-it}) V (1 \ot \nab^{it})$,
the proof of the next result is the same as the proof of proposition \ref{prop4.3}.

\begin{result} \label{res8.2}
We have for every $t \in \R$ that $((\tau_u)_t \ot (\si_u)_t)\de_u = \de_u \, (\si_u)_t$.
\end{result}

\bigskip

In the next part of this section, we prove that every element of
$\pi_u(A)$ is analytic with respect to $\tau_u$. The idea behind the
proof is the same as that of the proof of proposition \ref{prop7.1}.

\begin{lemma}  \label{lem8.1}
Consider $a \in A$ and $n \in \Z$. Then $\la(a)$ belongs to $D(Q^n)$
and $Q^n \la(a) = \nu^\frac{n i}{2}  \la(S^{2n}(a))$.
\end{lemma}
\begin{demo}
Choose $b \in A$. Then $\pi_r(b)$ belongs to $\Nfir$. We also know
that $\pi_r(b)$ belongs to ${\cal D}((\tau_r)_{-i})$ and
$(\tau_r)_{-i}(\pi_r(b)) = \pi_r(S^2(b))$ (see proposition 5.5 of
\cite{Kus}). So we get that $(\tau_r)_{-i}(\pi_r(b))$ belongs to
$\Nfir$.

Because $\la_r((\tau_r)_t(x)) = \nu^{\frac{t}{2}} \, Q^{it} \la_r(x)$
for every $x \in \Nfir$, this implies that $\la_r(\pi_r(b))$ belongs
to $D(Q)$ and
$$\nu^{-\frac{i}{2}} \, Q \la_r(\pi_r(b)) = \la_r\bigl((\tau_r)_{-i}(\pi_r(b))\bigr)
= \la_r\bigl(\pi_r(S^2(b))\bigr)$$

Because $\la_r(\pi_r(c))=\la(c)$ for every $c \in A$, we see that $\la(b)$ belongs to $D(Q)$
and $Q \la(b) = \nu^\frac{i}{2} \, \la(S^2(b))$.

The result follows from this result by induction.
\end{demo}

\begin{corollary}
Consider $a \in A$. Then $\la(a)$ is analytic with respect to $Q$.
\end{corollary}

\bigskip

\begin{proposition}
Consider $a \in A$. Then $\pi_u(a)$ is analytic with respect to $\tau_u$ and $(\tau_u)_{ni}
(\pi_u(a)) = \pi_u(S^{-2n}(a))$ for every $n \in \Z$.
\end{proposition}
\begin{demo}
Choose $b,c \in A$ and put $d = (\io \od \vfi)(\de(c^*)(1 \ot b))$.
Again, lemma \ref{lem3.1} implies that $\pi_u(d) =$ \newline $(\io
\ot \om_{\la(b),\la(c)})(V)$. Therefore, we get for every $t \in \R$ that
$$(\tau_u)_t(\pi_u(d)) = (\io \ot \om_{Q^{it} \la(b), Q^{it} \la(c)})(V)
\ . \hspace{1cm} (*)$$
Because $\la(b),\la(c)$ are analytic with respect to $Q$, this implies immediately that
$\pi_u(d)$ is analytic with respect to $\tau_u$.

Furthermore, equality (*) implies also that
\begin{eqnarray*}
& & (\tau_u)_i(\pi_u(d))
= (\io \ot \om_{Q^{i i} \la(b) , Q^{i \, \overline{i}} \la(c)})(V)
= (\io \ot \om_{Q^{-1} \la(b) , Q \la(c)})(V) \\
& & \spat = \nu^{-i} \, (\io \ot \om_{\la(S^{-2}(b)), \la(S^2(c))})(V)
\end{eqnarray*}
where we used lemma \ref{lem8.1} in the last equality. Hence, lemma
\ref{lem3.1} implies that
\begin{eqnarray*}
(\tau_u)_i(\pi_u(d))
& = &  \nu^{-i} \, \pi_u\bigl(\,(\io \od \vfi)(\,\de(S^2(c)^*)
(1 \ot S^{-2}(b))\,)\,\bigr) \\
& = &  \nu^{-i} \, \pi_u\bigl(\,(\io \od \vfi)(\,\de(S^{-2}(c^*))
(1 \ot S^{-2}(b))\,)\,\bigr) \\
& = & \nu^{-i} \, \pi_u\bigl(\,(\io \ot \vfi)(\,(S^{-2} \od S^{-2})
(\de(c^*)(1 \ot b))\,)\,\bigr) \ .
\end{eqnarray*}
By corollary 8.19 of \cite{Kus}, we know that $\vfi S^{-2} = \mu^{-1} \, \vfi =
\nu^i \, \vfi$. This implies that
$$ (\tau_u)_i(\pi_u(d))
= \pi_u\bigl(\,S^{-2}(\,(\io \od \vfi)(\de(c^*)(1 \ot b))\,)\,\bigr)
= \pi_u(S^{-2}(d)) \  .$$
From this all, we conclude for all $x \in A$ that $\pi_u(x)$ is analytic with respect to
$\tau_u$ and that $(\tau_u)_i(\pi_u(x))= \pi_u(S^{-2}(x))$. The general result follows by
induction.
\end{demo}

\bigskip

In the following part of this section, we will also transform $R_r$ to a $^*$-automorphism
on $A_u$. The idea behind this procedure will be the same as the idea in section
\ref{art4}, but the details are somewhat different (for instance, $\vfi_r$ is not
relatively invariant under $R_r$).

\medskip

At the end of section 7 of \cite{Kus}, we introduced the anti-unitary
involution $\hat{J}$ on $H$. This is in fact the modular conjugation
of the right invariant weight on $(\ah_r,\deh_r)$. Proposition 7.15 of
\cite{Kus} guarantees that $(\hat{J} \ot J)W^*(\hat{J} \ot J) = W$.
This relation was then used to prove that $R_r(x) = \hat{J} x^*
\hat{J}$ for every $x \in A_r$.

\medskip

As before, the relation $(\hat{J} \ot J)W^*(\hat{J} \ot J)
= W$ will imply that $J \ah_r J = \ah_r$ and $J M(\ah_r) J = M(\ah_r)$.

Then we can use $J$ to implement the anti-unitary antipode $\hat{R}_r$
on $(\ah_r,\deh_r)$, i.e. $\hat{R}_r(x) = J x^* J$ for every $x \in
\ah_r$.
We will have of course also that $\hat{R}_r$ is an involutive
anti-$^*$-automorphism of $\ah_r$ and that $\flip(\hat{R}_r \ot
\hat{R}_r)\deh_r = \deh_r \, \hat{R}_r$.

\medskip

The proof of the following lemma is similar to the proof of lemma
\ref{lem4.1} .

\begin{lemma}  \label{lem8.2}
Consider $\om_1,\om_2 \in B_0(H)^*$ such that $(\io \ot \om_1)(V) =
(\io \ot \om_2)(V)$. Then $(\io \ot \om_1(J.^*J))(V)
= (\io \ot \om_2(J.^*J))(V)$.
\end{lemma}

As before, this justifies the following definition :

\begin{proposition}
There exists a unique linear mapping $F$ from $\pi_u(A)$ into $A_u$ such that $F((\io \ot
\om_{p,q})(V))$ $= (\io \ot \om_{J q , J p})(V)$ for every $p,q \in \la(A)$.
\end{proposition}

\begin{result}
The mapping $F$ is antimultiplicative.
\end{result}
\begin{demo}
Choose $p_1,q_1,p_2,q_2 \in \la(A)$. Define the mapping $\om \in B_0(H)^*$ such that
$\om(x) =$ \newline $\langle W (x \ot 1) W^* (p_1 \ot p_2) , q_1 \ot q_2 \rangle$ for
every $x \in B_0(H)$. As usual, we have that
$$(\io \ot \om_{p_1,q_1})(V) \, (\io \ot \om_{p_2,q_2})(V)
= (\io \ot \om)(V) \ .$$
Using lemma \ref{lem8.2}, this gives us that
$$F\bigl(\,(\io \ot \om_{p_1,q_1})(V) \, (\io \ot \om_{p_2,q_2})(V)\,\bigr)
= (\io \ot \om(J .^* J))(V) \ .$$

We have for every $x \in \ah_r$ that
\begin{eqnarray*}
& & \om(J x^* J) = \langle W (J x^* J \ot 1) W^* (p_1 \ot p_2) , q_1 \ot q_2
\rangle = \langle
W (\hat{R}_r(x) \ot 1) W^* (p_1 \ot p_2) , q_1 \ot q_2 \rangle \\
& & \spat = \langle \deh_r(\hat{R}_r(x)) (p_1 \ot p_2) , q_1 \ot q_2 \rangle
= \langle \flip\bigl((\hat{R}_r \ot \hat{R}_r)(\deh_r(x))\bigr)\, (p_1 \ot p_2) ,
q_1 \ot q_2 \rangle \\ & & \spat = \langle (\hat{R}_r \ot \hat{R}_r)(\deh_r(x))
(p_2 \ot p_1) , q_2 \ot q_1 \rangle
= \langle (J \ot J) \deh_r(x)^* (J \ot J) (p_2 \ot p_1) , q_2 \ot q_1 \rangle \\
& & \spat = \langle \deh_r(x) (J q_2 \ot J q_1) , J p_2 \ot J p_1 \rangle
=  \langle W (x \ot 1) W^* J q_2 \ot J q_1 , J p_2 \ot J p_1 \rangle \ .
\end{eqnarray*}
Because $V$ belongs to $M(A_u \ot \ah_r)$ this equality implies that
$$(\io \ot \om(J .^* J))(V) = (\io \ot \om_{J q_2, J p_2})(V) \,
(\io \ot \om_{J q_1, J p_1})(V) \ .$$
Hence,
\begin{eqnarray*}
F\bigl(\,(\io \ot \om_{p_1,q_1})(V) \, (\io \ot \om_{p_2,q_2})(V)\,\bigr)
& = & (\io \ot \om_{J q_2, J p_2})(V) \, (\io \ot \om_{J q_1, J p_1})(V) \\
& = & F((\io \ot \om_{p_2,q_2})(V)) \, F((\io \ot \om_{p_1,q_1})(V)) \ .
\end{eqnarray*}
The result follows by linearity.
\end{demo}

\begin{result}
The mapping $F$ is selfadjoint.
\end{result}
\begin{demo}
Choose $p,q \in \la(A)$. Then $(\io \ot \om_{p,q})(V)^* = (\io \ot
\om_{q,p})(V^*)$. Because $p \in D(\nab^{-\frac{1}{2}})$ and $q \in
D(\nab^\frac{1}{2})$, corollary \ref{cor3.2} implies that
$$(\io \ot \om_{p,q})(V)^* = (\io \ot \om_{J \nab^{-\frac{1}{2}} p ,
J \nab^\frac{1}{2} q})(V) \ . $$
Hence, lemma \ref{lem8.2} implies that
$$F((\io \ot \om_{p,q})(V)^*)
= (\io \ot \om_{\nab^\frac{1}{2} q , \nab^{-\frac{1}{2}} p})(V) \ . $$
We also have that $J q$ belongs to $D(\nab^{-\frac{1}{2}})$ and
$\nab^{-\frac{1}{2}} J q = J \nab^\frac{1}{2} q$ and that
$J p$ belongs to $D(\nab^\frac{1}{2})$ and
$\nab^\frac{1}{2} J p = J \nab^{-\frac{1}{2}} p$.
So we see that
\begin{eqnarray*}
& & F((\io \ot \om_{p,q})(V)^*)
= (\io \ot \om_{J \nab^{-\frac{1}{2}} J q , J \nab^\frac{1}{2} J p})(V) \\
& & \spat = (\io \ot \om_{J p, J q})(V^*)
= (\io \ot \om_{J q , J p})(V)^*
= F((\io \ot \om_{p , q})(V))^* \ .
\end{eqnarray*}
The result follows by linearity.
\end{demo}

Therefore, proposition \ref{prop2.3} justifies the following
definition :

\begin{definition}
We define $R_u$ to be the anti-$^*$-homomorphism from $A_u$ into $A_u$
such that \newline $R_u((\io \ot \om_{p,q})(V)) = (\io \ot \om_{J q,
Jp})(V)$ for every $p,q \in \la(A)$.
\end{definition}

\begin{corollary}  \label{cor8.1}
We have that $(R_u \ot \io)(V) = (\io \ot (J .^*J))(V) = (\io \ot \hat{R}_r)(V)$.
\end{corollary}

\begin{corollary} \label{cor8.2}
We have for every $\om \in B_0(H)^*$ that $R_u((\io \ot \om)(V)) = (\io \ot \om(J.^*J))(V)$.
\end{corollary}

\medskip

Using this corollary, it is easy to see that the following holds.

\begin{proposition}
The mapping $R_u$ is an involutive anti-$^*$-automorphism on $A_u$.
\end{proposition}

Thanks to corollary \ref{cor8.1}, the antimultiplicativity of
$\hat{R}_r$ gives us the following commutation relation.

\begin{result} \label{res8.4}
The equality $\flip(R_u \ot R_u)\de_u = \de_u \, R_u$ holds.
\end{result}
\begin{demo}
By corollary \ref{cor8.1} and proposition \ref{prop3.2}, we have that
\begin{eqnarray*}
& & (\de_u \, R_u \ot \io)(V) = (\de_u \ot \hat{R}_r)(V)
= (\io \ot \io \ot \hat{R}_r)((\de_u \ot \io)(V)) \\
& & \spat = (\io \ot \io \ot \hat{R}_r)(V_{13} V_{23})
= (\io \ot \io \ot \hat{R}_r)((\flip \ot \io)(V_{23} V_{13})) \\
& & \spat = (\flip \ot \io)((\io \ot \io \ot \hat{R}_r)(V_{23} V_{13})) \ .
\end{eqnarray*}
Because $\hat{R}_r$ is antimultiplicative, this equation implies that
\begin{eqnarray*}
& & (\de_u \, R_u \ot \io)(V)
= (\flip \ot \io)\bigl((\io \ot \io \ot \hat{R}_r)(V_{13})
\,(\io \ot \io \ot \hat{R}_r)(V_{23}))\bigr) \\
& & \spat = (\flip \ot \io)(\,(R_u \ot \io \ot \io)(V_{13})
\, (\io \ot R_u \ot \io)(V_{23})\,)
= (\flip \ot \io)((R_u \ot R_u \ot \io)(V_{13} V_{23})) \\
& & \spat = (\flip(R_u \ot R_u) \ot \io)((\de_u \ot \io)(V))
= (\flip(R_u \ot R_u)\de_u \ot \io)(V) \ .
\end{eqnarray*}
Consequently, we get for every $\om \in B_0(H)^*$ that
$\de_u\bigl(R_u((\io \ot \om)(V))\bigr)
= \flip\bigl(\,(R_u \ot R_u)(\,\de_u((\io \ot \om)(V))\,)\,\bigr)$. The result follows.
\end{demo}

\bigskip

It is not difficult to prove that $R_u$ and $\tau_u$ commute :

\begin{result} \label{res8.3}
We have for every $t \in \R$ that $(\tau_u)_t \, R_u = R_u \, (\tau_u)_t$.
\end{result}
\begin{demo}
Choose $p,q \in H$. We know that $\nab^{it}$ and $J$ commute. Therefore,
\begin{eqnarray*}
& & R_u\bigl((\tau_u)_t((\io \ot \om_{p,q})(V))\bigr)
= R_u((\io \ot \om_{\nab^{it} p , \nab^{it} q})(V)) \\
& & \spat = (\io \ot \om_{J \nab^{it} q , J \nab^{it} p})(V)
= (\io \ot \om_{\nab^{it} J q , \nab^{it} J p})(V) \\
& & \spat = (\tau_u)_t((\io \ot \om_{J q, J p})(V))
= (\tau_u)_t\bigl(R_u((\io \ot \om_{p,q})(V))\bigr) \ .
\end{eqnarray*}
The result follows.
\end{demo}

At the end of this section, we prove the polar decomposition of the antipode on $A_u$.

\begin{theorem}
We have for every $a \in A$ that $R_u\bigl((\tau_u)_{-\frac{i}{2}}(\pi_u(a))\bigr)
= \pi_u(S(a))$.
\end{theorem}
\begin{demo}
Choose $b,c \in A$ and put $d = (\io \od \vfi)(\de(c^*)(1 \ot b))$. As
usual, we have that $\pi_u(d) = (\io \ot \om_{\la(b),\la(c)})(V)$. By
result \ref{res8.1}, this implies that
\begin{eqnarray*}
& & (\tau_u)_{-\frac{i}{2}}(\pi_u(d))
= (\io \ot \om_{\nab^{i (-\frac{i}{2})} \la(b) ,
\nab^{i (\,\overline{-\frac{i}{2}}\,)} \la(c)})(V) \\
& & \spat = (\io \ot \om_{\nab^\frac{1}{2} \la(b) , \nab^{-\frac{1}{2}} \la(c)})(V) \\
\end{eqnarray*}
Hence, using  corollary \ref{cor8.2}, we get that
\begin{eqnarray*}
& & R_u\bigl((\tau_u)_{-\frac{i}{2}}(\pi_u(d))\bigr)
=  (\io \ot \om_{J \nab^{-\frac{1}{2}} \la(c) , J \nab^\frac{1}{2} \la(b)})(V) \\
& & \spat = (\io \ot \om_{T^* \la(c) , T \la(b)})(V)
= (\io \ot \om_{\la(\rho(c^*)),\la(b^*)})(V) \ ,
\end{eqnarray*}
where we used the remarks after proposition \ref{prop1.1} in the last
equality. Using lemma \ref{lem3.1} once again, we see that
$$ R_u\bigl((\tau_u)_{-\frac{i}{2}}(\pi_u(d))\bigr)
= \pi_u\bigl(\,(\io \od \vfi)(\,\de(b)(1 \ot \rho(c^*))\,)\,\bigr) $$
Therefore, equation \ref{eq1.3} implies that
$$  R_u\bigl((\tau_u)_{-\frac{i}{2}}(\pi_u(d))\bigr)
= \pi_u\bigl((\io \od \vfi)((1 \ot c^*)\de(b))\bigr)
= \pi_u\bigl(\,S(\,(\io \od \vfi)(\de(c^*)(1 \ot b))\,)\,\bigr) = \pi_u(S(d)) \ .$$

\medskip

The theorem follows by linearity.
\end{demo}

\section{The left and right Haar weight on $(A_u,\de_u)$}

In this section, we introduce the left and right Haar weight on $A_u$.
For this, we will use the left and right Haar weight on $A_r$ and the
bridge mapping $\pi$. Once we have our weights, the proofs of
properties about these weights are completely analogous as the proofs
of their reduced counterparts.

\begin{definition}
We define the weight $\vfi_u = \vfi_r  \, \pi$, then $\vfi_u$ is a
densely defined lower semi-continuous weight on $A_u$ such that
$\pi_u(A) \subseteq \Mfiu$ and $\vfi_u(\pi_u(a)) = \vfi(a)$ for every
$a \in A$.
\end{definition}

Although $\vfi_r$ is a faithful weight on $A_r$, the weight $\vfi_u$
does not have to be faithful. It will turn out that the quantum group
$(A_u,\de_u)$ satisfies the definition of Masuda, Nakagami \&
Woronowicz except for the faithfulness of the left Haar weight.

\begin{definition}   \label{def9.1}
We use the GNS-construction $(H,\la_r,\io)$ for $\vfi_r$ to define the
GNS-construction for $\vfi_u$.   We define the mapping $\la_u$ from
$\Nfiu$ into $H$ such that $\la_u(a) = \la_r(\pi(a))$ for every $a
\in \Nfiu$. Then $(H,\la_u,\pi)$ is a GNS-construction for $\vfi_u$
such that $\la_u(\pi_u(a)) = \la(a)$ for every $a \in A$.
\end{definition}

\medskip

\begin{proposition}
The weight $\vfi_u$ is a KMS-weight on  $A_u$ with modular group
$\si_u$.
\end{proposition}

This proposition follows immediately from the fact that $\pi \,
(\si_u)_t = (\si_r)_t \, \pi$ for every $t \in \R$ (remarks after
definition \ref{def7.1}).

Because $\vfi_u$ is not faithful, the modular group is not uniquely
determined. However, by imposing an extra condition, the modular group
will be uniquely determined (see proposition \ref{prop10.2}).

\bigskip

In notation \ref{not3.1},  we introduced the notation $\om_{v,w}'$ for
every $v,w\in H$. Then the  proofs of the next two essential results
are completely analogous as the proofs of proposition 6.9 and 6.10 of
\cite{Kus}.

\begin{proposition}
Let $N$ be a dense left ideal in $A_u$ such that $(\om_{v,v}' \ot
\io)\de_u(x)$ belongs to $N$ for every  $v \in H$ and $x \in N$. Then
$\pi_u(A) \subseteq N$.
\end{proposition}

\begin{proposition}
Let $N$ be a dense left ideal in $A_u$ such that $(\io \ot
\om_{v,v}')\de_u(x)$ belongs to $N$ for every  $v \in H$ and $x \in N$.
Then $\pi_u(A) \subseteq N$.
\end{proposition}

\bigskip

Then we can prove the rather important result that $\vfi_u$ is
completely determined by its values on $\pi_u(A)$ :

\begin{theorem}  \label{thm9.1}
The set $\pi_u(A)$ is a core for $\la_u$.
\end{theorem}

The proof of this theorem is completely the analogous as the proof of
theorem 6.12 of \cite{Kus}. This starts with lemma 6.4 of \cite{Kus}
and goes on to theorem 6.12 of \cite{Kus} itself.

\bigskip

Thanks to this result, we get the  following strong form of left
invariance (see the proof of theorem 6.13 of \cite{Kus}). For used
notations, we refer again to the appendix.

\begin{theorem}
Consider $x \in \Mfiu$. Then $\de_u(x)$ belongs to
$\overline{\cM}_{\io \ot \vfi_u}$ and $(\io \ot \vfi_u)\de(x) =
\vfi_u(x) \, 1$.
\end{theorem}

This theorem justifies to call $\vfi_u$ the left Haar weight of
$(A_u,\de_u)$.

\bigskip

Then we get also the following weak form of left invariance.

\begin{corollary}
Consider $x \in \Mfiu$ and $\om \in A_u^*$. Then $(\om \ot
\io)\de_u(x)$ belongs to $\Mfiu$ and \newline $\vfi_u((\om \ot \io)\de_u(x)) =
\vfi_u(x) \, \om(1)$.
\end{corollary}

\bigskip

Using the counit $\vep_u$, the next result is trivial. To prove the
analog of it in the reduced case, we had to do a lot more of work.

\begin{proposition}
Let $x$ be an element in $A_u^+$ such that $(\om \ot \io)\de_u(x)$
belongs to $\Mfiu^+$ for every \newline $\om \in (A_u)^*_+$. Then  $x$
belongs to $\Mfiu^+$.
\end{proposition}

Of course, the result on $A_r$ (see theorem 3.11 of \cite{JK4}) allows us to get a
stronger version of this proposition.

\begin{proposition}
Let $x$ be an element in $A_u^+$ such that $(\om_{v,v}' \ot
\io)\de_u(x)$ belongs to $\Mfiu^+$ for every $v \in H$.
Then  $x$ belongs to $\Mfiu^+$.
\end{proposition}

\bigskip

In the same way as we prove proposition 6.15 of \cite{Kus}, we can
prove the strong left invariance proposed in the definition of Masuda,
Nakagami \& Woronowicz.

\begin{proposition}
Consider $a,b \in \Nfiu$ and $\om \in A_u^*$ such that $\om R_u \,
(\tau_u)_{-\frac{i}{2}}$ is bounded and call $\th$ the unique element
in $A_u^*$ which extends $\om R_u \, (\tau_u)_{-\frac{i}{2}}$. Then
$b^* \, (\om \ot \io)\de_u(a)$ and $(\th \ot \io)(\de_u(b^*)) \,a$
belong to $\Mfiu$ and
$$\vfi_u(\,b^* \, (\om \ot \io)\de_u(a)\,)
= \vfi_u(\,(\th  \ot \io)(\de_u(b^*)) \,a \,) \ .$$
\end{proposition}

\vspace{1.5cm}

In definition 9.2 of \cite{Kus}, the right Haar weight $\psi_r$ on
$(A_r,\de_r)$ (which was denoted by $\psi$ there) was defined by using
the left Haar weight and the anti-unitary antipode $R_r$ : $\psi_r =
\vfi_r \, R_r$.

\medskip

We used this weight $\psi_r$ to define a non-zero positive right
invariant linear $\psi$ on the algebraic quantum group $(A,\de)$ (see
definition 9.7 of \cite{Kus}). We found that $\pi_r(A)$ is a subset of
$\Mpsir$ and defined $\psi$ in such a way that $\psi(a) =
\psi_r(\pi_r(a))$ for every $a
\in A$.

\medskip

\begin{definition}
We define the weight $\psi_u = \psi_r \, \pi$. So $\psi_u$ is a
densely defined lower semi-continuous weight on $A_u$ such that
$\psi_u(\pi_u(a)) = \psi(a)$ for every $a \in A$.
\end{definition}

Using the fact that $\psi_r = \vfi_r \, R$ and $\pi R_u = R_r  \,
\pi$, we get the following result :

\begin{result}
The equality $\psi_u = \vfi_u \, R_u$ holds.
\end{result}

This suggests the following definition :

\begin{definition}
We define the norm-continuous one-parameter group $\si_u'$ on $A_u$
such that $(\si_u')_t = R_u \, (\si_u)_{-t} \, R_u$ for every $t \in
\R$.
\end{definition}

\begin{corollary}
The weight $\psi_u$ is a KMS-weight on $A_u$ with modular group
$\si_u'$.
\end{corollary}

It is straightforward to check that $\pi \, (\si_u')_t = (\si_r')_t
\, \pi$ for every $t \in \R$.

\bigskip

Combining results \ref{res8.2} and \ref{res8.4}, we get immediately
the following one.

\begin{result} \label{res9.1}
We have for every $t \in \R$ that $((\si_u')_t \ot (\tau_u)_{-t})\de_u
= \de_u \, (\si_u')_t$
\end{result}
\bigskip

In the same way as for $\vfi_u$, we can prove the following important
result about $\psi_u$ :

\begin{theorem}
Let $(H_{\psi_u},\la_{\psi_u},\pi_{\psi_u})$ be a GNS-construction for
$\psi_u$. Then $\pi_u(A)$ is a core for $\la_{\psi_u}$.
\end{theorem}

\bigskip

As usual, the left invariance of $\vfi_u$ is transferred to the right
invariance of $\psi_u$ :

\begin{theorem}
Consider $x \in \Mpsiu$. Then $\de_u(x)$ belongs to
$\overline{\cM}_{\psi_u \ot \io}$ and $(\psi_u \ot \io)\de_u(x) =
\psi_u(x) \, 1$.
\end{theorem}

This theorem justifies to call $\psi_u$ the right Haar weight of
$(A_u,\de_u)$.

\medskip

\begin{corollary}
Consider $x \in \Mpsiu$ and $\om \in A_u^*$. Then $(\io \ot
\om)\de(x)$ belongs to $\Mpsiu$ and \newline $\psi_u((\io \ot \om)\de_u(x)) =
\psi_u(x) \, \om(1)$.
\end{corollary}

\bigskip

We will end this section with a natural formula for $V$ in terms of the GNS-construction
of $\vfi_u$. As before, we refer to the appendix for the used notations.

\begin{proposition} \label{prop9.1}
Consider $a \in A_u$ and $b \in \Nfiu$. Then $\de(b)(a \ot 1)$ belongs
to $\cN_{\io \ot \vfi_u}$ and \newline $V (\io \ot
\la_u)(\de(b)(a \ot 1))
= a \ot \la_u(b)$.
\end{proposition}
\begin{demo}
First take a sequence $(a_n)_{n=1}^\infty$ in $A$ such that
$(\pi_u(a_n))_{n=1}^\infty$ converges to $a$. Because $\pi_u(A)$ is a
core for $\la_u$, we get the existence of a sequence
$(b_n)_{n=1}^\infty$ in $A$ such that $(\pi_u(b_n))_{n=1}^\infty$
converges to $b$ and $\bigl(\la_u(\pi_u(b_n))\bigr)_{n=1}^\infty$
converges to $\la_u(b)$.

This implies immediately that the sequence
$$\bigl(\, \de_u(\pi_u(b_n))(\pi_u(a_n) \ot 1) \, \bigr)_{n=1}^\infty
\rightarrow \de_u(b)(a \ot 1) \ . $$

Using the definition of $V$, we get for every $n \in \N$ that
$\de_u(\pi_u(b_n))(\pi_u(a_n) \ot 1)
= (\pi_u \ot \pi_u)(\de(b_n)(a_n \ot 1)) \in \pi_u(A) \od \pi_u(A) \subseteq
D(\io \ot \la_u)$ and
\begin{eqnarray*}
& & (\io \ot \la_u)\bigl(\de_u(\pi_u(b_n))(\pi_u(a_n) \ot 1))
= (\io \ot \la_u)\bigl((\pi_u \ot \pi_u)(\de(b_n)(a_n \ot 1))\bigr) \\
& & \spat =  (\pi_u \od \la)(\de(b_n)(a_n \ot 1))
= V^* \, (\pi_u(a_n) \ot \la(b_n)) =  V^* \, \bigl(\pi_u(a_n) \ot \la_u(\pi_u(b_n))\bigr)
\end{eqnarray*}
This implies that the sequence
$$\bigl(\, (\io \ot \la_u)(\,\de_u(\pi_u(b_n))(\pi_u(a_n) \ot 1)\,) \, \bigr)_{n=1}^\infty
\rightarrow V^* (a \ot \la_u(b)) \ . $$
Hence, the closedness of $\io \ot \la_u$ implies that $\de_u(b)(a \ot
1)$ belongs to $D(\io \ot \la_u)$ and $$(\io \ot \la_u)(\de_u(b)(a
\ot 1)) = V^* \, (a \ot \la_u(b)) \ . $$
\end{demo}

\section{Invariance properties of bi-\cst-isomorphisms and group-like elements}
\label{art10}

In this section, we prove some relative invariance properties of
bi-\cst-isomorphism and group-like elements. The results are
essentially the same  as the results of section 7 of \cite{Kus}.
However, some proofs have to be altered because the left Haar weight
on $(A_u,\de_u)$ is not necessarily  faithful.

\bigskip

The proof of the following proposition is the same as the proof of
proposition 7.1 of \cite{Kus}.

\begin{proposition}
Consider $^*$-automorphisms $\al$ and $\be$ on $A_u$ such that
$(\be \ot \al)\de_u  = \de_u \, \al$. Then there exists a strictly positive
number $r$ such that $\vfi_u \, \al = r \, \vfi_u$.
\end{proposition}

As before (see lemma \ref{lem4.2}), the mapping $\be$ above will
automatically satisfy $(\be \ot \be)\de_u = \de_u \, \be$. Therefore,

\begin{corollary}
Consider $^*$-automorphisms $\al$ and $\be$ on $A_u$ such that $(\be
\ot \al)\de_u = \de_u \, \al$. Then there exist  strictly positive numbers $r$ and $s$
such that $\vfi_u \, \al = r \, \vfi_u$ and $\vfi_u \, \be = s \,
\vfi_u$.
\end{corollary}

In a similar way as in proposition \ref{prop4.2}, we get  the
following result :

\begin{proposition} \label{prop10.3}
Consider $^*$-automorphisms $\al$ and $\be$ on $A_u$ such that $(\al
\ot \be)\de_u = \de_u \, \al$. Then there exists a strictly positive number $r$
such that $\vfi_u \, \al = r \, \vfi_u$ and $\vfi_u \, \be = r \,
\vfi_u$.
\end{proposition}

\bigskip

We want to prove a version of proposition 7.2 of \cite{Kus} on the
level of $A_u$. In the proof of proposition 7.2 of \cite{Kus}, we use
that $\vfi_r$ is faithful whereas $\vfi_u$ does not have to be
faithful. Therefore, we have to go around it in another way.

\vspace{1.5cm}

In section \ref{art4}, we lifted certain $^*$-automorphisms from the
reduced to the universal level. In the next part of this section, we
will consider the reverse process and show that both processes are
each other inverse.

\bigskip

Therefore,  consider $^*$-automorphisms $\al,\be$ on $A_u$ such that
$(\al \ot \be)\de_u = \de_u \, \al$. By  proposition \ref{prop10.3},
there exists a strictly positive number $r$ such that $\vfi_u \, \al =
r \, \vfi_u$ and $\vfi_u \, \be = r \, \vfi_u$.

\medskip

Define unitary operators $u,v \in B(H)$ such that $u \la_u(a) =
r^{-\frac{1}{2}} \, \la_u(\al(a))$ and $v \la_u(a) = r^{-\frac{1}{2}}
\, \la_u(\be(a))$. Then $\pi(\al(a)) = u
\pi(a) u^*$ and $\pi(\be(a)) = v \pi(a) v^*$ for every $a \in A$.
This implies immediately that
$u A_r u^* = A_r$ and $v A_r v^* = A_r$.

\medskip

\begin{definition} \label{def10.1}
Define the $^*$-automorphisms $\al_r$ and $\be_r$ on $A_r$ such that
$\al_r(x) = u x u^*$ and $\be_r(x) = v x v^*$ for every $x \in A_r$.
\end{definition}

Then $\pi \al = \al_r \, \pi$ and $\pi \be = \be_r \, \pi$. This
implies immediately that $(\al_r \ot \be_r)\de_r = \de_r \, \al_r$. So
we can use the results of section \ref{art4}.

\medskip

Because $\pi \al = \al_r \, \pi$, $\pi \be = \be_r \, \pi$ and $\vfi_u
= \vfi_r \, \pi$, we get also that $\vfi_r \, \al_r = r \, \vfi_r$ and $\vfi_r
\, \be_r = r \, \vfi_r$. Using definition \ref{def9.1}, we see also  that $u \la_r(a)
=  r^{-\frac{1}{2}} \, \, \la_r(\al_r(a))$ and $v \la_r(a) = r^{-\frac{1}{2}} \,
\la_r(\be_r(a))$ for every $a \in \Nfir$.

\medskip

With the notations used above, we get the following lemma.

\begin{lemma}
We have that $(\al \ot \io)(V) = (1 \ot u^*) V (1 \ot v)$
and $(\be \ot \io)(V) = (1 \ot v^*) V (1 \ot v)$.
\end{lemma}
\begin{demo}
Choose $a,b \in A$ and $p,q \in \Nfiu$. Then
\begin{eqnarray*}
\langle (\al \ot \io)(V) \, (a \ot \la_u(p)) , b \ot \la_u(q) \rangle
& =  & \al\bigl( \langle V (\al^{-1}(a) \ot \la_u(p)) , \al^{-1}(b)
\ot \la_u(q) \rangle \bigr) \\
& = &  \al\bigl( \langle \al^{-1}(a) \ot \la_u(p) , V^* (\al^{-1}(b)
\ot \la_u(q)) \rangle \bigr) \ . \hspace{1cm} \text{(a)}
\end{eqnarray*}
Using proposition \ref{prop9.1}, we get that $\de_u(q)(\al^{-1}(b)
\ot 1)$ belongs to $\cN_{\io \ot \vfi_u}$ and $V^* (\al^{-1}(b)  \ot \la_u(q))
= (\io \ot \la_u)(\de_u(q)(\al^{-1}(b) \ot 1))$.

So we see that $(\al^{-1}(b^*) \ot 1) \de_u(q^*) (\al^{-1}(a) \ot p)$
belongs to $\cM_{\io \ot \vfi_u}$ and
\begin{eqnarray*}
(\io \ot \vfi_u)((\al^{-1}(b^*) \ot 1) \de_u(q^*) (\al^{-1}(a) \ot p))
& = & \langle \al^{-1}(a) \ot \la_u(p) , (\io \ot \la_u)(\de_u(q)
(\al^{-1}(b) \ot 1))\rangle \\
& = & \langle \al^{-1}(a) \ot \la_u(p) , V^* (\al^{-1}(b)  \ot \la_u(q)) \rangle \ .
\end{eqnarray*}
Because $\vfi_u \, \be = r \, \vfi_u$, this implies that $(\al \ot
\be)((\al^{-1}(b^*) \ot 1) \de_u(q^*) (\al^{-1}(a) \ot p))$ belongs to
$\cM_{\io \ot \vfi_u}$ and
\begin{eqnarray*}
& & (\io \ot \vfi_u)\bigl((\al \ot \be)((\al^{-1}(b^*) \ot 1) \de_u(q^*)
(\al^{-1}(a) \ot p))\bigr) \\
& & \spat =  r \,\, \al\bigl( (\io \ot \vfi_u)((\al^{-1}(b^*) \ot 1)
\de_u(q^*) (\al^{-1}(a) \ot p))\bigr) \\
& & \spat = r \, \, \al\bigl(\langle \al^{-1}(a) \ot \la_u(p) ,
V^* (\al^{-1}(b)  \ot \la_u(q)) \rangle \bigr) \\
& & \spat =  r \,\,  \langle (\al \ot \io)(V) \, (a \ot \la_u(p)) ,
b \ot \la_u(q) \rangle
\end{eqnarray*}
where we used equation (a) in the last equality.

\medskip

Therefore, the equality $(\al_u \ot \be_u) \de_u = \de_u \, \al_u$
implies that  $(b^* \ot 1)\de_u(\al(q)^*) (a \ot \be(p))$ belongs to
$\cM_{\io \ot \vfi_u}$ and
$$(\io \ot \vfi_u)((b^* \ot 1)\de_u(\al(q)^*) (a \ot \be(p))
=  r \,\, \langle (\al \ot \io)(V)\,(a \ot \la_u(p)) , b \ot \la_u(q)
\rangle \ . \hspace{1cm} (b)$$

\medskip

We have also that $\be(p)$ and $\al(q)$ belong to $\Nfiu$ and
$\la_u(\be(p)) = r^\frac{1}{2} \, v \la_u(p)$ and $\la_u(\al(q)) =
r^\frac{1}{2} \, u \la_u(q)$. Therefore, proposition \ref{prop9.1}
implies that $\de_u(\al(q))(b \ot 1)$ belongs to $\cN_{\io \ot
\vfi_u}$ and
$$(\io \ot \la_u)(\de_u(\al(q))(b \ot 1)) = V^* (b \ot \la_u(\al(q)))
= r^\frac{1}{2} \, V^* (1 \ot u) (b \ot \la_u(q)) \ .$$
This implies that
\begin{eqnarray*}
& & (\io \ot \vfi_u)((b^* \ot 1)\de_u(\al(q)^*) (a \ot \be(p))
=  \langle a \ot \la_u(\be(p)) ,
(\io \ot \la_u)(\de_u(\al(q))(b \ot 1))  \rangle \\
& & \spat = r \, \langle (1 \ot v) (a \ot \la_u(p)) ,
V^* (1 \ot u) (b \ot \la_u(q)) \rangle
= r \, \langle (1 \ot u^*) V (1 \ot v) (a \ot \la_u(p)) ,
b \ot \la_u(q) \rangle
\end{eqnarray*}
Comparing this with (b), we get that
$$\langle (\al \ot \io)(V) \, (a \ot \la_u(p)) , b \ot \la_u(q) \rangle
= \langle (1 \ot u^*) V (1 \ot v) \, (a \ot \la_u(p)) ,
b \ot \la_u(q) \rangle \ . $$
So we see that $(\al \ot \io)(V) = (1 \ot u^*) V (1 \ot v)$.
The other equality is proven in a similar way.
\end{demo}

\bigskip\bigskip

Using this lemma and referring to results \ref{res4.4} and
\ref{res4.5}, we get the following proposition.

\begin{proposition} \label{prop10.4}
Consider $^*$-automorphisms $\al$,$\be$ on $A_u$ such that
$(\al \ot \be)\de_u = \de_u \, \al$.
Then we have that $(\al_r)_u = \al$ and $(\be_r)_u = \be$.
\end{proposition}

\medskip

The next proposition is even easier to prove.

\begin{proposition}
Consider $^*$-automorphisms $\al$,$\be$ on $A_r$ such that
$(\al \ot \be)\de_r = \de_r \, \al$.
Then we have that $(\al_u)_r = \al$ and $(\be_u)_r = \be$.
\end{proposition}

This follows immediately from the fact that $(\al_u)_r \, \pi = \pi \,
\al_u = \al \pi$ (see result \ref{res4.6}) and similarly for $\be$.

\vspace{1cm}

We can now use the above results to prove some kind of uniqueness
result.

\begin{proposition} \label{prop10.1}
Consider $^*$-automorphisms $\al_1,\al_2,\be_1,\be_2$ on $A_u$ such
that $(\al_1 \ot \be_1)\de_u = \de_u \, \al_1$, \newline $(\al_2 \ot
\be_2)\de_u = \de_u \, \al_2$ and $\pi \al_1 = \pi \al_2$.
Then $\al_1 = \al_2$ and $\be_1 = \be_2$.
\end{proposition}
\begin{demo}
By the results above, we have $^*$-automorphisms $(\al_1)_r$,
$(\al_2)_r$, $(\be_1)_r$ and $(\be_2)_r$ on $A_r$ such that $(\al_1)_r
\, \pi = \pi \al_1$, $(\al_2)_r \, \pi = \pi \al_2$, $(\be_1)_r \, \pi
= \pi \be_1$  and $(\be_2)_r \, \pi = \pi \be_2$.

By assumption, we have that $(\al_1)_r \, \pi = \pi \al_1 =
\pi \al_2 = (\al_2)_r \, \pi$, which implies that $(\al_1)_r =
(\al_2)_r$.

So, proposition \ref{prop10.4} implies that $\al_1 = ((\al_1)_r)_u =
((\al_2)_r)_u
= \al_2$. Because
$$(\al_1 \ot \be_1) \de_u = \de_u \, \al_1 = \de_u  \al_2 = (\al_2 \ot \be_2)\de_u =
(\al_1 \ot \be_2)\de_u \ ,$$ we get also that $\be_1 = \be_2$.
\end{demo}

Using the fact that $\flip(R_u \ot R_u)\de_u = \de_u \, R_u$, we get
also the following result.

\begin{corollary} \label{cor10.1}
Consider $^*$-automorphisms $\al_1,\al_2,\be_1,\be_2$ on $A_u$ such
that $(\be_1 \ot \al_1)\de_u = \de_u \, \al_1$, \newline $(\be_2 \ot
\al_2)\de_u = \de_u \, \al_2$ and $\pi \al_1 = \pi \al_2$. Then $\al_1 = \al_2$
and $\be_1 = \be_2$.
\end{corollary}

\bigskip

We will now give a first application of these uniqueness results.

\begin{proposition}
Consider $^*$-automorphisms $\al$ and $\be$ on $A_u$ such that $(\be
\ot \al)\de_u = \al$. \newline Then $\al \, (\si_u)_t = (\si_u)_t \, \al$ ,
$\be \, (\si_u)_t = (\si_u)_t \, \be$ and $\be \, (\tau_u)_t = (\tau_u)_t \, \be$ for
every $t \in \R$.
\end{proposition}
\begin{demo}
Fix $t \in \R$. Then there exist $^*$-automorphisms $\al'$, $\be'$ on
$A_r$ such that $\al' \, \pi = \pi \, \al$ and $\be' \, \pi = \pi \,
\be$ (This follows from definition \ref{def10.1} and the use of the map $R_u$).
Then $(\be' \ot \al') \de_r = \de_r \, \al'$.

From proposition 7.2 of \cite{Kus}, we get that $\al' \, (\si_r)_t =
(\si_r)_t \, \al'$. So
$$\pi \, \al \, (\si_u)_t = \al' \, (\si_r)_t \, \pi = (\si_r)_t \, \al' \, \pi
 = \pi \, (\si_u)_t \, \al \ .$$
At the same time, we have that
$$(\be \, (\tau_u)_t \ot \al \, (\si_u)_t)\de_u = \de_u  \, \al \, (\si_u)_t \hspace{1.5cm}
\text{ and } \hspace{1.5cm} ((\tau_u)_t \, \be \ot (\si_u)_t \, \al)\de_u
= \de_u \, (\si_u)_t \, \al$$
Therefore, corollary \ref{cor10.1} implies that $\al \, (\si_u)_t = (\si_u)_t \, \al$ and
$\be \, (\tau_u)_t = (\tau_u)_t \, \be$.

\medskip

We have also that $(\be \ot \be)\de_u = \de_u \, \be$. So we get in a
similar way that $\be \, (\si_u)_t = (\si_u)_t \, \be$.
\end{demo}

\begin{proposition}
Consider a $^*$-automorphism $\be$ on $A_u$ such that $(\be \ot
\be)\de_u = \de_u \, \be$. Then $\be \, R_u = R_u \, \be$.
\end{proposition}
\begin{demo}
Again, we get the existence of a $^*$-automorphism $\be_r$ on $A_r$
such that $\be_r \, \pi = \pi \, \be$. Then $(\be_r \ot \be_r)\de_r =
\de_r \, \be_r$, which implies that $\be_r \, R_r
= R_r \, \be_r$ by proposition 7.4 of \cite{Kus}.
Hence $\pi \,  R_u \, \be \,  R_u = R_r \, \be_r \, R_r \, \pi = \be_r
\, \pi = \pi \, \be$.

Because we also have that $(R_u \, \be \, R_u \ot R_u \, \be \, R_u)\de_u  =
\de_u \, R_u \, \be \, R_u$, we get that $R_u \, \be \, R_u = \be$ by proposition
\ref{prop10.1}.
\end{demo}

\bigskip

We want also to prove a version of proposition 7.5 of \cite{Kus} on
the level of $A_u$. Again, we have to find a way around the (possible)
non-faithfulness of $\vfi_u$.

\medskip

The prove of the next result is completely analogous as the proof of
lemma 7.11 of \cite{Kus}.

\begin{result}
Consider an element $x \in M(A_u)$ such that $\de_u(x) = x \ot 1$ or
$\de_u(x) = 1 \ot x$. Then $x$ will belong to $\C \, 1$.
\end{result}

\medskip

\begin{proposition}
Consider a unitary element $v \in M(A_u)$ such that $\de_u(v) =  v \ot
v$. Then there exists a strictly positive number $\lambda$ such that
$(\si_u)_t(v) = \lambda^{it} \, v$ for every $t \in \R$. Moreover,
$(\tau_u)_t(v) = v$ for every $t \in \R$.
\end{proposition}
\begin{demo}
Choose $t \in \R$. Because $((\tau_u)_t \ot (\tau_u)_t)\de_u = \de_u
\, (\tau_u)_t$, we get that $\de_u((\tau_u)_t(v)) = (\tau_u)_t(v) \ot (\tau_u)_t(v)$.

Because $(\pi \ot \pi)\de_u  = \de_r  \, \pi$, we have that $\de_r(\pi(v)) = \pi(v) \ot
\pi(v)$. Hence, proposition 7.5 of \cite{Kus} implies that $(\tau_r)_t(\pi(v)) = \pi(v)$.
So we get that $\pi((\tau_u)_t(v)) = (\tau_r)_t(\pi(v)) = \pi(v)$. Hence, proposition
\ref{prop6.1} implies that $(\tau_u)_t(v) = v$.

So we get that
\begin{eqnarray*}
& & \de_u(v^* (\si_u)_t(v)) = \de_u(v)^* \de_u((\si_u)_t(v))
= (v^*  \ot v^*) \, ((\tau_u)_t \ot (\si_u)_t)(\de_u(v)) \\
& & \spat = (v^* \ot v^*) \, ((\tau_u)_t(v) \ot (\si_u)_t(v))
= (v^* \ot v^*) \, (v \ot (\si_u)_t(v)) = 1 \ot v^* (\si_u)_t(v) \ .
\end{eqnarray*}
Therefore, the previous result implies the existence of a complex number $\lambda_t$ such
that $v^* (\si_u)_t(v) = \lambda_t \, 1$, so $(\si_u)_t(v) = \lambda_t
\, v$.

\medskip

We get easily the existence of a strictly positive number $\lambda$ such that
$\lambda^{it} = \lambda_t$  for every $t \in \R$.
\end{demo}

\bigskip

The proof of the next result is completely analogous as the proof of proposition 7.7 of
\cite{Kus}.

\begin{proposition}
Consider a unitary element $v$ in $M(A_u)$ such that $\de_u(v) = v \ot
v$. Then there exists a unique unitary element $x \in M(A)$ such that
$\pi_u(x\,a) = v \, \pi_u(a)$ and $\pi_u(a \, x) = \pi_u(a) \, v$ for
every $a \in A$. We have moreover that $\de(x) = x \ot x$.
\end{proposition}

Using this proposition, the proof of the following result is also
similar to the proof of proposition 7.8 of \cite{Kus}.

\begin{proposition}
Consider a unitary element $v$ in $M(A_u)$ such that $\de_u(v) = v \ot
v$. Then $R(v) = v^*$.
\end{proposition}

\bigskip

Using these results about unitary  elements, the proofs of the
following two results can be copied from proposition 7.9 of \cite{Kus}
and 7.10 of \cite{Kus}.

\begin{proposition}   \label{prop10.5}
Let $\al$ be a strictly positive  element affiliated with $A_u$ such
that $\de_u(\al) = \al \ot \al$. Then there exists a unique strictly
positive number $\lambda$ such that $(\si_u)_t(\al) = \, \lambda^t
\al$ for every $t \in \R$. We have moreover that $(\tau_u)_t(\al) = \al$
for every $t \in \R$ and $R_u(\al) = \al^{-1}$.
\end{proposition}

\vspace{1.5cm}

At the end of this section, we want to  mention some interesting
consequences of the uniqueness results in this section : proposition
\ref{prop10.1} and corollary \ref{cor10.1}.

\medskip

The first one follows immediately from proposition \ref{prop10.1} .

\begin{proposition}
Consider $t \in \R$ and let $\al$ be a $^*$-automorphism on $A_u$ such
that $(\al \ot \al)\de_u = \de_u \, \al$ and $\pi \, \al = (\tau_r)_t
\, \pi$. Then $(\tau_u)_t = \al$.
\end{proposition}

\medskip

The next consequence deals with the uniqueness of the modular group of
$\vfi_u$ under some extra condition.

\begin{proposition} \label{prop10.2}
Consider a norm continuous one-parameter group $\al$ on $A_u$ such
that $\vfi_u$ is KMS with respect to $\al$ and such that for every $t
\in \R$, there exists a $^*$-automorphism $\be_t$ on $A_u$ such that
$(\be_t \ot \al_t)\de_u = \de_u \, \al_t$. Then we have that $\al =
\si_u$.
\end{proposition}

If $\vfi_u$ is KMS with respect to $\al$, we have that $\pi \, \al_t = (\si_r)_t \, \pi$
for every $t \in \R$ (see \cite{JK1}). We know also that $\pi \, (\si_u)_t = (\si_r)_t \,
\pi$ for every $t \in \R$. Then the proposition follows immediately from corollary
\ref{cor10.1} .

\section{The modular function of the quantum group $(A_u,\de_u)$}

Consider the modular function $\sde_r$ of the quantum group
$(A_r,\de_r)$ (see definition 8.2 of \cite{Kus}, where it was denoted
by $\sde$). Again, we want to transform this modular function of
$(A_r,\de_r)$ to a modular function of $(A_u,\de_u)$. For this, we
will use the results of section \ref{art6}.

\medskip

Let us fix $t \in \R$. Using proposition 8.6 of \cite{Kus}, we have
that $\sde_r^{it}$ is a unitary element in $M(A_r)$ such that
$$\de_r(\sde_r^{it}) = \sde_r^{it} \ot \sde_r^{it} \ .$$

By definition \ref{def6.1}   and proposition \ref{prop6.1} , we get a
unique unitary element $(\sde_r^{it})_u$ in $M(A_u)$ such that
\begin{equation}
\pi((\sde_r^{it})_u) = \sde_r^{it} \hspace{1.5cm}  \text{ and }
\hspace{1.5cm} \de_u((\sde_r^{it})_u) = (\sde_r^{it})_u \ot
(\sde_r^{it})_u   \label{equ1}
\end{equation}

Using  corollary \ref{cor6.1}, we have moreover for every $\om \in
B_0(H)$ and $t \in \R$ that
\begin{equation}
(\io \ot \om)(V) \, (\sde_r^{it})_u = (\io \ot \sde_r^{-it} \om
\sde_r^{it})(V) \label{equ2}
\end{equation}

\bigskip

From this, we get that as usual the following results :
\begin{enumerate}
\item We have for every $a \in A_u$ that the mapping $ \R \rightarrow
A_u : t \mapsto a \, (\sde_r^{it})_u $ is norm-continuous.
\item The mapping $\R
\rightarrow M(A_u) : t \mapsto (\sde_r^{it})_u$ is a group homomorphism.
\end{enumerate}

\medskip

Therefore, the following definition is justified.

\begin{definition}
We define $\sde_u$ to be the unique strictly positive element
affiliated with $A_u$ such that $\sde_u^{it} = (\sde_r^{it})_u$ for
every $t \in \R$.
\end{definition}

Then $\sde_u$ is determined by the following two properties (see
equation \ref{equ1} and the remarks before it) :

\begin{proposition}
We have that $\de_u(\sde_u) = \sde_u \ot \sde_u$ and $\pi(\sde_u) =
\sde_r$.
\end{proposition}

\bigskip

Also the following holds by equation \ref{equ2} :

\begin{result} \label{res11.1}
We have for every $t \in \R$ and $\om \in B_0(H)^*$ that $(\io
\ot \om)(V) \, \sde_u^{it} = (\io \ot \sde_r^{-it} \om \sde_r^{it})(V)$.
\end{result}

\begin{corollary}
We have that $V(\sde_u \ot \sde_r) = (1 \ot \sde_r) V$.
\end{corollary}

\bigskip

The next proposition is an immediate consequence of proposition
\ref{prop10.5}.

\begin{proposition}
We have for every $t \in \R$ that $(\tau_u)_t(\sde_u) = \sde_u$.
Furthermore, $R_u(\sde_u) = \sde_u^{-1}$.
\end{proposition}

\bigskip

In proposition 8.17 of \cite{Kus}, we found that $(\si_r)_t(\sde_r) =
\nu^{-t} \, \sde_r$. As usual, we want to obtain this equation on the
level of $A_u$.

\begin{proposition} \label{prop11.1}
We have for every $t \in \R$ that $(\si_u)_t(\sde_u) = \nu^{-t}
\, \sde_u$.
\end{proposition}
\begin{demo}
Choose $s \in \R$. Take $\om \in B_0(H)^*$. Using result \ref{res11.1}
and the formula after definition \ref{def7.1}, we get that
\begin{eqnarray*}
& & (\si_u)_t((\io \ot \om)(V)) \, (\si_u)_t(\sde_u^{is})
= (\si_u)_t((\io \ot \om)(V)\,\sde_u^{is}) \\
& & \spat = (\si_u)_t((\io \ot \sde_r^{-is} \om \sde_r^{is})(V))
= (\io \ot P^{it} \sde_r^{-is} \om \sde_r^{is} \nab^{-it})(V) \ .
\end{eqnarray*}
By proposition 8.17 of \cite{Kus}, we know that
$(\si_r)_t(\sde_r^{is}) = \nu^{-ist} \, \sde_r^{is}$, which implies
that $\nab^{it} \sde_r^{is} \nab^{-it} =  \nu^{-ist} \, \sde_r^{is}$.

By proposition 7.9 of \cite{Kus}, we have also that
$(\tau_r)_{-t}(\sde_r^{-is}) = \sde_r^{-is}$, which by corollary 9.21
of \cite{Kus} implies that $(K_r)_t(\sde_r^{-is}) =
\sde_r^{-it}\, (\tau_r)_{-t}(\sde_r^{-is}) \, \sde_r^{it} =
\sde_r^{-is}$. Remembering that $K_r$ is implemented by $P$, we get
that $P^{it} \sde_r^{-is} P^{-it} = \sde_r^{-is}$.

So we see that $$(\si_u)_t((\io \ot \om)(V)) \, (\si_u)_t(\sde_u^{is})
= \nu^{-ist} \, (\io \ot \sde_r^{-is} P^{it} \om \nab^{-it} \sde_r^{is})(V) \ .$$
Using result \ref{res11.1}
and the formula after definition \ref{def7.1} once more, this implies that
$$ (\si_u)_t((\io \ot \om)(V)) \, (\si_u)_t(\sde_u^{is})
= \nu^{-ist} \, (\io \ot P^{it} \om \nab^{-it})(V) \, \sde_u^{is}
= \nu^{-ist} \, (\si_u)_t((\io \ot \om)(V)) \, \sde_u^{is} \ .$$
It follows that $(\si_u)_t(\sde_u^{is}) = \nu^{-ist} \, \sde_u^{is}$.

\medskip

From this all, we get the desired equality.
\end{demo}

\begin{corollary}
We have for every $t \in \R$ that $(\si_u')_t(\sde_u) = \nu^{-t} \,
\sde_u$.
\end{corollary}

\bigskip

The uniqueness results (proposition \ref{prop10.1} and corollary
\ref{cor10.1}) will prove useful once again :

\begin{result}
We have for every $a \in A_u$ and $t \in \R$ that $(\si_u')_t(a) =
\sde_u^{it} \, (\si_u)_t(a) \, \sde_u^{-it}$ and \newline $(K_u)_t(a) =
\sde_u^{-it} \, (\tau_u)_{-t}(a) \, \sde_u^{it}$.

\end{result}
\begin{demo}
Take $t \in \R$ and define the $^*$-automorphisms $\al$ and $\be$ on
$A_u$ such that $\al(a) = \sde_u^{it} \, (\si_u)_t(a) \,
\sde_u^{-it}$ and $\be(a) = \sde_u^{it} \, (K_u)_t(a) \, \sde_u^{-it}$
for every $a \in A_u$.

By proposition 9.9 of \cite{Kus}, we have for every $a \in A_u$ that
$$\pi(\al(a)) = \pi(\sde_u^{it} \, (\si_u)_t(a) \ \sde_u^{-it}) = \sde_r^{it}
\,(\si_r)_t(\pi(a)) \, \sde_r^{-it} = (\si_r')_t(\pi(a)) = \pi((\si_u')_t(a)) \ .$$
Using result \ref{res7.2}, we have also for every $a \in A_u$ that
$$\de_u(\al(a)) = \de_u(\sde_u^{it} \, (\si_u)_t(a) \, \sde_u^{-it})
= (\sde_u^{it} \ot \sde_u^{it})\, ((\si_u)_t \ot (K_u)_t)(\de_u(a)) \,
(\sde_u^{-it} \ot \sde_u^{-it}) = (\al \ot \be)\de(a) \ .$$ On the other hand,
we know by result \ref{res9.1} also that
$$\de_u \, (\si_u')_t = ((\si_u')_t \ot (\tau_u)_{-t})\de_u \ . $$
Hence, proposition \ref{prop10.1} implies that $(\si_u')_t = \al$ and $(\tau_u)_{-t}=\be$.
The proposition follows.
\end{demo}

\bigskip

Again, $\psi_u$ is absolutely continuous with respect to $\vfi_u$ and
a Radon Nikodym derivative is given by $\sde_u$.

\medskip

First, we need some notations. Consider a \cst-algebra $B$ and a
KMS-weight $\Upsilon$ on $B$ with modular group $\Sigma$. Let $\al$ be
a strictly positive element affiliated with $C$ such that there exists
a strictly positive number $\lambda$ such that $\Sigma_t(\al) =
\lambda^t \, \al$ for all $t \in \R$.

Then it is possible to define the KMS-weight $\Upsilon_\al =
\Upsilon(\al^\frac{1}{2} \, . \, \al^\frac{1}{2})$. For a precise definition,
see \cite{JK1}.

\bigskip

\begin{proposition}
We have the equality $\psi_u = (\vfi_u)_{\sde_u}$.
\end{proposition}
\begin{demo}
Remember that $\psi_u = \psi_r \, \pi$, $\vfi_u = \vfi_r \, \pi$ and
$\pi(\sde_u) = \sde_r$.

By theorem 9.18 of \cite{Kus}, we  know that $\psi_r  = (\vfi_r)_{\sde_r}$. Hence,

$$ (\vfi_u)_{\sde_u} = (\vfi_r \, \pi)_{\sde_u} =
(\vfi_r)_{\pi(\sde_u)} \, \pi =(\vfi_r)_{\sde_r} \, \pi = \psi_r \,
\pi = \psi_u \ .$$
\end{demo}

Because $\vfi_u$ does not have to be faithful, the Radon Nikodym
derivative is not uniquely determined. But, as before, imposing an
extra condition on the Radon Nikodym derivative implies the
uniqueness.

\medskip

Remember from proposition \ref{prop10.5} that if $\al$ is a strictly positive element
affiliated with $A_r$ satisfying $\de_u(\al) = \al \ot \al$, there
exists a strictly positive number $\lambda$ such that $(\si_u)_t(\al)
= \lambda^t \, \al$ for every $t \in \R$ so that $(\vfi_u)_\al$ is defined.

\begin{proposition}
Let $\al$ be a strictly positive element affiliated with $A_u$ such
that $\de_u(\al) = \al \ot \al$ and $(\vfi_u)_\al = \psi_u$. Then
$\al$ will be equal to $\sde_u$.
\end{proposition}
\begin{demo}
Becaue $(\vfi_u)_{\sde_u} = \psi_u = (\vfi_u)_\al$, we get that
$\pi(\al) = \pi(\sde_u)$. Because moreover $\de_u(\al) = \al \ot \al$
and $\de_u(\sde_u) = \sde_u \ot \sde_u$, proposition \ref{prop6.1}
implies that $\al = \sde_u$.
\end{demo}

\bigskip

In the next part of this section, we will prove the universal variant
of proposition 8.12 of \cite{Kus}. Remember from definition 8.13 of \cite{Kus}
that $\sde^z$ is defined on the algebraic level for every complex number $z$.

\begin{proposition}
Consider $a \in A$ and $z \in \C$. Then $\pi_u(a)$ belongs to
${\cal D}(\sde_u^z)$ and $\sde_u^z \, \pi_u(a) = \pi_u(\sde^z \, a)$.
\end{proposition}
\begin{demo}
Choose $b,c \in A$ and put $d = (\io \od \vfi)(\de(c^*)(1 \ot b))$.
Then $\pi_u(d) = (\io \ot \om_{\la(b),\la(c)})(V)$.
By result \ref{res11.1}, we have for every $t \in \R$ that
$$\pi_u(d) \, \sde_u^{it} = (\io \ot \om_{\sde_r^{-it} \la(b), \sde_r^{-it}
\la(c)})(V) \ .$$

Choose $y \in \C$. By the previous equality, lemma 8.8 of \cite{Kus}
and definition 8.13 of \cite{Kus}, we get immediately that $\pi_u(d)
\, \sde_u^{iy}$ is bounded and
$$\overline{\pi_u(d) \, \sde_u^{iy}}
= (\io \ot \om_{\sde^{-iy} \la(b), \sde^{-i\, \overline{y}} \la(c)})(V)
= (\io \ot \om_{\la(\sde^{-iy} b), \la(\sde^{-i\, \overline{y}} c)})(V)$$
Using lemma \ref{lem3.1} once again, we get that
\begin{eqnarray*}
\overline{\pi_u(d) \, \sde_u^{iy}}
& = & \pi_u\bigl( (\io \od \vfi)(\de(c^* \sde^{i y})(1 \ot
\sde^{-iy}b))\bigr)
= \pi_u\bigl((\io \od \vfi)(\de(c^*) (\sde^{iy} \ot \sde^{iy})
(1 \ot \sde^{-iy} b))\bigr) \\
& = &  \pi_u((\io \od \vfi)(\de(c^*)(1 \ot b)) \, \sde^{iy})
= \pi_u(d \, \sde^{iy})
\end{eqnarray*}

This implies for every $p \in A$ and every $x \in \C$ that $\pi_u(a)
\, \sde_u^x$ is bounded and $\overline{\pi_u(p) \, \sde_u^x} = \pi_u(p \, \sde^x)$.

\medskip

So we find in particular that $\pi_u(a^*) \sde_u^{\overline{z}}$ is
bounded and $\overline{\pi_u(a^*) \sde_u^{\overline{z}}} = \pi_u(a^*
\sde^{\overline{z}})$. This implies that $\pi_u(a)$ belongs to
${\cal D}(\sde_u^z)$ and $$\sde_u^z \, \pi_u(a) =
\bigl(\,\overline{\pi_u(a^*)\,\sde_u^{\overline{z}}}\,\bigr)^*
= \pi_u(a^* \sde^{\overline{z}})^* = \pi_u(\sde^z \, a) \ .$$
\end{demo}

So we get in particular that $\sde_u^{it} \, \pi_u(A) \subseteq
\pi_u(A)$ for every $t \in \R$. Using this fact we can prove the
following result in the same way as proposition 8.14 of \cite{Kus}.

\begin{proposition}
Consider $z \in \C$, then the set $\pi_u(A)$ is a core for $\sde_u^z$.
\end{proposition}

\bigskip

Using the theory of regular \cst-valued weights (see \cite{JK2}), we get an analytic
version of the algebraic equality  $(\vfi \od \io)\de(a) = \vfi(a) \, \sde$ where $a$ is
an element of $A$. (for used notations, see section \ref{app3}).

\medskip

The proof of the next proposition is the same as the proof of lemma 8.15 of \cite{Kus}.

\begin{proposition}
Consider $a \in \Nfiu$ and $b \in {\cal D}(\sde_u^\frac{1}{2})$. Then $\de(a)(1 \ot b)$
belongs to $\cN_{\vfi \ot \io}$ and
$$ \langle (\la_u \ot \io)(\de(a)(1 \ot b)) , (\la_u \ot \io)(\de(a)(1 \ot b)) \rangle =
\vfi_u(a^* a) \, (\sde_u^\frac{1}{2} b)^* (\sde_u^\frac{1}{2} b) \ \ .$$
\end{proposition}

\medskip

\begin{theorem}
Let $y$ be an element in $\Mfiu^+$. Then $\de_u(y)$ belongs to
$\hat{\cM}_{\vfi_u \ot \io}$ and \newline $(\vfi_u \ot \io)\de_u(y) =
\vfi_u(y)
\, \sde_u$.
\end{theorem}
\begin{demo} We refer to section \ref{app3} for used notations.

By definitions \ref{defap1} and \ref{defap2}, the previous  proposition implies for every
$p \in \Nfiu$ that $\de_u(p)$ belongs to $\tilde{\cN}_{\vfi_u \ot \io}$, that
${\cal D}(\sde_u^\frac{1}{2})
\subseteq D\bigl((\la_u \ot \io)(\de_u(p))\bigr)$ and that
$$\langle (\la_u \ot \io)(\de_u(p)) \, q , (\la_u \ot \io)(\de_u(p))\, q \rangle =
\vfi_u(p^* p) \, (\sde_u^\frac{1}{2} q)^* (\sde_u^\frac{1}{2} q) \hspace{1cm} \text{(a)} $$
for every $q \in {\cal D}(\sde_u^\frac{1}{2})$.

\medskip\medskip

Fix $c \in \Nfiu$ for the rest of this proof. Using theorem
\ref{thm9.1}, we get the existence of a sequence $(c_n)_{n=1}^\infty$
in $A$ such that $(\pi_u(c_n))_{n=1}^\infty$ converges to $c$ and
$\bigl(\la_u(\pi_u(c_n))\bigr)_{n=1}^\infty$ converges to $\la_u(c)$.

\medskip

Take $b \in A$. By the previous part, we know already that $\pi_u(b)$
is an element of $D\bigl((\la_u \ot \io)(\de_u(c))\bigr)$.

\medskip

Take $x \in D\bigl((\la_u \ot \io)(\de_u(c))\bigr)$. Using definitions \ref{defap1} and
\ref{defap2} once more, we see that $\de_u(c)(1 \ot x)$ belongs to $\cN_{\vfi_u \ot
\io}$ and $(\la_u \ot \io)(\de_u(c)(1 \ot x)) = (\la_u \ot \io)(\de_u(c)) \, x$.

\medskip

Take $a \in A$. Then we have also that $(\pi_u \od
\pi_u)(\de(a)(1 \ot b))$ belongs to and $\cN_{\vfi_u \ot
\io}$ and
$$(\la_u \ot \io)\bigl((\pi_u \od \pi_u)(\de(a)(1 \ot b))\bigr)
= (\la_u \ot \io)\bigl(\de_u(\pi_u(a))(1 \ot \pi_u(b))\bigr)
= (\la_u \ot \io)\bigl(\de_u(\pi_u(a))\bigr) \, \pi_u(b) \ . $$

\medskip

By section \ref{alg}, there exist $e \in A$ such that $\de(a)(1 \ot b)
= \de(a)(1 \ot b) (e \ot 1)$.

We know that $\pi_u(e)$ belongs to ${\cal D}((\si_u)_{i})$ and that
$(\si_u)_i(\pi_u(e)) = \pi_u(\rho^{-1}(e))$.

This implies that
\begin{eqnarray*}
& & \langle (\la_u \ot \io)\bigl(\de_u(\pi_u(a))\bigr) \, \pi_u(b) ,
(\la_u \ot \io)(\de_u(c)) \, x \rangle \\
& & \spat = \langle (\la_u \ot \io)\bigl((\pi_u \od \pi_u)
(\de(a)(1 \ot b))\bigr) , (\la_u \ot
\io)(\de_u(c)(1 \ot x)) \rangle \\
& & \spat = (\vfi_u \ot \io)\bigl([\de_u(c)(1 \ot x)]^* \,
(\pi_u \od \pi_u)(\de(a)(1 \ot b)) \bigr)\\
& & \spat = (\vfi_u \ot \io)( [\de_u(c)(1 \ot x)]^* \,
(\pi_u \od \pi_u)(\de(a)(1 \ot b)) \,(\pi_u(e) \ot 1)) \\
& & \spat = (\vfi_u \ot \io)\bigl((\pi_u(\rho^{-1}(e)) \ot 1) \,
[\de_u(c)(1 \ot x)]^* \,
(\pi_u \od \pi_u)(\de(a)(1 \ot b))\bigr) \\
& & \spat = (\vfi_u \ot \io)\bigl( (\pi_u(\rho^{-1}(e)) \ot x^*) \,
\de_u(c^*) \, (\pi_u \od
\pi_u)(\de(a)(1 \ot b))\bigr) \ . \hspace{1cm} \text{(b)}
\end{eqnarray*}

Notice that $\pi_u(\rho^{-1}(e))^* \ot x$ and $(\pi_u \od
\pi_u)(\de(a)(1 \ot b))$ belong to $\cN_{\vfi_u \ot \io}$. This
implies that the sequence
$$\bigl(\,(\vfi_u \ot \io)(\,(\pi_u(\rho^{-1}(e)) \ot x^*) \,
\de_u(\pi_u(c_n)^*) \, (\pi_u \od \pi_u)(\de(a)(1 \ot b))\,) \,
\bigr)_{n=1}^\infty$$ converges to
$$(\vfi_u \ot \io)(\,(\pi_u(\rho^{-1}(e)) \ot x^*) \, \de_u(c^*) \,
(\pi_u \od \pi_u)(\de(a)(1 \ot b))\,) \ . \hspace{1cm} \text{(c)}$$
We have also for every $n \in \N$ that $(\pi_u \od \pi_u)((\rho^{-1}(e)
\ot 1) \de(c_n^* a)(1 \ot b))$ belongs to $\cM_{\vfi_u \ot \io}$ which by result
\ref{resap1} implies that
\begin{eqnarray*}
& & (\vfi_u \ot \io)(\,(\pi_u(\rho^{-1}(e)) \ot x^*) \,
\de_u(\pi_u(c_n)^*) \, (\pi_u \od \pi_u)(\de(a)(1 \ot b))\,) \\
& & \spat = (\vfi_u \ot \io)(\,(1 \ot x^*)(\pi_u \od
\pi_u)((\rho^{-1}(e) \ot 1) \de(c_n^* a)(1 \ot b))\,) \\
& & \spat = x^* \,\, (\vfi_u \ot \io)(\,(\pi_u  \od \pi_u)(
(\rho^{-1}(e) \ot 1) \de(c_n^* a)(1 \ot b))\,) \\ & &
\spat = x^* \,\, (\vfi \od \pi_u)((\rho^{-1}(e) \ot 1) \de(c_n^* a)(1 \ot b) )
= x^* \,\, (\vfi \od \pi_u)(  \de(c_n^*a) (1 \ot b)(e \ot 1) ) \\
& & \spat = x^* \,\, (\vfi \od \pi_u)(  \de(c_n^*a) (1 \ot b))
= \vfi(c_n^* a) \,\,  x^*  \pi_u(\sde \, b)
= \langle \la_u(\pi_u(a)) , \la_u(\pi_u(c_n)) \rangle \,\, x^*  \pi_u(\sde \, b)
\end{eqnarray*}
were the algebraic version was used in the second last equality. This
implies that the sequence
$$\bigl(\,(\vfi_u \ot \io)(\,(\pi_u(\rho^{-1}(e)) \ot x^*) \, \de_u(\pi_u(c_n)^*) \,
(\pi_u \od \pi_u)(\de(a)(1 \ot b))\,) \, \bigr)_{n=1}^\infty$$ converges to $\langle
\la_u(\pi_u(a)) , \la_u(c) \rangle \,\, x^*  \pi_u(\sde \, b) \ .$ Combining this with
(c), we see that
$$(\vfi_u \ot \io)\bigl( (\pi_u(\rho^{-1}(e)) \ot x^*) \, \de_u(c^*) \,
(\pi_u \od \pi_u)(\de(a)(1 \ot b)) \bigr)
= \langle \la_u(\pi_u(a)) , \la_u(c) \rangle \,\, x^*  \pi_u(\sde \, b) \ .$$
Hence, equality (b) implies that
$$ \langle (\la_u \ot \io)\bigl(\de_u(\pi_u(a))\bigr) \, \pi_u(d) ,
(\la_u \ot \io)(\de_u(c)) \, x \rangle
= \langle \la_u(\pi_u(a)) , \la_u(c) \rangle \,\, x^*  \pi_u(\sde \, b) \ .$$

\medskip\medskip

From this last equality, we get for every $n \in \N$ that
$$ \langle (\la_u \ot \io)\bigl(\de_u(\pi_u(c_n))\bigr) \, \pi_u(b) ,
(\la_u \ot \io)(\de_u(c)) \, x \rangle
= \langle \la_u(\pi_u(c_n)) , \la_u(c) \rangle \,\, x^*  \pi_u(\sde \, d) \ .$$

We have chosen the sequence $(c_n)_{n=1}^\infty$ in such a way  that
$\bigl(\la_u(\pi_u(c_n))\bigr)_{n=1}^\infty$ converges to $\la_u(c)$.

By equality (a), we have for every $n \in \N$ that
$$\| (\la_u \ot \io)\bigl(\de_u(\pi_u(c_n))\bigr) \, \pi_u(b)
- (\la_u \ot \io)\bigl(\de_u(c))\bigr) \, \pi_u(b) \|^2
= \|  \la_u(\pi_u(c_n)) - \la_u(c) \|^2 \,\, \|\sde_u^\frac{1}{2} \, \pi_u(b)\|^2 \ . $$
This implies  that the sequence
$ \bigl(\,(\la_u \ot \io)\bigl(\de_u(\pi_u(c_n))\bigr) \, \pi_u(d)\, \bigr)_{n=1}^\infty$
converges to $(\la_u \ot \io)(\de_u(c)) \, \pi_u(d)$.

So we can conclude that
$$ \langle (\la_u \ot \io)(\de_u(c)) \, \pi_u(b) , (\la_u \ot \io)(\de_u(c)) \, x \rangle
= \langle \la_u(c) , \la_u(c) \rangle \,\, x^*  \pi_u(\sde \, b) \ .$$

\medskip

This last equality implies immediately that $(\la_u \ot \io)(\de_u(c)) \, \pi_u(b)$
belongs to $D((\la_u \ot \io)(\de_u(c))^*)$ and that
$$(\la_u \ot \io)(\de_u(c))^*\bigl((\la_u \ot \io)(\de_u(c)) \, \pi_u(b)\bigr)
= \langle \la_u(c) , \la_u(c) \rangle \,\,  \pi_u(\sde \, b)
= \vfi_u(c^* c) \,\, \sde_u \, \pi_u(b) \ . $$

So we see that $\pi_u(b)$ belongs to $D\bigl((\la_u \ot \io)(\de_u(c))^*
(\la_u \ot \io)(\de_u(c))\bigr)$ and
$$(\la_u \ot \io)(\de_u(c))^* (\la_u \ot \io)(\de_u(c)) \, \pi_u(b)
= \vfi_u(c^* c) \,\, \sde_u \, \pi_u(b) \, \hspace{1cm} \text{(d)}$$

\medskip

Because $\sde_u$ is affiliated to $A_u$, we know that
$1 + \vfi_u(c^* c) \, \sde_u$ has dense range in $A_u$.
Because $\pi_u(A)$ is a core for $\sde_u$, this implies that $(1 + \vfi_u(c^* c)
\, \sde_u) \, \pi_u(A)$ is dense in $A_u$.

Hence, equation (d) implies that $\bigl(1 + (\la_u \ot \io)(\de_u(c))^* (\la_u \ot
\io)(\de_u(c))\bigr) \, \pi_u(A)$ is dense in $A_u$. By the remarks after definition
\ref{defap2} and definition \ref{defap3}, this implies that $(\la_u \ot
\io)(\de_u(c))$ is regular so $\de_u(c)$ belongs to $\hat{\cN}_{\vfi_u \ot \io}$

Because of equation (d) and the fact that $\pi_u(A)$ is a core of
$\sde_u$, we see that
$$\vfi_u(c^* c) \,\, \sde_u \subseteq (\la_u \ot
\io)(\de_u(c))^* (\la_u \ot \io)(\de_u(c)) \ .$$
But both are selfadjoint, so they must be equal. By definition \ref{defap4}, we get that
$$(\vfi_u \ot \io)(\de_u(c^* c)) = (\la_u \ot \io)(\de_u(c))^* (\la_u \ot \io)(\de_u(c))
= \vfi_u(c^* c) \,\, \sde_u \ . $$
\end{demo}

\section{The universal corepresentation of $(A_u,\de_u)$}
\label{art12}

Remember from section \ref{alg} that we constructed an (algebraic) universal
corepresentation $\aU \in M(A \od \ah)$ for $(A,\de)$ in \cite{JK3}. This element served as
a link between corepresentations of $(A,\de)$ and $^*$-homomorphisms on $\ah$.

\medskip

In this section, we will use this algebraic universal corepresentation to define a
\cst-algebraic universal corepresentation of $(A_u,\de_u)$.

\bigskip

It is clear that everything we did for $A$, we can also do for $\ah$. So we will get a
universal \cst-algebraic quantum group $(\ah_u,\deh_u)$ in the same way as we got the
universal \cst-algebraic quantum group $(A_u,\de_u)$.

We will denote the canonical embedding from $\ah$ into $\ah_u$ by $\pih_u$. We also define
$\pih$ as the unique $^*$-homomorphism from $\ah_u$ into $\ah_r$ such that $\pih \! \circ
\! \pih_u = \pih_r$.

\medskip

Remembering that $\aU$ is a unitary element in $M(A \od \ah)$, we can give give the
following definition.

\begin{definition}
We define the unitary element $\uU$ in $M(A_u \ot \ah_u)$ in such a way that
$$\uU \, (\pi_u \od \pih_u)(x) = (\pi_u \od \pih_u)(\aU \, x)  \hspace{1.5cm}
\text{ and } \hspace{1.5cm} (\pi_u \od \pih_u)(x) \, \uU = (\pi_u \od \pih_u)(x \, \aU)$$
for every $x \in A \od \ah$.
\end{definition}

\medskip

Notice that we also have that $(\de \od \io)(\aU) = \aU_{13} \, \aU_{23}$ and
$(\io \od \deh)(\aU) = \aU_{12} \, \aU_{23}$ which implies the following result :

\begin{proposition}
We have that $(\de \ot \io)(\uU) = \uU_{13} \, \uU_{23}$ and $(\io \ot \deh)(\uU)
= \uU_{12} \, \uU_{23}$.
\end{proposition}

\medskip

The definition of $\uU$ together with proposition \ref{prop1.3} and proposition
\ref{prop3.1} will also immediately imply that $(\pi \od \pih)(\uU) = W$ and $(\io \ot
\pih)(\uU) = V$.

\medskip

Consider $v,w \in H$.

In notation \ref{not3.1}, we defined $\om_{v,w}' \in A_u^*$ such
that $\om_{v,w}'(x) = \langle \pi(x) \, v , w \rangle$ for every $x \in A_u$.

In a similar way, we define $\hat{\om}_{v,w}' \in \ah_u^*$ such that $\hat{\om}_{v,w}'(x)
= \langle \pih(x) \, v , w \rangle$ for every $x \in \ah_u$.

\medskip

Then we have the following equalities :

\begin{result}
Consider $a,b \in A$, then the following two equalities hold :
\begin{enumerate}
\item $(\io  \ot  \hat{\om}_{\la(a),\la(b)}')(\uU) = \pi_u\bigl((\io \od \vfi)(\de(b^*)
(1 \ot a))\bigr)$
\item $(\om_{\la(a),\la(b)}' \ot \io)(\uU) = \pih_u(a \vfi b^*)$
\end{enumerate}
\end{result}
\begin{demo}
\begin{enumerate}
\item We get immediately that $\hat{\om}_{\la(a),\la(b)}' = \om_{\la(a),\la(b)} \, \pih$.
So, using the fact that $(\io \ot \pih)(\uU) = V$ and lemma \ref{lem3.1}, we get that
\begin{eqnarray*}
& & (\io  \ot  \hat{\om}_{\la(a),\la(b)}')(\uU)
= (\io \ot \om_{\la(a),\la(b)})((\io \ot \pih)(\uU)) \\
& & \spat = (\io \ot \om_{\la(a),\la(b)})(V) = \pi_u\bigl((\io \od \vfi)(\de(b^*)
(1 \ot a))\bigr) \ .
\end{eqnarray*}
\item Choose $c,d \in A$. Take $\th \in \ah$.
Then
\begin{eqnarray*}
(\om_{\la(c d),\la(b)}' \ot \io)(\uU) \, \pih_u(\th)
& = & (\om_{\pi_r(c) \la(d) , \la(b)}' \ot \io)(\uU) \, \pih_u(\th) \\
& = & (\om_{\la(d),\la(b)}' \ot \io)\bigl(\uU (\pi_u(c) \ot \pih_u(\th))\bigr) \ .
\end{eqnarray*}
Hence, the definition of $\uU$ implies that
\begin{eqnarray*}
(\om_{\la(c d),\la(b)}' \ot \io)(\uU) \, \pih_u(\th)
& = & (\om_{\la(d),\la(b)}' \ot \io)\bigl((\pi_u \od \pih_u)(\aU (c \ot \th))\bigr) \\
& = & \pi_u\bigl( (d \vfi b^* \od \io)(\aU (c \ot \th)) \bigr) \hspace{2.5cm} \text{(*)}
\end{eqnarray*}
By result 6.8 of \cite{JK3}, we know that $(c d \vfi b^* \ot \io)(\aU) = c d
\vfi b^*$. Looking at remark 4.4 of \cite{JK3}, this means that
$$ (d \vfi b^* \od \io)(\aU (c \ot \th)) =
 (c d \vfi b^* \od \io)(\aU) \, \th = ( c d \vfi b^* ) \, \th $$
Substituting this in equation (*), we get that
$$ (\om_{\la(c d),\la(b)}' \ot \io)(\uU) \, \pih_u(\th)
= \pi_u\bigl(( c d \vfi b^* ) \, \th\bigr) = \pi_u(c d \vfi b^*) \, \pi_u(\th) \ .$$
So we see that
$$ (\om_{\la{c d},\la(b)}' \ot \io)(\uU) = \pi_u(c d \vfi b^*) \ .$$
Now the equality follows because $A^2 = A$.
\end{enumerate}
\end{demo}

\medskip

\begin{corollary}
We have the following two properties :
\begin{enumerate}
\item We have that $A_u$ is a subset of the closure of $\{\, (\io \ot \om)(\uU)
\mid \om \in \ah_u^* \,\}$ in $M(A_u)$.
\item We have that $\ah_u$ is a subset of the closure of $\{\, (\om \ot \io)(\uU)
\mid \om \in A_u^* \,\}$ in $M(\ah_u)$.
\end{enumerate}
\end{corollary}

\bigskip

Now we can prove easily the universality property of $\uU$ (this kind of universality was
introduced for compact quantum groups in \cite{PW}).

\begin{theorem}
Consider a \cst-algebra $C$ and a unitary corepresentation $\cU$ of $(A_u,\de_u)$ on $C$.
Then there exists a unique non-degenerate$^*$-homomorphism $\th$ from $\ah_u$ into $M(C)$
such that $(\io \ot \th)(\uU) = \cU$.
\end{theorem}
\begin{demo}
Uniqueness follows immediately from the previous corollary. We turn to the existence.

\medskip

In proposition \ref{prop6.3}, we got a dense $^*$-algebra $C_\cU$ of $C$ such that $\cU \,
(\pi_u(A) \od C_\cU) = \pi_u(A) \od C_\cU$ and $(\pi_u(A) \od C_\cU) \, \cU = \pi_u(A) \od
C_\cU$.

We defined moreover in definition \ref{def6.2} the element $\hat{\cU} \in M(A \od C_\cU)$
such that $(\pi_u \od \io)(x) \, \cU = (\pi_u \od \io)(x \, \hat{\cU})$ and $\cU \, (\pi_u
\od \io)(x) = (\pi_u \od \io)(\hat{\cU} \, x)$ for every $x \in A \od C_\cU$.
Then $\hat{\cU}$ is a unitary corepresentation of $(A,\de)$ on $C_\cU$.

\medskip

In proposition 6.11 of \cite{JK3}, we proved the existence of an algebraically
non-degenerate $^*$-homomorphism $\eta$ from $\ah$ into $M(C_\cU)$ such that $(\io \od
\eta)(\aU) = \hat{\cU}$.

Combining lemma \ref{lem2.1} and the universality property of $\ah_u$, we get the
existence of a unique non-degenerate $^*$-homomorphism $\th$ from $\ah_u$ into $M(C)$ such
that $\th(\pi_u(\om)) \, c = \eta(\om) \, c$ and $c \, \th(\pi_u(\om)) = c \, \eta(\om)$
for every $c \in C_\cU$ and $\om \in \ah$.

\medskip

Choose $a \in A$, $\om \in \ah$ and $c \in C$. Then
\begin{eqnarray*}
(\io \ot \th)(\uU) \, (\pi_u(a) \ot \th(\pi_u(\om))\, c \,)
& = & (\io \ot \th)\bigl(\uU (\pi_u(a) \ot  \pi_u(\om))\bigr) \, (1 \ot c) \\
& = & (\io \ot \th)\bigl((\pi_u \od \pih_u)(\aU (a \ot \om))\bigr) \, (1 \ot c)
\end{eqnarray*}
where we used the definition of $\uU$ in the last equality. So we get that
$$  (\io \ot \th)(\uU) \, (\pi_u(a) \ot \th(\pi_u(\om))\, c \,)
= (\pi_u \od \eta)( \aU (a  \ot \om)) \, (1 \ot c) $$
If we use the fact that $(\io \od \eta)(\aU) = \hat{\cU}$, this implies that
$$ (\io \ot \th)(\uU) \, (\pi_u(a) \ot \th(\pi_u(\om))\, c \,)
= (\pi_u \od \io)( \hat{\cU} (a \ot \eta(\om) \, c\,) )
= \cU \, (\pi_u(a) \ot \eta(\om) \, c\,) $$
where we used the definition of $\hat{\cU}$ in the last equality. Using the fact that
$\eta(\ah) \, C_\cU = C_\cU$ (which is dense in $C$), we get that
$$ (\io \ot \th)(\uU) = \cU \ . $$
\end{demo}

So the previous theorem guarantees that the mapping which sends a non-degenerate
$^*$-homomorphism $\th : \ah_u \rightarrow M(C)$ into the corepresentation $(\io \ot
\th)(\uU)$, is a bijection between the set of non-degenerate $^*$-homomorphisms from $\ah_u$
into $M(C)$ and the set of unitary corepresentations of $(A_u,\de_u)$ on $C$.

\bigskip

In \cite{VD1}, A. Van Daele proved that we can identify $(A,\de)$ and $(\ahh,\dhh)$. In
\cite{JK3}, we used this identification to prove that $\flip(\aU)$ is the universal
algebraic corepresentation of $(\ah,\deh)$.

This will imply that $\flip(\uU)$ can be used in the same way to get a bijection from the
set of non-degenerate $^*$-homomorphisms on $A_u$ and the set of unitary corepresentations
of $(\ah_u,\deh_u)$.

\section{Appendix: some information about weights} \label{app}

In this section, we will collect some necessary information and conventions about weights
and slicing with weights.

\subsection{Weights}

In this first section, we give some information about weights. The  standard reference for
lower semi-continuous weights is \cite{Comb}. We start of with some standard notions
concerning lower semi-continuous weights.

\medskip

Consider a $C^*$-algebra $A$  and a densely defined lower semi- continuous weight $\varphi$
on $A$. We will use the following notations:
\begin{itemize}
\item ${\cal M}^+_\varphi = \{\, a \in A^+ \mid \varphi(a) < \infty  \,\} $
\item ${\cal N}_\varphi = \{\, a \in A \mid \varphi(a^*a) < \infty \,\} $
\item ${\cal M}_\varphi = \text{span\ } {\cal M}^+_\varphi
= {\cal N}_\varphi^* {\cal N}_\varphi$ .
\end{itemize}

\medskip

A GNS-construction of $\varphi$ is by definition a triple
$(H_\varphi,\pi_\varphi,\Lambda_\varphi)$ such that
\begin{itemize}
\item $H_\varphi$ is a Hilbert space
\item $\Lambda_\varphi$ is a linear map from ${\cal N}_\varphi$ into
$H_\varphi$ such that
\begin{enumerate}
\item  $\Lambda_\varphi({\cal N}_\varphi)$ is dense in $H_\varphi$
\item   We have for  every $a,b \in {\cal N}_\varphi$, that
$\langle \Lambda_\varphi(a),\Lambda_\varphi(b) \rangle
= \varphi(b^*a) $
\end{enumerate}
Because $\varphi$ is lower semi-continuous, $\Lambda_\varphi$ is closed.
\item $\pi_\varphi$ is a non-degenerate representation of $A$ on
$H_\varphi$ such that $\pi_\varphi(a)\,\Lambda_\varphi(b) = \Lambda_\varphi(ab)$ for
every $a \in M(A)$ and $b \in {\cal N}_\varphi$.  (The non-degeneracy of
$\pi_\varphi$ is a consequence of the lower semi-continuity of $\varphi$.)
\end{itemize}

\bigskip

The following concepts play a central role in the theory of lower
semi-continuous weights.

\begin{definition}
We define the sets $${\cal F}_\varphi = \{\, \omega \in A^*_+ \mid
\omega \leq \varphi \,\}$$ and
$${\cal G}_\varphi = \{\, \alpha\,\omega \mid
\omega \in {\cal F}_\varphi \, , \, \alpha \in  \,\, ]0,1[ \,\,\} \ . $$
\end{definition}

The advantage of ${\cal G}_\vfi$ over ${\cal F}_\vfi$ is the fact that ${\cal G}_\varphi$
is a directed subset of ${\cal F}_\varphi$ (under the normal order on $A^*_+$). A proof of
this fact can be found in \cite{Q-V} or \cite{Verd}. There is also proof for this fact for
a slightly more general case in proposition 2.13 of \cite{JK2}.

\medskip

The most important result concerning lower semi-continuous weights is the following one
(proved in \cite{Comb}).

\begin{theorem}
We have for every $x \in A^+$ that
$$ \varphi(x) = \sup_{\omega \in {\cal F}_\varphi} \omega(x) \ . $$
\end{theorem}

This can also be stated in terms of the directed set $\cG_\vfi$ as follows:

\begin{proposition}
Consider an element $x \in A^+$. Then
\begin{enumerate}
\item The element $x$ belongs to $\Mfi^+$ $\Leftrightarrow$  The net
$\bigl(\om(x)\bigr)_{\om \in \cG_\vfi}$ is convergent in $\R^+$.
\item If $x$ belongs to $\Mfi^+$, then the net
$\bigl(\om(x)\bigr)_{\om \in \cG_\vfi}$converges to $\vfi(x)$.
\end{enumerate}
\end{proposition}

Then  we get immediately that the net $\bigl(\om(x)\bigr)_{\om \in \cG_\vfi}$ converges to
$\vfi(x)$ for every $x \in \Mfi$.

\bigskip\medskip

Any lower semi-continuous weight $\vfi$ has a natural extension to a weight on $M(A)$.
Remember that every $\om \in A^*$ has a unique extension $\overline{\om}$ to $M(A)$ which
is strictly continuous and we put $\om(x) = \overline{\om}(x)$ for every $x \in M(A)$.

\begin{definition}
We define the weight $\overline{\vfi}$ on $M(A)$ such that $$\overline{\vfi}(x)=
 \sup_{\omega \in {\cal F}_\varphi} \omega(x)  $$
for every $x \in M(A)^+$.
\end{definition}

We define ${\overline{\cal M}}_\vfi = {\cal M}_{\overline{\vfi}}$ and
${\overline{\cal N}}_\vfi = {\cal N}_{\overline{\vfi}}$. For any
$x \in {\overline{\cal M}}_\vfi$, we put $\vfi(x) = \overline{\vfi}(x)$.

\medskip

\begin{proposition}
Consider an element $x \in M(A)^+$. Then
\begin{enumerate}
\item The element $x$ belongs to $\overline{\cM}_\vfi^+$ $\Leftrightarrow$  The net
$\bigl(\om(x)\bigr)_{\om \in \cG_\vfi}$ is convergent in $\R^+$.
\item If $x$ belongs to $\overline{\cM}_\vfi^+$, then the net
$\bigl(\om(x)\bigr)_{\om \in \cG_\vfi}$converges to $\vfi(x)$.
\end{enumerate}
\end{proposition}

As a consequence, we have for every $x \in \overline{\cM}_\vfi$ that the the net
$\bigl(\om(x)\bigr)_{\om \in \cG_\vfi}$converges to $\vfi(x)$.

\bigskip

Consider a GNS-construction $(H_\vfi,\lafi,\pifi)$ of a densely defined lower
semi-continuous weight $\vfi$. Then there is a natural way to get a GNS-construction of
$\overline{\vfi}$ :

\begin{proposition}
The mapping $\lafi : A \rightarrow H$ is closable for the strict topology on $M(A)$ and
the norm topology on $H$. We denote this closure by $\overline{\la}_\vfi$. Then
$(H_\vfi,\overline{\la}_\vfi,\overline{\pi}_\vfi)$ is a GNS-construction for
$\overline{\vfi}$.
\end{proposition}

In particular, we have that $D(\overline{\la}_\vfi) = \overline{\cN}_\vfi$ and we put
$\lafi(x) = \overline{\la}_\vfi(x)$ for every $x \in \overline{\cN}_\vfi$.

\bigskip\bigskip

Now, we will introduce the class of KMS-weights. All weights used in this paper, will
belong to this class. For some more details, we refer to \cite{JK1}.

\begin{definition}
Consider a \cst-algebra $A$ and a weight $\vfi$ on $A$, we say that $\vfi$ is a KMS-weight
on A if and only if $\vfi$ is a  densely defined lower semi-continuous weight on A such
that there exists  a norm-continuous one-parameter group $\si$ on $A$ satisfying the
following properties:
\begin{enumerate}
\item $\vfi$ is invariant under $\si$: $\vfi \si_t = \vfi$ for every $t \in \R$.
\item We have for every $a \in {\cal D}(\si_\frac{i}{2})$ that
      $\vfi(a^* a) = \vfi(\si_\frac{i}{2}(a) \si_\frac{i}{2}(a)^*)$.
\end{enumerate}
The modular group $\si$ is called a modular group for $\vfi$.
\end{definition}

If  the weight $\vfi$ is faithful, then the one-parameter group $\si$ is uniquely
determined and is called the modular group of $\si$.  This is not the usual definition of
a KMS-weight (see \cite{Comb1}), but we prove in \cite{JK1} that this definition is
equivalent with the usual one.

\medskip

In the next proposition,  we formulate some basic properties of KMS-weights. Therefore,
let $(H_\vfi,\lafi,\pifi)$ be a GNS-construction for $\vfi$.

\begin{proposition}
Consider a \cst-algebra $A$, and a KMS-weight $\vfi$ on $A$ with modular group $\si$.
Then:
\begin{enumerate}
\item There exists a unique anti-unitary operator $J$ on $H_\vfi$ such
      that $J \lafi(x) = \lafi(\si_\frac{i}{2}(x)^*)$ for every $x \in
      \Nfi \cap {\cal D}(\si_\frac{i}{2})$.
\item Let $a \in {\cal D}(\si_\frac{i}{2})$ and $x \in \Nfi$. Then $x a$
belongs to $\Nfi$ and $\lafi(x a) = J \pi_\vfi(\si_\frac{i}{2}(a))^* J \, \lafi(x)$.
\item Let $a \in {\cal D}(\si_{-i})$ and $x \in \Mfi$. Then $a x$ and $x \si_{-i}(a)$
belong to $\Mfi$ and $\vfi(a x) = \vfi(x \si_{-i}(a))$.
\end{enumerate}
\end{proposition}

\subsection{The tensor products of weights} \label{app2}

It the next part we will quickly say something about the tensor product of two
KMS-weights.

Therefore, consider two \cst-algebras $A$ and $B$. Let $\vfi$ be a KMS-weight on $A$ with
modular group $\si$ and $\psi$ a KMS-weight on $B$ with modular group $\tau$.

\medskip

\begin{definition}
We define the tensor product weight $\vfi \ot \psi$ on $A \ot B$ in such a way that
$$(\vfi \ot \psi)(x) = \sup \, \{ \, (\om \ot \th)(x) \mid
\om \in \cF_\vfi , \eta \in \cF_\psi \, \} $$
for every $x \in (A \ot B)^+$. Then $\vfi \ot \psi$ is a densely defined lower
semi-continuous weight on $A \ot B$ such that $\cM_\vfi \od \cM_\psi \subseteq \cM_{\vfi
\ot \psi}$ and $(\vfi \od \psi)(a \ot b) =
\vfi(a) \, \psi(b)$ for $a \in \cM_\vfi$ and $b \in \cM_\psi$.
\end{definition}

Of course, this definition is also possible for lower semi-continuous weights and is done
in \cite{Q-V} ans \cite{Verd}.

\bigskip

Consider now a GNS-construction $(H_\vfi,\la_\vfi,\pi_\vfi)$ for $\vfi$ and a
GNS-construction $(H_\psi,\la_\psi,\pi_\psi)$ for $\psi$.

\begin{proposition}
The mapping $\la_\vfi \od \la_\psi : \cN_\vfi \od
\cN_\psi \rightarrow H_\vfi \ot H_\psi$ is closable and we denote the closure of it by
$\la_\vfi \ot \la_\psi$. Then $(H_\vfi \ot H_\psi,\la_\vfi \ot \la_\psi,\pi_\vfi \ot
\pi_\psi)$ is a GNS-construction for $\vfi \ot \psi$.
\end{proposition}

The closability of $\la_\vfi \od \la_\psi$ is also true for lower semi-continuous weights.
However the KMS-condition is used to prove the second part of this proposition (For this,
it would also be sufficient that $\vfi$ and $\psi$ satisfy a weaker condition, the
so-called regularity condition, see \cite{Q-V} and \cite{Verd}).

\medskip

Using the results of \cite{JK1}, it is not so difficult  to prove that $\vfi \ot \psi$ is
a KMS-weight on $A \ot B$ with a modular group $\si \ot \tau$ defined in such a way that
$(\si \ot \tau)_t = \si_t \ot \tau_t$ for every $t \in \R$.

\medskip

We would also like to mention that all of the results (except the lasr one) are also true
in the case of so-called regular \cst-valued weights (see section 8 of \cite{JK2}).

\subsection{Slicing with weights}   \label{app3}

In a last  part we say something about slicing with weights. For a more detailed
exposition, we refer to \cite{JK2}, \cite{Q-V} and \cite{Verd}. The main reference is
\cite{JK2}. A lot of the results cannot be found literally in \cite{JK2}, mainly because
we are in the very special case where $\io$ is a $^*$-homomorphism.

\medskip

However, we believe this part of the appendix is readable without knowing the precise
results of \cite{JK2} and that it gives a fairly good idea of the problems in \cite{JK1}.

The reader is encouraged to fill in the details of this appendix, knowing that the most
difficult problems are resolved in \cite{JK1} (mainly in section 7 and section 8).

\bigskip

For the rest of this subsection, we fix two \cst-algebras  A  and B and a KMS-weight
$\vfi$ on $B$. Let $(H_\vfi,\lafi,\pifi)$ be a GNS-construction for $\vfi$.

\medskip

\begin{definition}
We define the map $\io \ot \vfi$ from within $(A \ot B)^+$ into $A^+$ as follows:
\begin{itemize}
\item We define the set  ${\cal M}_{\io \ot \vfi}^+ = \{\, x \in (A \ot B)^+ \mid
      \text{\ \ the net \ \ } \bigl( \,(\io \ot \om)(x)\,\bigr)_{\om \in {\cal G}_\vfi}
      \text{\ \ is norm convergent } $ \newline
      in $ A \, \}$
\item The mapping $\io \ot \vfi$ will have as domain the set ${\cal M}_{\io \ot \vfi}^+$
      and for any $x \in {\cal M}_{\io \ot \vfi}^+$, we have by definition that the net
      $\bigl( \,(\io \ot \om)(x)\,\bigr)_{\om \in {\cal G}_\vfi}$ converges
      to $(\io \ot \vfi)(x)$.
\end{itemize}
\end{definition}

It is clear that ${\cal M}_{\io \ot \vfi}^+$ is a dense hereditary cone in $(A \ot B)^+$.
Furthermore we define the following sets:
\begin{enumerate}
\item We define the ideal ${\cal N}_{\io \ot \vfi} = \{\, x  \in A \ot B \mid x^* x
\text{\ \ belongs       to \ \ } {\cal M}_{\io \ot \vfi}^+ \,\}$.
\item We define the $^*$-algebra ${\cal M}_{\io \ot \vfi} = \text{span \ }
 {\cal M}_{\io \ot \vfi}^+ = {\cal N}_{\io \ot \vfi}^* {\cal N}_{\io \ot \vfi}$.
\end{enumerate}
Of course, there exist a unique linear map $\psi$ from ${\cal M}_{\io \ot \vfi}$ to $A$
which extend $\io \ot \vfi$.

For any $x \in {\cal M}_{\io \ot \vfi}$, we put $(\io \ot
\vfi)(x) = \psi(x)$.

\medskip

It is then clear that $\bigl( \,(\io \ot \om)(x)\,\bigr)_{\om \in {\cal G}_\vfi}$
converges to $(\io \ot \vfi)(x)$ for every $x \in {\cal M}_{\io \ot \vfi}$.

\bigskip

We also have for every $x \in {\cal M}_{\io \ot \vfi}$ and $\th \in A^*$ that $(\th \ot
\io)(x)$
belongs to $\Mfi$ and $$\vfi\bigl((\th \ot \io)(x)\bigr) = \th\bigl((\io \ot
\vfi)(x)\bigr) \ .$$

\bigskip

We are now going to describe some sort of GNS-construction for $\io \ot \vfi$. It is
possible to prove that the mapping $\io \od \lafi$ from $A \od \Nfi$ into $A \od H$ is
closable (as a mapping from the \cst-algebra  $A \ot B$ into the Hilbert-C$^*$-module $A
\ot H$). We define $\io \ot \lafi$ to be the closure of this mapping $\io \od \lafi$.

\medskip

It is also possible to prove that $D(\io \ot \lafi)$ is a subset of ${\cal N}_{\io \ot \vfi}$
and that
$$(\io \ot \vfi)(y^* x) = \langle (\io \ot \lafi)(x) , (\io \ot \lafi)(y) \rangle $$
for every $x,y \in D(\io \ot \lafi)$.

\medskip

All that is said about $\io \ot \vfi$ until now goes through for lower
semi-continuous weights (and can be found in \cite{Q-V} and
\cite{Verd}). Because we assumed that $\vfi$ is also KMS, we have also
the rather non-trivial result that $D(\io \ot \lafi) = {\cal N}_{\io
\ot \vfi}$. This last result is also true if $\vfi$ would obey a
weaker condition called regularity (see \cite{JK2}).

\bigskip

We have also that $(\io \ot \pifi)(x) (\io \ot \lafi)(a) = (\io \ot \lafi)(x a)$ for every
$x \in M(A \ot B)$ and $a \in \cN_{\io \ot \vfi}$.

\bigskip

It is not difficult to check the following results.

\begin{result} \label{resap1}
Consider an element $a \in M(A)$, then we have the following properties.
\begin{enumerate}
\item We have for every $x \in \cN_{\io \ot \vfi}$ that $x (a \ot 1)$ belongs to
$\cN_{\io \ot \vfi}$ and $(\io \ot \lafi)(x) \, a = (\io \ot \lafi)(x(a \ot 1))$.
\item We have for every $x \in \cM_{io \ot \vfi}$ that $x (a \ot 1)$ and $(a \ot 1) x$
belong to $\cM_{\io \ot \vfi}$ and that
$$ (\io \ot \vfi)(x(a \ot 1)) = (\io \ot \vfi)(x) \, a \hspace{1.5cm}
\text{ and } \hspace{1.5cm}
(\io \ot \vfi)((a \ot 1)x) = a \, (\io \ot \vfi)(x) \ . $$
\end{enumerate}
\end{result}

\bigskip\medskip

For the rest of this section, we have to rely on \cite{JK2}.

\medskip

We would like to have an extension of $\io \ot \vfi$ to $M(A \ot B)$. This is done in the
following way:

\begin{definition}
We define the map $\overline{\io \ot \vfi}$ from within $M(A \ot B)^+$ into $M(A)^+$ as
follows:
\begin{itemize}
\item We define the set  ${\overline{\cal M}}_{\io \ot \vfi}^+ = \{\, x \in M(A \ot B)^+
\mid \text{\ \ the net \ \ } \bigl( \,(\io \ot \om)(x)\,\bigr)_{\om \in {\cal G}_\vfi}
\text{\ \ is strictly convergent}$ \newline  in $M(A) \, \}$.
\item The mapping $\overline{\io \ot \vfi}$ will have as domain the set
${\overline{\cal M}}_{\io \ot \vfi}^+$
and for any $x \in {\overline{\cal M}}_{\io \ot \vfi}^+$, we have  by definition that the
net
$\bigl( \,(\io \ot \om)(x)\,\bigr)_{\om \in {\cal G}_\vfi}$ converges strictly
to $\bigl(\,\overline{\io \ot \vfi}\,\bigr)(x)$.
\end{itemize}
\end{definition}

This definition is in fact not entirely correct because it depends on $\vfi$ and not on
$\io \ot \vfi$. It is possible to give a definition in terms of the mapping
$\io \ot \vfi$ and our definition would then be a proposition.

\medskip

We should also mention that the symbol $\io \ot \vfi$ in \cite{JK2} means something
different than is in this paper. In fact, $\io \ot \vfi$ in \cite{JK2} denotes the
restriction of our map $\overline{\io \ot \vfi}$ to $\overline{\cM}_{\io \ot \vfi}^+ \cap
A^+$. This difference is however not fundamental.

\medskip

The next proposition reveals a nice feature about $\overline{\io \ot \vfi}$.

\begin{proposition}
Consider $x \in M(A \ot B)^+$. Then $x$ belongs to ${\overline{\cal M}}_{\io \ot \vfi}^+ $
if and only if the net \newline $\bigl( \,a^* (\io \ot \om)(x) a \,\bigr)_{\om \in {\cal
G}_\vfi}$ is norm convergent for every $a \in A$.
\end{proposition}

\medskip

It is clear that ${\overline{\cal M}}_{\io \ot \vfi}^+$ is a hereditary cone in $M(A \ot
B)^+$. Furthermore we define the following sets:
\begin{enumerate}
\item We define the left ideal ${\overline{\cal N}}_{\io \ot \vfi} = \{\, x \in M(A \ot B)
\mid x^* x \text{\ \ belongs to \ \ }
{\overline{\cal M}}_{\io \ot \vfi}^+ \,\}$.
\item We define the $^*$-algebra ${\overline{\cal M}}_{\io \ot \vfi} = \text{span \ }
{\overline{\cal M}}_{\io \ot \vfi}^+ = {\overline{\cal N}}_{\io \ot \vfi}^*
{\overline{\cal N}}_{\io \ot \vfi}$.
\end{enumerate}

Of course, there exist a unique linear map $\overline{\psi}$ from ${\overline{\cal
M}}_{\io \ot \vfi}$ to $M(A)$ which extend $\overline{\io \ot \vfi}$.

For any $x \in {\overline{\cal M}}_{\io \ot \vfi}$, we put $(\io \ot \vfi)(x) =
\overline{\psi}(x)$.

\medskip

It is then clear that $\bigl( \,(\io \ot \om)(x)\,\bigr)_{\om \in {\cal G}_\vfi}$
converges strictly to $(\io \ot \vfi)(x)$ for every $x \in {\overline{\cal M}}_{\io \ot
\vfi}$.

\bigskip

We also have for any $x \in {\overline{\cal M}}_{\io \ot \vfi}$ and $\th \in A^*$ that
$(\th \ot \io)(x)$ belongs to ${\overline{\cal M}}_\vfi$ and $$\vfi\bigl((\th \ot
\io)(x)\bigr) = \th\bigl((\io \ot \vfi)(x)\bigr) \ .$$

\bigskip

We want also to get something like a GNS-construction for $\overline{\io \ot \vfi}$.

\medskip

Notice that $\cK(A , A \ot H)$ can be turned into a Hilbert-\cst-module by defining the
scalar product on $\cK(A,A \ot H)$ in such a way that $\langle x , y \rangle = y^* x$ for
every $x,y \in \cK(A,A \ot H)$.

We can also define a mapping $\Upsilon$ from $A \ot H$ into $\cK(A,A \ot H)$ such that
$\Upsilon(v)\,a = v\,a$ for every $a \in A$ and $v \in A \ot H$. Then this mapping is a
unitary transformation form $A \ot H$ to $\cK(A,A \ot H)$.

We will use this mapping $\Upsilon$ to identify $\cK(A,A \ot H)$ with $A \ot H$.

\medskip

So we can consider $\io \ot \lafi$ as a closed linear mapping from $\cN_{\io \ot \vfi}$
into $\cK(A,A \ot H)$.

\vspace{1.5mm}

One can prove that $\io \ot \lafi$ is closable for the strict topology on $M(A \ot B)$ and
the strong-$^*$-topology on $\cL(A,A \ot H)$ and  we denote this closure by $\overline{\io
\ot \lafi}$.

So  $\overline{\io \ot \lafi}$ is a strictly-strongly$^*$ closed linear map from within
$M(A \ot B)$ into $\cL(A,A \ot H)$.

\medskip

We can prove (but this is not straightforward)  that $\overline{\cN}_{\io \ot
\vfi}= D\bigl(\,\overline{\io \ot \lafi}\,\bigr)$.

\medskip

For every $x \in \overline{\cN}_{\io \ot\vfi}$, we define $(\io \ot \lafi)(x) =
\bigl(\,\overline{\io \ot \lafi}\,\bigr)(x)$, so $(\io \ot \lafi)(x) \in \cL(A,A \ot H)$.

\bigskip

We have moreover for every $x,y \in \overline{\cN}_{\io \ot \vfi}$ that
$$ (\io \ot \vfi)(y^* x) = (\io \ot \lafi)(y)^* (\io \ot \lafi)(x)$$

\bigskip

We have also that $(\io \ot \pifi)(x) (\io \ot \lafi)(a) = (\io \ot \lafi)(x a)$ for every
$x \in M(A \ot B)$ and $a \in \overline{\cN}_{\io \ot \vfi}$.

\bigskip

It is not difficult to check the following results.

\begin{result}
Consider an element $a \in M(A)$, then we have the following properties.
\begin{enumerate}
\item We have for every $x \in \overline{\cN}_{\io \ot \vfi}$ that $x (a \ot 1)$ belongs to
$\overline{\cN}_{\io \ot \vfi}$ and $(\io \ot \lafi)(x) \, a = (\io \ot \lafi)(x (a
\ot1))$.
\item We have for every $x \in \overline{\cM}_{io \ot \vfi}$ that  $x (a \ot 1)$
and $(a \ot 1) x$ belong to $\overline{\cM}_{\io \ot \vfi}$ and that
$$ (\io \ot \vfi)(x(a \ot 1)) = (\io \ot \vfi)(x) \, a \hspace{1.5cm} \text{ and }
\hspace{1.5cm} (\io \ot \vfi)((a \ot 1) x) = a \, (\io \ot \vfi)(x) \ . $$
\end{enumerate}
\end{result}

\bigskip

We would like to end this part with some remarks
\begin{itemize}
\item  Consider $x \in \overline{\cN}_{\io \ot \vfi}$, then $(\io \ot \vfi)(x^* x)$
belongs to $A$ if and only $(\io \ot \lafi)(x)$ belongs to $\cK(A,A \ot H)$.
\item Consider $x \in (A \ot B)^+$, then $x$ belongs to $\cM_{\io \ot \vfi}^+$
if and only if $x$ belongs to $\overline{\cM}_{\io \ot \vfi}^+$ and $(\io \ot \vfi)(x)$
belongs to $A^+$.
\end{itemize}

\bigskip\bigskip

We are also interested in extending $\io \ot \vfi$ even further and let it take values in
the set of elements affiliated with $A$.

\bigskip

Let $x \in M(A \ot B)$ and define the set
$$D_x = \{ \, a \in A \mid x (a \ot 1) \in \overline{\cN}_{\io \ot \vfi} \text{ and }
(\io \ot \vfi)((a^* \ot 1) x^* x (a \ot 1)) \in A \, \} $$

Then we have for every $a \in A$ that $a$ belongs to $D_x$ if and only if the net
$\bigl(\, a^* \, (\io \ot \om)(x^* x) \, a \, \bigr)_{\om \in \cG_\vfi}$ is convergent in
$A$.

\medskip

\begin{definition}  \label{defap1}
We define the set $$\tilde{\cN}_{\io \ot \vfi} = \{ \, x \in M(A \ot B) \mid D_x \text{ is
dense in } A \, \}$$
\end{definition}

\medskip

\begin{definition}  \label{defap2}
Consider $x \in \tilde{\cN}_{\io \ot \vfi}$. Then we define the mapping $(\io \ot
\lafi)(x)$ from $D_x$  into $A \ot H$ as follows : \ \
Let $a \in D_x$. Then we know that $x (a \ot 1)$ belongs to $\overline{\cN}_{\io \ot
\vfi}$ and $(\io \ot \lafi)(x (a \ot 1))$ belongs to $\cK(A,A \ot H)$.
So there exists a unique element $v \in A \ot H$ such that $(\io \ot \lafi)
(x (a \ot 1))\, b = v \, b $ for every $b \in A$
and we define $(\io \ot \lafi)(x)(a) = v$.
\end{definition}

It is then possible to prove that $(\io \ot \lafi)(x)$ is a closed densely defined linear
operator from within $A$ into $A \ot H$ such that $(\io \ot \lafi)(x)^*$ is densely
defined.

\medskip

\begin{definition} \label{defap3}
We define the set $$\hat{\cN}_{\io \ot \vfi} = \{ \, x \in \tilde{\cN}_{\io \ot \vfi} \mid
(\io \ot \lafi)(x) \text{ is regular } \, \}  \ . $$
\end{definition}

\begin{definition}  \label{defap4}
We define the set $$\hat{\cM}_{\io \ot
\vfi}^+ = \{ \, x \in M(A \ot B)^+ \mid x^\frac{1}{2} \in \hat{\cN}_{\io \ot \vfi} \ .$$

For every $x \in \hat{\cM}_{\io \ot \vfi}^+$, we put $(\io \ot \vfi)(x) = (\io \ot
\lafi)(x^\frac{1}{2})^* (\io \ot \lafi)(x^\frac{1}{2})$, so $(\io \ot \vfi)(x)$ is a
positive element affiliated with $A$.
\end{definition}

Consider $y \in M(A \ot B)$. Then one can prove that $y \in \hat{\cN}_{\io \ot \vfi}$ if
and only if $y^* y \in \hat{\cM}_{\io \ot \vfi}^+$ and we have in this case that $(\io \ot
\vfi)(y^* y) = (\io \ot \lafi)(y)^* (\io \ot \lafi)(y)$

\bigskip\bigskip

It goes without saying that everything goes also through for the slicing $\vfi \ot \io$.


\begin{thebibliography}{VD}

\bibitem{Abe} {\sc E. Abe},
Hopf Algebras. {\it Cambridge University Press} (1977).


\bibitem{B-S}  {\sc S. Baaj  \&  G. Skandalis},
Unitaires multiplicatifs et dualit\'e pour les produits crois\'es de
\cst-alg\`ebres. {\it Ann. scient. \'{E}c. Norm. Sup., $4{}^e$ s\'{e}rie, t. 26}
(1993), 425--488.

\bibitem{Comb} {\sc F. Combes},
Poids sur une \cst-alg\`ebre.
{\it J. Math. pures et appl.} {\bf 47} (1968), 57--100.


\bibitem{Comb1} {\sc F. Combes},
Poids associ\'e \`a une alg\`ebre hilbertienne \`a gauche.
{\it Compos. Math.} {\bf 23} (1971), 49--77.




\bibitem{ER} {\sc  E.G. Effros \&  Z.-J. Ruan},
Discrete Quantum Groups I. The Haar Measure. {\it Int. J. of Math.} (1994).

\bibitem{E}  {\sc  M. Enock \&  J.-M. Schwartz},
Kac Algebras and Duality of Locally Compact Groups.
{\it Springer-Verlag, Berlin}  (1992).

\bibitem{Drab}{\sc B. Drabant \& A. Van Daele},
Pairing and Quantum Double of Multiplier Hopf
Algebras. {\it Preprint K.U.Leuven}  (1996)


\bibitem{GL} {\sc  E.C. Gootman \& A.J. Lazar},
Quantum Groups and Duality.
{\it Reviews in Math. Physics} {\bf 5} No. 2 (1993),
417--451

\bibitem{Kus} {\sc J. Kustermans \& A. Van Daele},
\cst-algebraic quantum groups arising from algebraic quantum groups. (1996)
\ To appear in {\it International Journal of Mathematics}.



\bibitem{JK1} {\sc J. Kustermans}, A construction procedure for
KMS-weights on \cst-algebras. In preparation.

\bibitem{JK3} {\sc J. Kustermans}, Examining the dual
of an algebraic quantum group. {\it Preprint Odense Universitet} (1997).

\bibitem{JK4} {\sc J. Kustermans}, A natural extension of
a left invariant lower semi-continuous weight. {\it Preprint Odense Universitet} (1997).

\bibitem{JK2} {\sc J. Kustermans}, Regular \cst-valued weights on \cst-algebras.
{\it Preprint K.U. Leuven} (1997).






\bibitem{Lan} {\sc C. Lance},
Hilbert $C^*$-modules, a toolkit for operator algebraists. Leeds. (1993).

\bibitem{MasNak} {\sc T. Masuda \& Y. Nakagami} A von Neumann Algebra
Framework for the Duality
of Quantum Groups. {\it Publications of the RIMS Kyoto University}
{\bf 30} (1994) , 799--850

\bibitem{Pe-Tak} {\sc G.K. Pedersen \& M. Takesaki},
The Radon-Nikodym theorem for von Neumann algebras.
{\it Acta Math.} {\bf 130} (1973), 53--87.


\bibitem{PW}  {\sc  P. Podle\'s \&  S.L. Woronowicz}, Quantum Deformation
of the Lorentz Group.
{\it Commun. Math. Phys.} {\bf 130} (1990), 381--431.


\bibitem{Q-V} {\sc J. Quaegebeur \& J. Verding},
A construction for weights on $C^*$-algebras. Dual weights for $C^*$-crossed
products. {\it Preprint K.U. Leuven} (1994).


\bibitem{Ver} {\sc J. Quaegebeur \& J. Verding},
Left invariant weights and the left regular corepresentation
for locally compact quantum semi-groups. {\it Preprint K.U. Leuven} (1994).


\bibitem{Stra} {\sc S. Stratila \& L. Zsido},
Lectures on von Neumann algebras. {\it Abacus Press, Tunbridge Wells, England} (1979).

\bibitem{Tay}  {\sc D.C. Taylor}, The Strict Topology for Double Centralizer Algebras.
{\it Trans. Am. Math. Soc} {\bf 150} (1970), 633 -- 643

\bibitem{Tak}  {\sc  M. Takesaki},
Theory of Operator Algebras I.
{\it Springer-Verlag, New York} (1979).



\bibitem{VD1}   {\sc A. Van Daele}, An Algebraic Framework for Group Duality.
(1996) \ To appear in {\it Advances of Mathematics}.

\bibitem{VD2}   {\sc  A. Van Daele}, Dual Pairs of Hopf ${}^*$-algebras.
{\it Bull. London Math. Soc.} {\bf 25} (1993), 209--230.


\bibitem{VD4}   {\sc  A. Van Daele}, Discrete Quantum Groups.
{\it Journal of Algebra} {\bf 180} (1996), 431--444.

\bibitem{VD5}   {\sc   A. Van Daele}, The Haar Measure on a Compact
Quantum Group. {\it Proc. Amer. Math. Soc.} {\bf 123}(1995), 3125-3128

\bibitem{VD6}   {\sc   A. Van Daele}, Multiplier Hopf Algebras.
{\it Trans. Am. Math. Soc.} {\bf 342} (1994), 917--932.

\bibitem{Verd} {\sc J. Verding}, Weights on \cst-algebras.
{\it Phd-thesis.  K.U. Leuven} (1995)

\bibitem{Wor1}   {\sc  S.L. Woronowicz},  Compact matrix pseudogroups.
{\it Commun. Math. Phys.} {\bf 111}  (1987),  613--665.

\bibitem{Wor2}   {\sc  S.L. Woronowicz},  Compact quantum groups.
{\it Preprint Warszawa} (1993).

\bibitem{Wor3}  {\sc  S.L. Woronowicz}, From multiplicative unitaries to quantum
groups. {\it Preprint Warszawa} (1995).

\bibitem{Wor4}  {\sc  S.L. Woronowicz}, Pseudospaces, pseudogroups
and Pontriagin duality. {\it Proceedings of the International Conference on Mathematical
Physics, Lausanne} (1979),  407--412.


\bibitem{Wor6}  {\sc  S.L. Woronowicz}, Unbounded elements affiliated with $C^*$-algebras
and non-compact quantum groups.
{\it Commun. Math. Phys.} {\bf  136} (1991),  399--432.


\end{thebibliography}
\end{document}